\documentclass[11pt]{article}
\usepackage{preprint}
\usepackage{hyperref}
\usepackage{url}
\usepackage{booktabs}
\usepackage{graphicx}
\usepackage{fontawesome5}
\usepackage{tabularx}
\usepackage{makecell}
\usepackage{algorithm}
\usepackage[noend]{algorithmic}
\usepackage{arydshln}
\usepackage{amssymb}
\usepackage{nicefrac}   
\usepackage{enumitem}
\usepackage[textsize=footnotesize]{todonotes}
\usepackage{agafonov}

\newtheorem{theorem}{Theorem}[section]
\newtheorem{lemma}[theorem]{Lemma}
\newtheorem{proposition}[theorem]{Proposition}
\newtheorem{corollary}[theorem]{Corollary}
\theoremstyle{definition}
\newtheorem{example}[theorem]{Example}
\newtheorem{definition}[theorem]{Definition}
\theoremstyle{plain}

\title{Optimal Parameter-Free Second-order Acceleration with Inexact Hessians under Generalized Smoothness}

\author{%
\begin{tabular}{c@{\hspace{3em}}c@{\hspace{3em}}c}
Artem Agafonov$^{1}$ & Aslan Ashabokov$^{2}$ & Alexander D'yakonov \\[6pt]
Martin Tak\'a\v{c}$^{1}$ & Alexander Gasnikov$^{2}$ & Dmitry Kamzolov
\end{tabular}\\[10pt]
{\normalsize $^{1}$MBZUAI\qquad $^{2}$MIRAI}}

\newsavebox{\backtrackingbox}
\newsavebox{\basicnewtonbox}
\newlength{\pairedalgorithmheight}

\begin{document}

\maketitle
\begin{abstract}
  We study second-order convex optimization under generalized Hessian smoothness, where the local Hessian Lipschitz constant may grow linearly with the gradient norm. For inexact Hessians, we develop adaptive algorithms, including an accelerated Monteiro-Svaiter-type method whose complexity matches the lower bound for second-order methods with $\delta$-inexact Hessians in the classical Lipschitz-Hessian setting. Our methods are parameter-free in the sense that they require no problem-specific parameters as input, including the Hessian inexactness level $\delta$ and smoothness constants. Instead, their backtracking adapts a single regularization parameter that jointly captures Hessian inexactness and local Hessian smoothness. The resulting rates split into two regimes. Far from the solution, the gradient-dependent smoothness term yields geometric convergence. Near the solution, the rate is governed by the ordinary Hessian smoothness and the inexactness level $\delta$, recovering the classical Lipschitz-Hessian guarantees. 
\end{abstract}

\addtocontents{toc}{\protect\setcounter{tocdepth}{-1}}
\section{Introduction}
We study second-order methods for the problem
$\min_{x\in \R^d} f(x)$, 
where $f \in \mathcal{C}^2$ is a convex function. We assume that a solution $x^{\ast}\in \R^d$ exists and we denote $f^{\ast}\eqdef f(x^{\ast})$. Given an initial point $x_0$, we define $R \eqdef \|x_0 - x^{\ast}\|$ as a distance to the solution. The classical theory assumes exact Hessians that are Lipschitz
continuous~\citep{nesterov2006cubic,nesterov2018lectures}:
\begin{equation*}
  \norm{\nabla^2 f(y)-\nabla^2 f(x)}\leq L\norm{y-x},
  \qquad x,y\in\R^d.
  \tag{H}\label{as:h}
\end{equation*}
In this setting, optimal methods find an $\e$-solution in
$\cO\bigl((LR^3/\e)^{2/7}\bigr)$
iterations~\citep{kovalev2022first,carmon2022optimal}, matching the lower
bound~\citep{arjevani2019oracle}. We consider two extensions of this
theory.\\[2pt]
\textit{Generalized smoothness.}
In first-order optimization, $(L_0,L_1)$-smoothness lets the local
Lipschitz-gradient constant grow with the gradient
norm~\citep{zhang2019gradient,crawshaw2022robustness}. We study its
second-order analogue~\citep{xie2024trust,jerad2026fast}: for some
$M_0\geq 0$, $M_1\geq0$, and $\rho>0$,
\begin{equation*}
  \norm{\nabla^2 f(y)-\nabla^2 f(x)}
  \leq\ls M_0+M_1\norm{\grad f(x)}\rs\norm{y-x},
  \qquad x,y\in\R^d:\ \norm{y-x}\leq\rho.
  \tag{LH}\label{as:lh}
\end{equation*}
It reduces to~\eqref{as:h} for $M_1=0$ and $\rho=+\infty$, and it covers
problems such as Poisson regression and $\ell_p$ regression with
$p\geq4$, for which~\eqref{as:h} fails
(see examples in Section~\ref{sec:examples}).
\\[2pt]
\textit{Inexact Hessians.}
Exact Hessians are often too expensive. In practice, methods use approximations
obtained by subsampling~\citep{kohler2017sub,xu2020newton}, quasi-Newton
updates~\citep{kamzolov2026accelerated}, or lazy
updates~\citep{doikov2023second}. We assume that the approximation $H(x)$
is symmetric and satisfies
\begin{equation*}
  \norm{H(x)-\nabla^2 f(x)}\leq\delta.
  \tag{IH}\label{as:inexact_hessian}
\end{equation*}

\textbf{Parameter-free methods.\footnote{%
    By parameter-free, we mean requiring no problem constants as input and converging for any initialization.}} 
Both extensions make assumptions more realistic, but they add constants
that a method has to know: \eqref{as:lh} involves $M_0$, $M_1$, and
$\rho$, and~\eqref{as:inexact_hessian} adds $\delta$. These constants are
rarely known. Moreover, under~\eqref{as:lh} they are not unique. Poisson
regression satisfies~\eqref{as:lh} for every $\rho>0$, with $M_0$ and
$M_1$ growing in $\rho$ (Example~\ref{ex:poisson_lh_not_gh}). A method
that takes them as input must fix $\rho$ in advance, whereas the
guarantees of a parameter-free method hold for the best $\rho$. 
\\[2pt]
To the best of our knowledge, even under~\eqref{as:h}, no parameter-free method attains the optimal complexity with inexact Hessians. The method of~\citet{carmon2022optimal} is optimal and adaptive, but its analysis requires exact Hessians. The optimal method of~\citet{chen2026optimal} requires $L$, $\delta$, and $R$ as input. Under generalized smoothness, existing convex guarantees are non-accelerated~\citep{semenov2026gradient,jerad2026fast}. In view of these results, we ask two questions:
\begin{center}
  \textit{Can a parameter-free
  method attain the optimal complexity with inexact Hessians?}\\[0pt]
  \textit{Can such a method be extended to generalized Hessian
  smoothness?}
\end{center}

\textbf{Our main result.} We answer both questions affirmatively with a single method
(Algorithm~\ref{alg:ahpe_damped}). Using the Monteiro--Svaiter
condition~\citep{monteiro2013accelerated} as an acceptance test, we do not estimate the Hessian inexactness and the local Lipschitz constant of the Hessian
separately, but rather adapt them as a single regularization parameter. Under~\eqref{as:h}
and~\eqref{as:inexact_hessian}, the proposed method needs $\cO\ls (LR^3/\e)^{2/7}+\sqrt{\delta R^2/\e}\rs$  Hessian evaluations, up to an additive logarithmic term independent of
$\e$. This matches the lower bound~\citep{agafonov2024advancing}.
Under~\eqref{as:lh}, the same method converges linearly far from the
solution and then at the same optimal rate with $L=M_0$. Table~\ref{tab:convex_comparison} compares our main methods with closest competitors and the lower bound.

\begin{table}
  \centering
  \caption{Comparison of our methods with the fastest known second-order convex methods and the lower bound with inexact Hessians}
  \label{tab:convex_comparison}
  \small
  \setlength{\arrayrulewidth}{0.4pt}
  \renewcommand{\tabularxcolumn}[1]{m{#1}}
  \newcommand{\padcell}[1]{\makecell{\rule{0pt}{13pt}#1\rule[-6pt]{0pt}{0pt}}}
  \begin{tabularx}{\textwidth}{@{}
    >{\centering\arraybackslash}m{4.6cm} %
    >{\centering\arraybackslash}m{1.5cm} %
    >{\centering\arraybackslash}m{1.2cm} %
    >{\centering\arraybackslash}X %
  @{}}
    \toprule
    Method & Assump.
      & \makecell[c]{Param.\\Free}
      & Hessian complexity \\
    \midrule
    \padcell{Lower bound\\
      {~\citep{agafonov2024advancing}}}
      & \makecell[c]{~\eqref{as:h}\\~\eqref{as:inexact_hessian}} & \scalebox{.7}{\faMinus}
      & $\Omega\bigl((\tfrac{LR^3}{\e})^{2/7}
          +\sqrt{\tfrac{\delta R^2}{\e}}\bigr)$ \\
    \midrule

    \padcell{Optimal MS Acceleration\\
      {~\citep{carmon2022optimal}}}
      & ~\eqref{as:h} & \cmark
      & $\cO\bigl((\tfrac{LR^3}{\e})^{2/7}\bigr)$ \\
    \hdashline

    \padcell{Optimal inexact Newton\\
      {~\citep{chen2026optimal}}}
      & \makecell[c]{~\eqref{as:h}\\~\eqref{as:inexact_hessian}} & \xmark
      & $\cO\bigl((\tfrac{LR^3}{\e})^{2/7}
          +\sqrt{\tfrac{\delta R^2}{\e}}\bigr)$ \\
    
    \hdashline
    
    \padcell{Gradient-Regularized Newton\\
      {~\citep{semenov2026gradient}}}
      & \makecell[c]{~\eqref{as:gh}\\~~~~~\eqref{as:inexact_hessian}$^{a}$} & \cmark
      & $\widetilde{\cO}\bigl(
          \sqrt{\tfrac{M_0D^3}{\e}}+\sqrt{M_1}D
          +\tfrac{\delta D^2}{\e}\bigr)^{b}$ \\
    \hdashline
    \padcell{AR2-LS
      {~\citep{jerad2026fast}}}
      &~\eqref{as:lh} & \cmark
      & $\widetilde{\cO}\bigl(\e^{-1/2}\bigr)$ \\
    \hdashline
    \padcell{Adaptive regularized Newton\\
      {(This work, Thm.~\ref{thm:basic_pd_ms_convergence})}}
      & \makecell[c]{~\eqref{as:lh}\\~\eqref{as:inexact_hessian}} & \cmark
      & $\cO\bigl(\sqrt{\tfrac{M_0D^3}{\e}}+\tfrac{\delta D^2}{\e}\bigr)
          +\widetilde{\cO}\bigl(1+\sqrt{M_1}D+\tfrac{D}{\rho}\bigr)$ \\
    \hdashline
    \padcell{Inexact adaptive A-NPE\\
      {(This work, Thm.~\ref{thm:ahpe_convergence}})}
      & \makecell[c]{~\eqref{as:lh}\\~\eqref{as:inexact_hessian}} & \cmark
      & $\cO\bigl((\tfrac{M_0R^3}{\e})^{1/3}
            +\sqrt{\tfrac{\delta R^2}{\e}}\bigr)+\widetilde{\cO}\bigl(1+\sqrt{M_1}R+\tfrac{R}{\rho}\bigr)$\\
    \hdashline
    \padcell{Inexact adaptive \\ damped A-NPE\\
      {(This work, Thm.~\ref{thm:ahpe_damped_convergence})}}
      & \makecell[c]{~\eqref{as:lh}\\~\eqref{as:inexact_hessian}} & \cmark
      & $\cO\bigl((\tfrac{M_0R^3}{\e})^{2/7}
            +\sqrt{\tfrac{\delta R^2}{\e}}\bigr)
        +\widetilde{\cO}\bigl(
            1+M_1^{1/3}R^{2/3}+(\tfrac{R}{\rho})^{2/3}\bigr)$\\
    \bottomrule
  \end{tabularx}
  \begin{minipage}{\textwidth}
    \scriptsize
    $D=\operatorname{diam}\{f\leq f(x_0)\}<\infty$,
    $G=\sup_{f\leq f(x_0)}\norm{\nabla f}$.
    $\widetilde{\cO}$ hides logarithmic factors.
    Parameter-free methods require no problem parameters.
    \textsuperscript{a}~$H_k\succeq0$ is also assumed.
    \textsuperscript{b}~Combines Corollary~2 and Appendix~F.6 of~\citet{semenov2026gradient}.
  \end{minipage}
\end{table}

\subsection{Related Work}

\textbf{Beyond classical smoothness.}
In first-order optimization, empirical observations on neural networks have motivated $(L_0,L_1)$-smoothness assumption~\citep{zhang2019gradient,crawshaw2022robustness}, which is studied for both convex~\citep{koloskova2023revisiting,li2023convex,tyurin2024toward,gorbunov2025methods} and nonconvex~\citep{koloskova2023revisiting,chen2023generalized,vankov2025optimizing} objectives. In particular, convex accelerated under $(L_0,L_1)$-smoothness has recently received considerable attention~\citep{gorbunov2025methods,vankov2025optimizing,tyurin2025near}, including parameter and line-search free methods~\citep{borodich2026nesterov}. Second-order methods under generalized smoothness have been developed in
nonconvex~\citep{xie2024trust,fang2026parameter} and
convex~\citep{semenov2026gradient,jerad2026fast}
settings, some under the global variant of~\eqref{as:lh}
\begin{equation*}
    \|\nabla^2 f(y)-\nabla^2 f(x)\|
    \leq
    \bigl(M_0+M_1\|\nabla f(x)\|\bigr)\|y-x\|,
    \qquad x,y\in\mathbb{R}^d.
    \tag{GH}\label{as:gh}
\end{equation*}

\textbf{Inexact Hessians.}
Under~\eqref{as:h}, many methods with inexact Hessians augment cubic
regularization with a quadratic term that compensates for the Hessian
error~\citep{ghadimi2017second,agafonov2023inexact,agafonov2024advancing}.
These methods take both $L$ and $\delta$ as input.
\citet{kamzolov2026accelerated} adapt the quadratic coefficient but still
require a cubic coefficient $\M\geq2L$. \citet{antonakopoulos2022extra}
propose Extra-Newton, a parameter-free but suboptimal method.
\citet{semenov2026gradient} propose an adaptive but non-accelerated
method for generalized Hessian smoothness with inexact Hessians. Recently,~\citet{chen2026optimal} proposed an optimal method matching the lower bound~\citep{agafonov2024advancing}, but their method is non-adaptive.\\[2pt]
\textbf{Adaptive second-order methods.}
With exact Hessians and under~\eqref{as:h}, adaptive cubic regularization estimates the cubic coefficient by backtracking~\citep{nesterov2006cubic,cartis2011adaptive,grapiglia2020tensor}.
Gradient-regularized Newton methods replace the cubic subproblem with a
linear system~\citep{mishchenko2023regularized,doikov2023gradient}, and
their adaptive variants tune the regularization coefficient in the same
way~\citep{doikov2024super}. \citet{carmon2022optimal} obtain an optimal
adaptive method by searching over the regularization of a Newton step
using Monteiro--Svaiter acceleration~\citep{monteiro2013accelerated}.

\subsection{Contributions}

\begin{itemize}[noitemsep,topsep=0pt,leftmargin=12pt]
  \item \textbf{Far/near rate separation under~\eqref{as:lh}.} We analyze an adaptive regularized Newton method and three accelerated variants. Far from the solution, the methods converge linearly at a rate that depends on $M_1$ and $\rho$. Near the solution, they recover the standard rates under~\eqref{as:h} (see Table~\ref{tab:convex_comparison}).
  \item \textbf{Inexact Hessians and adaptivity.} The three methods listed in Table~\ref{tab:convex_comparison} allow inexact Hessians and remain parameter-free. A single regularization parameter in the Newton step accounts for both the Hessian error and the local smoothness, and methods adapt it by backtracking on the Monteiro-Svaiter condition.
  \item \textbf{Optimal adaptive acceleration.} Algorithm~\ref{alg:ahpe_damped} combines this backtracking with the damped acceleration scheme of~\citet{carmon2022optimal}. Under~\eqref{as:lh} and~\eqref{as:inexact_hessian}, it converges linearly during the first $\widetilde{\cO}\bigl(1+M_1^{1/3}R^{2/3}+(R/\rho)^{2/3}\bigr)$ iterations and then needs $\cO\bigl((M_0R^3/\e)^{2/7}+\sqrt{\delta R^2/\e}\bigr)$ more iterations. Under~\eqref{as:h}, the first phase reduces to an additive logarithmic term, and the bound matches the lower bound~\citep{agafonov2024advancing}. The number of gradient evaluations has the same order up to additive terms at most logarithmic in $1/\e$.
  \item \textbf{Revisited analysis of optimal acceleration.} Our convergence proof is new, and we discuss its differences from the analyses of~\citet{carmon2022optimal} and~\citet{chen2026optimal} in Section~\ref{sec:related_proofs}.
  \item \textbf{Recovering local quadratic rate.} For strongly convex functions, our adaptive regularized Newton converges globally linearly and does not accumulate Hessian errors. If the errors vanish (e.g., with lazy Hessian updates), it converges superlinearly. Under a stronger, relative accuracy condition on the Hessian approximation, it also converges locally quadratically (Section~\ref{sec:basic_pd_ms_strong}).
\end{itemize}

\section{Preliminaries}
\label{sec:preliminaries}

Throughout the paper, $\norm{\cdot}$ denotes the Euclidean norm for vectors
and the induced operator norm for matrices. For both conditions~\eqref{as:gh} and~\eqref{as:lh}, we define $L(x) \eqdef M_0+M_1\norm{\grad f(x)}$.
When using $H(x)$ satisfying~\eqref{as:inexact_hessian}, we assume
w.l.o.g. that it is symmetric, since replacing it by symmetric matrix $(H(x) + H^T(x))/2$ does not increase $\delta$. We set
$\log_+t\eqdef\max\lb0,\log t\rb$ and
$\log_{b,+}t\eqdef\max\lb0,\log_b t\rb$ for $t>0$ and $b>1$. Under~\eqref{as:gh} we take  $\rho=+\infty$, and use $a/(+\infty)=0$ for finite $a\geq0$ and $a/0=+\infty$ for $a>0$. In particular, $1/\sqrt{M_1}=+\infty$ if $M_1=0$.

\subsection{Generalized Hessian Smoothness}
\label{sec:examples}

\textbf{Relations between assumptions.}
Functions satisfying~\eqref{as:h}, such as logistic regression,
also satisfy~\eqref{as:gh} and~\eqref{as:lh}. Examples below show that the inclusions $\eqref{as:h}\subset\eqref{as:gh}\subset\eqref{as:lh}$ are strict.
In particular, setting $x=x^{\ast}$ in~\eqref{as:gh} gives $\|\nabla^2 f(y)\|
  \leq \|\nabla^2 f(x^{\ast})\| + M_0\|y-x^{\ast}\|$.
Thus, the Hessian norm can grow at most linearly with the distance from
$x^{\ast}$, as under~\eqref{as:h}. However,~\eqref{as:gh} is more general than~\eqref{as:h}, since the coefficient controlling the Hessian variation at
$x$ may grow with $\|\nabla f(x)\|$. The local condition~\eqref{as:lh} does not
imply this global linear bound and therefore covers a broader class of
functions. 

\begin{example}
\label{ex:gh_not_h}
The function
$
  f(x)=x^2+\int_0^x(x-t)\sin(t^2)\mathrm{d}t
$
is $1$-strongly convex and satisfies~\eqref{as:gh} with $M_0=M_1=2$, but
it does not satisfy~\eqref{as:h}.
\end{example}

\begin{example}
\label{ex:cosh_lh_not_gh}
The $1$-strongly convex function $f(x)=\cosh x-1$ satisfies~\eqref{as:lh} for every $\rho>0$ with $M_0=\sinh\rho$ and
$M_1=e^\rho$, but it does not satisfy~\eqref{as:gh}.
\end{example}

\begin{example}
\label{ex:exp_lh_not_gh}
The convex function $f(x)=e^x-x-1$ satisfies~\eqref{as:lh} for every
$\rho>0$ with $M_0=M_1=e^\rho$, but it does not satisfy~\eqref{as:gh}.
\end{example}

\begin{example}
\label{ex:quartic_lh_not_gh}
For $\mu\geq0$, consider
$
  f(x)=\frac{\mu}{2}\norm{x}^2+\frac{1}{4}\norm{x}^4
$
for $x \in \R^d$.
The function is convex for $\mu=0$ and $\mu$-strongly convex for $\mu>0$.
For every $\rho>0$, it satisfies~\eqref{as:lh} with $M_0=3\rho+4$ and
$M_1=2$, but it does not satisfy~\eqref{as:gh}.
\end{example}

\begin{example}[Poisson regression]
\label{ex:poisson_lh_not_gh}
Let $a_i\in\R^d$ with $[a_i]_1=1$ and $b_i\in\mathbb{N}$, $i=1,\ldots,n$. The Poisson regression loss $f(x)=\sum_{i=1}^n\ls e^{\la a_i,x\ra}-b_i\la a_i,x\ra\rs$ is convex and satisfies~\eqref{as:lh} for every $\rho>0$ with $M_0=\norm{b}_1c^3e^{c\rho}$ and $M_1=c^3e^{c\rho}$, where $c=\max_i\norm{a_i}$, but it does not satisfy~\eqref{as:gh}.
\end{example}

\begin{example}[$\ell_p$ regression]
\label{ex:lp_lh_not_gh}
Let $p>3$, $A\in\R^{n\times d}$, $A\neq0$, and $b\in\R^n$. The $\ell_p$ regression loss $f(x)=\tfrac1p\norm{Ax-b}_p^p$ is convex and satisfies~\eqref{as:lh} for every $\rho>0$ with 
\begin{equation*}
  M_0=(p-1)(p-2)\norm{A}^3\ls1+2\norm{b}_p+\norm{A}\rho\rs^{p-3}, \qquad  M_1=3\sqrt{n}M_0/\sigma_+(A),
\end{equation*}
where $\sigma_+(A)$ is the smallest positive singular value of $A$, but it does not satisfy~\eqref{as:gh}.
\end{example}

For $f\in C^3(\R^d)$, condition~\eqref{as:lh} implies the pointwise bound~\citep{xie2024trust}
\begin{equation*}
  \norm{\nabla^3 f(x)}_{\mathrm{op}}
  \eqdef \sup_{\norm{u}=1}\norm{\nabla^3 f(x)[u]}
  \leq M_0+ M_1\norm{\grad f(x)},
  \qquad x\in\R^d.
  \tag{TD}\label{as:td}
\end{equation*}

\begin{proposition}
\label{prop:td_to_lh}
Suppose that, for some $L_0>0$ and $L_1\geq0$,
$\norm{\nabla^2 f(x)} \leq L_0+L_1\norm{\grad f(x)}$.
If~\eqref{as:td} holds with constants $\widehat M_0$ and
$\widehat M_1$, then, for every $\rho>0$,~\eqref{as:lh} holds with
\begin{equation*}
    M_0(\rho)=\widehat M_0+
  \widehat M_1L_0\rho e^{L_1\rho},\qquad
  M_1(\rho)=\widehat M_1e^{L_1\rho}.
\end{equation*}
\end{proposition}

\subsection{Globally Convergent Second-order Methods under Inexact Hessians}

With an inexact Hessian, cubic regularization commonly uses the model
\begin{equation}
  \label{eq:cubic_model_main}
  \omega_x^{\M,\hdel}(s)
  \eqdef f(x)+\la\grad f(x),s\ra
  +\tfrac12\la H(x)s,s\ra
  +\tfrac{\hdel}{2}\norm{s}^2
  +\tfrac{\M}{6}\norm{s}^3,
\end{equation}
with $s = y - x$.
The quadratic term compensates for Hessian error and the cubic term controls the Taylor remainder. Under~\eqref{as:h}, $\hdel\geq\delta$ and $\M \geq L$ ensure model
convexity and majorization~\citep[Theorem 1]{agafonov2023inexact}. In Nesterov-type acceleration these coefficients also enter the estimating sequence~\citep{nesterov2008accelerating,ghadimi2017second}. 
\\[2pt] 
For a parameter-free method, a useful observation is that the first-order optimality condition of the cubic model depends on these two regularizers only through their sum:
\begin{equation*}
  (H(x) + \tau I)s = - \nabla f(x), \qquad \tau = \hdel + \tfrac{\M}{2}\|s\|.
\end{equation*}
This suggests treating $\tau$ as a singe regularization parameter, rather than adapting $\hdel$ and $M$ separately. The difficulty here is that $\tau$ depends on unknown step length $\|s\|$. It is a circular dependence, which arises in the accelerated Newton proximal extragradient (A-NPE) method~\citep{monteiro2013accelerated} and its optimal variants~\citep{kovalev2022first,carmon2022optimal,chen2026optimal}. In particular, it can be resolved by searching over the regularization parameter while solving only quadratic subproblems. This perspective motivates our approach: instead of separately adapting two regularization parameters, we search directly for a regularized Newton step that satisfies a verifiable Monteiro-Svaiter condition~\citep{monteiro2013accelerated}.
\begin{definition}
  \label{def:ms}
  A Monteiro-Svaiter (MS) oracle at $x$ returns $y = x + s$ and $\tau > 0$ such that
  \begin{equation*}
    \|\nabla f(x + s) + \tau s\| \leq \tau \|s\| / 2.
  \end{equation*}
\end{definition}

Under~\eqref{as:h}, inexact cubic step $y = x + \argmin_{s} \omega_{x}^{\M, \hdel}(s)$ with $\hdel\geq 2\delta$ and $\M\geq 2L$ is an instance of MS oracle with $\tau = \hdel + \M\|s\|/2$~\citep[Lemma 3.1]{chen2026optimal}.
Under exact Hessian the same Lipschitz assumption, the gradient-regularized Newton step~\citep{mishchenko2023regularized,doikov2023gradient} is also an MS oracle for $\tau=\sqrt{M\norm{\grad f(x)}}>0$ and $M\geq L$.

\section{Basic Method}

\begin{figure}[H]
  \setlength{\intextsep}{0pt}
\begin{lrbox}{\backtrackingbox}
\begin{minipage}{0.48\textwidth}
  \begin{algorithmic}[1]
    \STATE \textbf{Input:} $x\in\R^d$ and $\tau>0$.
    \LOOP
      \IF{$\PD(x,\tau)$}
        \STATE Solve $(H(x)+\tau I)s=-\grad f(x)$.
        \IF{$\MS(x,\tau,s)$}
          \RETURN $(x+s,\tau)$.
        \ENDIF
      \ENDIF
      \STATE Set $\tau:=2\tau$. \label{line:increase_tau}
    \ENDLOOP
  \end{algorithmic}
\end{minipage}
\end{lrbox}
\begin{lrbox}{\basicnewtonbox}
\begin{minipage}{0.48\textwidth}
  \begin{algorithmic}[1]
    \STATE \textbf{Input:} $x_0\in\R^d$ and $\eta_0 > 0$.
    \FOR{$k\geq0$}
      \STATE $(x_{k+1},\tau_k)\eqdef\MSBacktrack(x_k,\eta_k)$.
      \STATE Set $\eta_{k+1}\eqdef\frac{\tau_k}{2}
        \min\lb1,\frac{\norm{\grad f(x_{k+1})}}{\norm{\grad f(x_k)}}\rb$.\label{line:basic_pd_ms_eta_update}
    \ENDFOR
  \end{algorithmic}
\end{minipage}
\end{lrbox}
\setlength{\pairedalgorithmheight}{\dimexpr\ht\backtrackingbox+\dp\backtrackingbox\relax}
\ifdim\dimexpr\ht\basicnewtonbox+\dp\basicnewtonbox\relax>\pairedalgorithmheight
  \setlength{\pairedalgorithmheight}{\dimexpr\ht\basicnewtonbox+\dp\basicnewtonbox\relax}
\fi
\begin{minipage}[t]{0.48\textwidth}
  \vspace{0pt}
\begin{algorithm}[H]
  \caption{\strut $\MSBacktrack(x,\tau)$}
  \label{alg:pd_ms_oracle}
  \parbox[t][\pairedalgorithmheight][t]{\linewidth}{\usebox{\backtrackingbox}}
\end{algorithm}
\end{minipage}\hfill
\begin{minipage}[t]{0.48\textwidth}
  \vspace{0pt}
\begin{algorithm}[H]
  \caption{\strut Adaptive regularized Newton}
  \label{alg:basic_pd_ms}
  \parbox[t][\pairedalgorithmheight][t]{\linewidth}{\usebox{\basicnewtonbox}}
\end{algorithm}
\end{minipage}
\end{figure}

\label{sec:basic}
We begin by proposing a basic non-accelerated adaptive method under MS oracle. First, we define two tests for our backtracking procedure
\begin{equation}
    \PD(x,\tau): H(x)+\tau I\succ0, \quad
    \MS(x,\tau,s): \norm{\grad f(x+s)+\tau s}\leq\tfrac{\tau}{2}\norm{s}.
  \label{eq:pd_ms_tests}
\end{equation}
The tests do not require $\delta$, $L(x)$, or $\rho$.
The $\PD$ test is checked first, so the linear system is solved
only when $H(x)+\tau I$ is positive definite. Algorithm~\ref{alg:pd_ms_oracle} defines $\MSBacktrack$, which ensures that tests~\eqref{eq:pd_ms_tests} hold. The resulting Adaptive regularized Newton methods is listed as Algorithm~\ref{alg:basic_pd_ms}. 

At each outer iteration, the method evaluates one $H(x)$ 
and keeps it fixed throughout backtracking. For example, an SVD-based implementation can reuse the same factorization across all trial values of $\tau$.
The gradient $\grad f(x_k)$ is reused throughout backtracking, and
$\grad f(x_{k+1})$ is reused in the update of $\eta_{k+1}$ and at the next iteration.
The halving in line~\ref{line:basic_pd_ms_eta_update} prevents overregularization, which is important in our convergence analysis (see proof sketch below). Consequently, $\tau$ is allowed to both increase and decrease across outer iterations. The gradient-dependent scaling factor is not needed for convex convergence, but allows us to recover local quadratic convergence under exact Hessians. 

\subsection{Convex analysis}

Let $D=\max_{x\in\mathcal{L}(x_0)}\norm{x-x_0}$ be the radius of the level set
$\mathcal{L}(x_0)=\lb x\in\R^d:f(x)\leq f(x_0)\rb$.
Since Algorithm~\ref{alg:basic_pd_ms} is monotone (Lemma~\ref{lem:basic_pd_ms_step}),
all iterates remain in $\mathcal{L}(x_0)$. Define
\begin{equation}
  B\eqdef\max\lb1, 2D/\rho,6\sqrt{M_1}D\rb,
  \qquad
  \e_{\tr}\eqdef\max\lb{M_0D^3}/{B^2},{\delta D^2}/{B}\rb,
  \label{eq:basic_pd_ms_scales}
\end{equation}
where $\e_{\tr}$ marks the switch from the linear rate to the rate of inexact cubic Newton. 
\begin{theorem}
\label{thm:basic_pd_ms_convergence}
Suppose that $f$ is convex,~\eqref{as:inexact_hessian} holds, and
either~\eqref{as:gh} or~\eqref{as:lh} holds.
Let
$\ell_0\eqdef\max\lb1,\tfrac{2D\norm{\grad f(x_0)}}{\e_{\tr}},
      \tfrac{\eta_0D^2}{\sqrt6 B\e_{\tr}}\rb$.
For every $\e>0$, Algorithm~\ref{alg:basic_pd_ms} generates an iterate $x_N$
satisfying $f(x_N)-f^\ast\leq\e$ within the following number of iterations
\begin{equation}
  N=\cO\ls
    1+\log \ell_0
    +B\log_+\tfrac{f(x_0)-f^\ast}{\max\lb\e,\e_{\tr}\rb}
      +\sqrt{\tfrac{M_0D^3}{\e}}+\tfrac{\delta D^2}{\e}
    \rs.
  \label{eq:basic_pd_ms_complexity}
\end{equation}
\end{theorem}

Theorem~\ref{thm:basic_pd_ms_convergence} bounds the number of outer iterations $N$, which equals the number of inexact Hessian evaluations. In Appendix~\ref{sec:pd_ms_oracle}, we show that the total number of
gradient evaluations is $\cO(N+\log_+\tfrac{B \|\nabla f(x_0)\|}{D\eta_0})$,
matching $N$ up to a constant factor and an additive logarithmic term.\\[2pt] 
The complexity~\eqref{eq:basic_pd_ms_complexity} separates the dependence on $M_0, M_1$ and $\delta$, and reveals two accuracy regimes. If the target accuracy $\e \geq \e_\tr$, the method remains in the far regime and converges linearly, requiring $\cO(B+\log \ell_0
+B\log_+\tfrac{f(x_0)-f^\ast}{\e})$ iterations. If $\e < \e_\tr$, the method first reaches $\e_{\tr}$ at a linear rate, after which the complexity is determined by $M_0$ and $\delta$, matching the rate of inexact cubic Newton under~\eqref{as:h}. The additional term $\log \ell_0$ is independent of $\e$ and arises
from bounding the number of overregularized iterations (see the proof sketch). Let us introduce one of the main lemmas.

\begin{lemma}
\label{lem:ahpe_damped_movement_main}
Suppose that $f$ is convex,~\eqref{as:inexact_hessian} holds, and
either~\eqref{as:gh} or~\eqref{as:lh} holds. 
Let $(y,\tau) \eqdef \MSBacktrack(x, \eta)$.
If at least one trial is rejected (equivalently $\tau>\eta$), then 
\begin{equation*}
  \norm{s}>\min\lb\tfrac\rho2,\tfrac1{6\sqrt{M_1}}\rb~~\text{or}~~\tau\leq8\delta+24M_0\norm{s}
\end{equation*}
\end{lemma}

This shows that for iterations with rejected trials we have two regimes: either regularization matches the scale of inexact cubic Newton step (set $M_0 = L$ under~\eqref{as:h} and see example after Definition~\ref{def:ms}) or the step is determined by generalized smoothness constants $M_1, \rho$ and is sufficiently large.\\[4pt]
\textbf{Proof sketch of Theorem~\ref{thm:basic_pd_ms_convergence}.}
Denote $\Delta_k = f(x_k) - f(x^*)$, $g_k = \| \nabla f(x_k) \|$, and $s_k = x_{k+1} - x_k$. By $\MS$ test we show $\tau_k \|s_k\| / 2 \leq g_{k+1} \leq 3 \tau_k\|s_k\|/2$ and using convexity 
\begin{equation*}
  f(x_k) - f(x_{k+1}) \geq 2g_{k+1}^2 / (3\tau_k), \qquad g_{k+1} \leq 2 g_k/\sqrt{3}.
\end{equation*}

Thus, the method is monotone and progress of one step depends on how large the accepted regularization $\tau_k$ can be relative to $g_{k+1}$. We call iteration $k$ controlled if $\tau_k\leq\Theta_{\rm c}(g_{k+1})$, where $\Theta_{\rm c}(g)\eqdef\max\lb8\delta+\sqrt{48M_0g},4g/\rho,12\sqrt{M_1}g\rb$. The proof is based on separate counting of controlled and uncontrolled iterations.\\[2pt]
First, we show, that if at least one trial was rejected the iteration is controlled. By Lemma~\ref{lem:ahpe_damped_movement_main} and $\MS$ relation 
we have 
$\tau_k\leq\Theta_{\rm c}(g_{k+1})$. Let $U_k\eqdef \max\lb\tau_{k-1},\Theta_{\rm c}(\norm{\grad f(x_k)})\rb$. Then by the definition of $U_k$
we get $U_{k+1} \leq (2/\sqrt{3})U_k$ on every controlled iteration.\\[2pt]
An uncontrolled iteration has no rejected trials and therefore accepts $\tau_k=\eta_k$. By the update rule for $\eta_k$, we obtain $U_{k+1} = \tau_k \leq \tau_{k-1}/2 \leq U_k/2$. Therefore $\Theta_{\rm c}(g_m) \leq U_0(2/\sqrt{3})^{m_{\rm c}}2^{-(m-m_{\rm c})}$ after $m$ iteration of which $m_{\rm c}$ are controlled. Taking the logarithm yields 
\begin{equation*}
    m \leq 2 m_{\rm c} + \log_2(U_0 / \Theta_{\rm c}(g_m)).
\end{equation*}
Thus, to show the convergence, it is enough to count the number of controlled steps. For controlled step, using the boundness of $\tau_k$, convexity and monotonicity we show recurrence 
\begin{equation*}
  \Delta_{k} - \Delta_{k+1} \gtrsim \tfrac{\Delta_{k+1}}{B + \sqrt{M_0D^3/\Delta_{k+1}} + \delta D^2/\Delta_{k+1}}.
\end{equation*}
 When $B$ dominates in denominator ($\Delta_{k+1} \geq \e_{\tr}$) we get linear convergence, otherwise the rate is determined by $M_0$ and $\delta$ terms. To accumulate this progress, we introduce a scalar potential that increases by at least one on every controlled iteration and does not decrease on the others. This gives us a bound 
\begin{equation*}
  m_{\rm c} \lesssim B \log_+ {\Delta_0}/{\max\{\e, \e_\tr\}} + \sqrt{{M_0D^3}/{\e}} + {\delta D^2}/{\e} 
\end{equation*}
as long as $\Delta_m > \e$. Bounding the logarithmic term in upper bound on $m$ gives us the rate~\eqref{eq:basic_pd_ms_complexity}. The complete proof is given in Appendix~\ref{app:basic_method_proof}. 
\subsection{Strongly convex analysis}
\label{sec:basic_pd_ms_strong}

In this section, we denote $\delta_k\eqdef\norm{H(x_k)-\nabla^2f(x_k)}\leq\delta$ and $r_0\eqdef\min\lb\rho/2,1/(6\sqrt{M_1})\rb$.

\begin{theorem}
\label{thm:basic_pd_ms_strong_main}
Suppose that $f$ is $\mu$-strongly convex,~\eqref{as:inexact_hessian} holds,
and either~\eqref{as:gh} or~\eqref{as:lh} holds. Let
\begin{equation*}
  \overline\tau\eqdef\max\lb
    \eta_0,8\delta+48M_0\sqrt{{2(f(x_0)-f^\ast)}/\mu},
    {2(f(x_0)-f^\ast)}/{r_0^2}\rb.
\end{equation*}
Every completed iteration $k$ of Algorithm~\ref{alg:basic_pd_ms}
satisfies $\tau_k\leq\overline\tau$ and 
\begin{equation*}
  \tau_k\leq8\max_{0\leq j\leq k}2^{j-k}\delta_j
    +\overline\tau\ls 1-\tfrac{2\mu}{3\overline\tau+2\mu}\rs^k, \quad
  f(x_{k+1})-f^\ast
  \leq\ls1-\tfrac{2\mu}{3\tau_k+2\mu}\rs^2
    (f(x_k)-f^\ast).
\end{equation*}
\end{theorem}
Theorem~\ref{thm:basic_pd_ms_strong_main} guarantees global linear
convergence for any sequence $\delta_k \leq \delta$ and shows that Hessian errors do not accumulate. If $\delta_k\to 0$, then $\tau_k\to 0$, so the contraction factor in convergence rate
goes to zero and the method converges superlinearly. In particular, for exact Hessians, we obtain global superlinear
convergence as established for cubic Newton in~\citet{kamzolov2025optami}. This theorem also covers lazy Hessian updates~\citep{doikov2023second},
where an exact Hessian is reused for a fixed number of iterations.
Indeed, for any $m\geq1$, the choice
$H(x_k)=\nabla^2f(x_{m\lfloor k/m\rfloor})$
satisfies $\delta_k\to0$.\\[2pt]
The rate in Theorem~\ref{thm:basic_pd_ms_strong_main} bounds the number of outer iterations and, equivalently, Hessian computations. To bound the first-order oracle calls, let $N$ be the first index satisfying 
$f(x_N) - f^\ast \leq \e$. By Corollary~\ref{cor:basic_pd_ms_strong_gradient_complexity} it can be bounded as 
$\cO(N
+\log_{2,+}\ls\tfrac{\norm{\grad f(x_0)}}{\eta_0}
\max\lb\tfrac{M_0}{\mu}+\tfrac{\delta}{\sqrt{\mu\e}},
\tfrac1{r_0}\rb\rs)$.
For inexact Hessians, this bound guarantees only linear convergence
in gradient evaluations, even when $\delta_k\to0$ 
gives superlinear
convergence in Hessian evaluations.
For $\delta=0$, the additional logarithmic term
is independent of $\e$, and the gradient count is $\cO(N+1)$.
Thus, for exact Hessians, the method converges superlinearly in terms of both gradient and Hessian evaluations.\\[2pt]
Our next goal is to recover the local quadratic convergence typical
for second-order methods with exact Hessians. We show that Algorithm~\ref{alg:basic_pd_ms} preserves local quadratic and global superlinear convergence even with inexact Hessians.
Both rates hold in terms of gradient and Hessian evaluations under
a stronger accuracy assumption than~\eqref{as:inexact_hessian}~\citep[Examples~7,~8]{bellavia2020subsampled,bellavia2022adaptive,chen2022accelerating,semenov2026gradient}:
\begin{equation}
  \delta_k\leq\nu\norm{\grad f(x_k)},\qquad k\geq0.
  \label{as:basic_pd_ms_relative_global}
\end{equation}

\begin{theorem}
\label{thm:basic_pd_ms_relative_main}
Suppose that $f$ is $\mu$-strongly convex, either~\eqref{as:gh} or~\eqref{as:lh} holds, and~\eqref{as:basic_pd_ms_relative_global} holds for some $\nu>0$.
Algorithm~\ref{alg:basic_pd_ms} converges globally superlinearly.
Every completed iteration $k$ satisfies
\begin{equation*}
  \tau_k\leq
  \max\lb\tfrac{2^{-k}\eta_0}{\norm{\grad f(x_0)}},
    \tfrac{24M_0}{\mu}+8\nu,\tfrac2{r_0}\rb
  \norm{\grad f(x_k)},~ 
  f(x_{k+1})-f^\ast
  \leq\ls1-\tfrac{2\mu}{3\tau_k+2\mu}\rs^2
    \ls f(x_k)-f^\ast\rs.
\end{equation*}
Let $C\eqdef
  \tfrac{9\ls\norm{\nabla^2f(x^\ast)}+\mu\rs}{2\mu^2}
  \ls8\nu+\tfrac{24M_0}{\mu}\rs^2$.
Let for some $K_0\geq1$, $f(x_{K_0})-f^\ast
  \leq\min\lb\tfrac{\mu r_0^2}{8},\tfrac1{2C}\rb$.
Then, after finitely many additional iterations independent of $\e$,
the method converges quadratically:
\begin{equation*}
  f(x_{k+1})-f^\ast
  \leq C\ls f(x_k)-f^\ast\rs^2,
  \qquad k\geq  K\geq K_0.
\end{equation*}
\end{theorem}

By Corollary~\ref{cor:basic_pd_ms_relative_complexity},
$N\geq1$ completed iterations require at most $N$ Hessian evaluations
and $\cO(N+\log_{2,+}\ls\tfrac{\norm{\grad f(x_0)}}{\eta_0}
  \max\lb\tfrac{M_0}{\mu}+\nu,\tfrac1{r_0}\rb\rs)$
gradient evaluations.
The logarithmic term is independent of $\e$, so the global superlinear
iteration complexity also holds for gradient evaluations.
Since our method converges globally, the index $K_0$ is reachable. 
After iteration $K$, which  is defined in Theorem~\ref{thm:basic_pd_ms_quadratic} and independent of $\e$, the method requires  $\cO(\log\log(1/\e))$ additional Hessian
and gradient evaluations suffice to reach $f(x)-f^\ast\leq\e$.

\section{Acceleration under Generalized Smoothness}
\label{sec:ahpe_acceleration}
\subsection{Suboptimal Acceleration}

Before moving to the optimal acceleration scheme, we briefly discuss generalization of accelerated methods with the suboptimal rate $\cO(1/k^3)$ under~\eqref{as:h}. Inspired by~\citet{tyurin2025near}, in Appendix~\ref{sec:fully_adaptive} we propose an adaptive Nesterov-type accelerated method with exact Hessians. The main challenge is that under~\eqref{as:lh}, the regularization parameter $\M_k$ must dominate $L(w_{k+1})$ at the output $w_{k+1}$ of the cubic step~\eqref{eq:cubic_model_main}, even though this point is unknown when $\M_k$ is chosen. We refer to the proof sketch in Appendix~\ref{sec:nesterov_acceleration} for the key idea behind resolving this difficulty. Furthermore, this acceleration scheme uses the cubic model, which has two parameters: $\M$ and $\hdel$. To avoid the difficulty of adapting two parameters simultaneously, in Appendix~\ref{app:ahpe_acceleration} we also propose an adaptive A-NPE method with inexact Hessians under~\eqref{as:lh}, which attains the rate $\cO(1/k^3 + \delta/k^2)$ under~\eqref{as:h}.

\subsection{Adaptive Damped A-NPE Method}
\label{sec:adaptive_damped_ahpe}

\begin{figure}[t]
\setlength{\intextsep}{0pt}
\begin{algorithm}[H]
  \caption{Adaptive damped A-NPE with an inexact Hessian}
  \label{alg:ahpe_damped}
  \begin{algorithmic}[1]
    \STATE \textbf{Input:} $x_0\in\R^d$ and $\eta_0>0$.
    \STATE $(w_1,\eta_1)\eqdef\MSBacktrack(x_0,\eta_0)$.\label{line:ahpe_damped_initialization}
    \STATE Set
        $A_1\eqdef 1/\eta_1,~
        z_1\eqdef x_0-A_1\grad f(w_1)$.
    \FOR{$k\geq1$}
      \STATE Set
          $a'_k\eqdef\tfrac{1+\sqrt{1+4\eta_kA_k}}{2\eta_k},~
          A'_k\eqdef A_k+a'_k,~
          v_k\eqdef\tfrac{A_kw_k+a'_kz_k}{A'_k}$.
      \STATE $(y_k,\tau_k)\eqdef\MSBacktrack(v_k,\eta_k)$.
      \STATE Set
          $\gamma_k\eqdef\tfrac{\eta_k}{\tau_k},~
          a_{k+1}\eqdef\gamma_ka'_k,~
          A_{k+1}\eqdef A_k+a_{k+1}.$
          \label{line:weights_damped_ahpe}
      \STATE Set
      $
          w_{k+1}\eqdef\tfrac{(1-\gamma_k)A_kw_k+\gamma_kA'_ky_k}{A_{k+1}},~
          z_{k+1}\eqdef z_k-a_{k+1}\grad f(y_k).
      $
      \IF{$\tau_k=\eta_k$}
        \STATE Set $\eta_{k+1}\eqdef \eta_k / 2$.
      \ELSE
        \STATE Set $\eta_{k+1}\eqdef2\eta_k$.
      \ENDIF
    \ENDFOR
  \end{algorithmic}
\end{algorithm}
\end{figure}

In this section, we present the adaptive damped A-NPE method, listed as Algorithm~\ref{alg:ahpe_damped}. 
The method is parameter-free and attains the optimal convergence rate for convex optimization with inexact Hessians in the classical Lipschitz-Hessian setting.
It combines backtracking (Algorithm~\ref{alg:pd_ms_oracle}) with acceleration scheme of~\citet{carmon2022optimal}. The parameter $\eta_k$ is the guess regularization parameter for $\MSBacktrack$ procedure, which returns $\tau_k \geq \eta_k$.  If $\tau_k = \eta_k$ the method performs full step with $A_{k+1} = A_{k} + a'_{k}$, otherwise the update of $A_{k+1}$ is damped $A_{k+1} = A_{k} + (\eta_k/\tau_k)a'_{k}$. The next guess is set to $\eta_k/2$ when the first trial is accepted and to $2\eta_k$ otherwise.\\[2pt]
As in Algorithm~\ref{alg:basic_pd_ms}, each outer iteration uses at most one Hessian, which is reused throughout backtracking. We again maintain one regularization parameter, which adapts to both Hessian inexactness and local smoothness, without requiring them as an input. The intuition behind this acceleration scheme can be found in the proof sketch below. While the outer scheme of Algorithm~\ref{alg:ahpe_damped} matches the method of~\citet{carmon2022optimal}, our convergence analysis is different. We discuss the key differences below.\\[2pt]
For the convergence bounds of Algorithm~\ref{alg:ahpe_damped}, define
\begin{equation*}
  B\eqdef\max\lb1,(6\sqrt{M_1}R)^{2/3},\ls 2R/\rho\rs^{2/3}\rb,
  \quad
  \e_{\tr}\eqdef
  \max\lb48M_0R^3/B^{7/2},16\delta R^2/B^2\rb,
\end{equation*}
where $\rho=+\infty$ under~\eqref{as:gh}. Let $\overline\tau_0\eqdef\max\lb4\delta,2\sqrt{L(x_0)\norm{\grad f(x_0)}},2\norm{\grad f(x_0)}/\rho\rb$.

\begin{theorem}
\label{thm:ahpe_damped_convergence}
Suppose that $f$ is convex,~\eqref{as:inexact_hessian} holds, and
either~\eqref{as:gh} or~\eqref{as:lh} holds.
For any $\e>0$, the method generates an iterate $w_N$ satisfying
$f(w_N)-f^\ast\leq\e$ within
\begin{equation}
  \cO\ls
    B\ls1+\log_+\tfrac{\max\lb\eta_0,\overline\tau_0\rb R^2}{\max\lb\e,\e_{\tr}\rb}\rs
    +\ls\tfrac{M_0R^3}{\e}\rs^{2/7}
    +\sqrt{\tfrac{\delta R^2}{\e}}
  \rs
  \label{eq:ahpe_damped_outer_complexity}
\end{equation}
iterations. Moreover, there exists iterate $K\geq1$ with
$K=\cO\ls B\max\lb1,\log\tfrac{\max\lb\eta_0,\overline\tau_0\rb R^2}{\e_{\tr}}\rb\rs$ such that
\begin{gather}
  f(w_k)-f^\ast
  =\cO\ls\max\lb\eta_0,\overline\tau_0\rb R^2\exp\ls-\tfrac{c(k-1)}B\rs\rs,
  \qquad 1\leq k<K,
  \label{eq:ahpe_damped_far_final}\\
  f(w_{K+n})-f^\ast
  =\cO\ls\tfrac{M_0R^3}{(B+n)^{7/2}}+\tfrac{\delta R^2}{(B+n)^2}\rs,
  \qquad n=0,1,\ldots.
  \label{eq:ahpe_damped_near_final}
\end{gather}
\end{theorem}

Theorem~\ref{thm:ahpe_damped_convergence} bounds outer iterations, which equals the number of inexact Hessian evaluations. In Appendix~\ref{sec:ahpe_damped_oracle}, we show that the number of gradient evaluations is of the same order as~\eqref{eq:ahpe_damped_outer_complexity} up to two additive terms $\cO(\log_+\tfrac{\overline\tau_0}{\eta_0})$ and $\widetilde{\cO}(B(1+\log_+\tfrac{\max\lb\eta_0,\overline\tau_0\rb R^2}{\e})^{2/3})$, where $\widetilde{\cO}$ hides logarithmic factor independent of $\e$. By~\eqref{eq:ahpe_damped_far_final} and~\eqref{eq:ahpe_damped_near_final}, the method first converges linearly, and then converges in Lipschitz Hessian regime.\\[2pt]
Under~\eqref{as:h} the complexity is $\cO(1+\log_+\tfrac{\max\lb\eta_0,\overline\tau_0\rb}{LR+\delta}+(\tfrac{LR^3}{\e})^{2/7}+(\tfrac{\delta R^2}{\e})^{1/2})$, matching the optimal method of~\citet{chen2026optimal} and the corresponding lower bound~\citep{agafonov2024advancing} up to an additive term independent of $\e$. Gradient computations have additional term $\cO(\log_+\tfrac{\overline\tau_0}{\eta_0})$. To the best of our knowledge, this is the first second-order method to simultaneously achieve full adaptivity, optimal complexity, and robustness to inexact Hessians. The complete proof is given in Appendix~\ref{app:ahpe_damped}.\\[4pt]
\textbf{Proof sketch.} Let $s_i \eqdef y_i - v_i$ and $\Lambda(A)\eqdef\max\lb32\delta,(96M_0R)^{4/7}A^{-3/7},B^2/A\rb$. By the algorithm construction, the regularization guess $\eta_i$ is driven toward $\Lambda(A_i)$. $\Lambda(A)$ serves as an ideal regularization scale, which will allow us to recover the optimal Lipschitz-Hessian rate
near the solution and a geometric rate far from it.\\[2pt]
We next use the standard Monteiro-Svaiter potential argument. Let $\ECal_k = A_k (f(w_k) - f^\ast) + 1/2\|z_k - x^\ast\|^2$. Then, $3/8\sum_{i=1}^{k-1} A_i'\eta_i \|s_i\|^2 \leq \ECal_1 \leq R^2/2$. 
We first count two types of iterations:
  \begin{gather*}
  \cS_{k-1}^{>}\eqdef
  \lb i\in\{1,\ldots,k-1\}:\tau_i>\eta_i,
       \ \eta_i\geq\Lambda(A_i)/2\rb,\\
  \cS_{k-1}^{=}\eqdef
  \lb i\in\{1,\ldots,k-1\}:\tau_i=\eta_i,
       \ \eta_i\leq2\Lambda(A_i)\rb.
\end{gather*}
On $\cS_{k-1}^{>}$, at least one trial is rejected and $\eta_i$ is
sufficiently large. 
Lemma~\ref{lem:ahpe_damped_movement_main}, together with $A_i' \geq A_i$ lower-bounds the corresponding terms in the potential estimate:
\begin{equation*}
  \tfrac{3R^2}{16 \sqrt{A_k \Lambda(A_k)}} |\cS_{k-1}^{>}|  \leq \tfrac{3}{8}\textstyle{\sum}_{i\in \cS_{k-1}^{>}} A_i'\eta_i \|s_i\|^2 \leq \tfrac{R^2}{2}\quad \Rightarrow \quad |\cS_{k-1}^{>}| \leq \tfrac{8}{3}\sqrt{A_k\Lambda(A_k)}.
\end{equation*}
Here, the condition on $\eta_i$ yields sufficient potential decrease, similarly to how controlled iterations in the basic method yield sufficient function decrease. Note, that substituting $A_k \asymp R^2/\e$ and definition of $\Lambda(A_k)$ into $\sqrt{A_k\Lambda(A_k)}$ recovers the desired rate~\eqref{eq:ahpe_damped_outer_complexity}.

On $\cS_{k-1}^{=}$, no trial is rejected, so $\gamma_i=1$ and $A_{i+1} = A_i + a_i'$. 
The upper bound on $\eta_i=\tau_i$ prevents overregularization. Using the update in line~\ref{line:weights_damped_ahpe}, we have 

\begin{equation}
  \label{eq:counting_down}
  \tfrac{a_i'}{A_i} \geq \tfrac{1}{\sqrt{2A_i\Lambda(A_i)}}~\Rightarrow~\tfrac13 \leq \sqrt{A_i\Lambda(A_i)}\log\tfrac{A_{i+1}}{A_i} \leq \textstyle{\int}_{A_i}^{A_{i+1}}\sqrt{\tfrac{\Lambda(t)}{t}}dt.
\end{equation}
Thus, on these iterations the upper bound on $\eta_i$  guarantees fast growth of $A_i$. Summing the last inequality over $\cS_{k-1}^{=}$ and calculating integral gives the upper bound~\eqref{eq:ahpe_damped_outer_complexity}. Note that $\Lambda$ is chosen to match the progress on  $\cS_{k-1}^{>}$ and $\cS_{k-1}^{=}$.\\[2pt]
Iterations that are not in $\cS_{k-1}^{>} \cup \cS_{k-1}^{=}$ improve regularization guess. If $\tau_i=\eta_i$ and $\eta_i>2\Lambda(A_i)$, the algorithm halves
$\eta_i$. Otherwise, $\tau_i>\eta_i$ and $\eta_i<\Lambda(A_i)/2$,
and it doubles $\eta_i$. The last challenge is to bound the number of such iterations. Similarly to uncontrolled iterations in the basic method, we show 
\begin{equation*}
    k - 1 \leq 2 (|\cS_{k-1}^{>}| + |\cS_{k-1}^{=}|) + \log (A_k\Lambda(A_1)).
\end{equation*}
Combining iteration counts yields 
\begin{equation*}
  k\lesssim
  B\log_+\min\{A_k,A_{\tr}\}/A_1
  +\sqrt{A_k\Lambda(A_k)},
\end{equation*}
where $A_{\tr}\eqdef R^2/(2\e_\tr)$.
Since $f(w_k)-f^\ast\leq R^2/(2A_k)$, it suffices to reach
$A_\e\eqdef R^2/(2\e)$.
For $\e\geq\e_{\tr}$, we have $A_\e\leq A_{\tr}$ and
$\sqrt{A_\e\Lambda(A_\e)}=B$, so the complexity depends
logarithmically on the target accuracy.
For $\e<\e_{\tr}$, the logarithmic term is constant, while 
\begin{equation*}
  \sqrt{A_\e\Lambda(A_\e)} \simeq \left(M_0R^3 /\e\right)^{2/7}
+\sqrt{\delta R^2 / \e}.
\end{equation*}

\subsection{Related Proof Techniques}
\label{sec:related_proofs}

Algorithmically, our method builds upon the optimal MS acceleration of~\citet{carmon2022optimal}. Recently,~\citet{chen2026optimal} established the optimal rate for inexact Hessians using the schedule of~\citet{adil2024convex}, but their method is non-adaptive. Below, we discuss the key challenges of extending their analyses to an adaptive method under~\eqref{as:lh} and~\eqref{as:inexact_hessian}.\\[4pt]
\textbf{Comparison with~\citet{carmon2022optimal}.} They split iterations into two sets: $\cS^{\uparrow}=\{i:\tau_i>\eta_i\}$, on which the potential bound $\sum_{i\in\cS^{\uparrow}}A'_i\eta_i\norm{s_i}^2\lesssim R^2$ is used, and $\cS^{\downarrow}=\{i:\tau_i=\eta_i\}$, on which the growth of $A_i$ is given by the weight identity $\eta_i(a'_i)^2=A'_i$. In the latter case, they use the bound of type $\sqrt{A_{i+1}}-\sqrt{A_i}\geq 1 / (2\sqrt{\eta_i})$, which is enough for polynomial rate. The first difference is that we use log form in~\eqref{eq:counting_down} to get geometric far rate. \\[2pt]
From this point the proofs follow different ideas.~\citet{carmon2022optimal} transfer the growth of $A_i$ from $\cS^{\downarrow}$ to neighboring iterations of $\cS^{\uparrow}$. The reason is that on $\cS^{\uparrow}$ they have  $\norm{s_i}\gtrsim\tau_i/L$ because of rejection, which turns the potential bound into $\sum_{i\in\cS^{\uparrow}}A'_i\eta_i^3\lesssim L^2R^2$. Next, since both bounds are on $\cS^{\uparrow}$, they use reverse H\"older inequality to remove $\eta_i$ and get the rate $A_k\gtrsim k^{7/2}/(LR)$. \\[2pt]
Because of assumptions~\eqref{as:lh} and~\eqref{as:inexact_hessian} we do not have such lower bound on $\|s_i\|$. Instead, Lemma~\ref{lem:ahpe_damped_movement_main} gives three regimes for rejected steps: it is proportional to $\tau_i$ for $\delta\lesssim\tau_i\lesssim M_0r_0$, it is bounded below only by $r_0$ for larger $\tau_i$, and it may be arbitrary small for $\tau_i\lesssim\delta$. Thus, we introduce $\Lambda(A)$, whose three terms correspond to these three regimes, and follow the proof of Theorem~\ref{thm:basic_pd_ms_convergence} instead. $\Lambda(A_i)$ plays the role of $\Theta_{\rm c}(g)$, iterations of $\cS^{>}_{k-1}$ and $\cS^{=}_{k-1}$ (because of additional bound on guess $\eta_i$) make progress in the potential or show the growth of $A_i$.  The remaining iterations move $\eta_i$ toward $\Lambda(A_i)$, so their number is bounded by the multiplier on number of steps with progress and additional logarithm (similarly to uncontrolled iterations in the basic method).\\[4pt]
\textbf{Comparison with~\citet{chen2026optimal}.}
The cubic step is used with explicit regularization ($\tau_i=2\delta+M_0\norm{s_i}$), and the guess is $\eta_i\simeq\delta+(M_0R)^{4/7}A_i^{-3/7}$~\citep[Lemma~3.1 and Theorem~3.1]{chen2026optimal}. Their acceleration scheme keeps $\eta_i$ fixed until $A_i$ doubles and then resets it. Under~\eqref{as:h} this matches our $\Lambda(A)$ up to constants for $A\geq A_{\tr}$. However, the analysis is based on the fact that between two resets the guess is a known constant of order at least $\delta$, while an adaptive guess changes at every iteration and can be arbitrarily small. We resolve this challenge by counting only  iterations with $\eta_i\asymp\Lambda(A_i)$ and bounding the remaining as described above.

\section{Experiments}
\label{sec:experiments}

\begin{figure}[t]
  \centering
  \includegraphics[width=\textwidth]{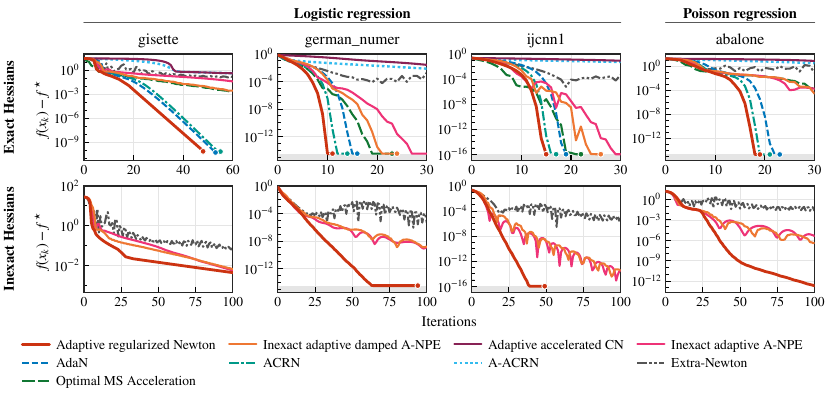}
  \caption{Top row: exact Hessians. Bottom row: Hessians computed on a fixed random 10\% of the samples, with exact gradients.}
  \label{fig:exp}
\end{figure}

We compare our methods with adaptive second-order baselines on logistic and Poisson regression (Figure~\ref{fig:exp}). We use datasets from LIBSVM~\citep{chang2011libsvm} \texttt{gisette}, \texttt{german\_numer} and \texttt{ijcnn1} for logistic regression, and \texttt{abalone} for Poisson regression. 
Methods start from $x_0=-\mathbf1$. %
We run Algorithms~\ref{alg:basic_pd_ms}, \ref{alg:ahpe_damped}, \ref{alg:fully_adaptive_cubic} and~\ref{alg:adaptive_ahpe} against AdaN~\citep{mishchenko2023regularized}, the adaptive cubic regularized Newton (ACRN), its accelerated version (A-ACRN)~\citep{kamzolov2026accelerated}, Extra-Newton~\citep{antonakopoulos2022extra}, and Optimal MS Acceleration~\citep{carmon2022optimal}\\[2pt]
\textbf{Parameters.}
All methods are adaptive, so we do not tune them. For a fair comparison, Algorithms~\ref{alg:basic_pd_ms},~\ref{alg:ahpe_damped} and~\ref{alg:adaptive_ahpe} start backtracking with $\eta_0=\sqrt{L(x_0)\norm{\nabla f(x_0)}}$ , Optimal MS Acceleration takes $\lambda'_0=\eta_0$, $\sigma=1/2$, and AdaN takes $H_0=L(x_0)/2$, which it doubles before its first trial. For the step-size parameters of Extra-Newton, we take $\gamma=2\norm{x_0}$ as a proxy for the diameter $D$, following~\citet[Remark~3.1]{antonakopoulos2022extra}, and $\beta_0=L(x_0)\gamma^2\norm{\nabla f(x_0)}$. Thus, these six methods perform the same first trial step. ACRN, A-ACRN and Algorithm~\ref{alg:fully_adaptive_cubic} start from the standard theoretical choice of the cubic regularization parameter, $\M_0=L(x_0)$.\\[2pt]
\textbf{Inexact Hessians.}
We consider a centralized distributed setup with $f(x)=\frac1m\sum_{j=1}^m f_j(x)$, where $f_j(x)$ is the local loss stored on machine $j$. The server also stores local data and computes the Hessian $H(x)=\nabla^2 f_1(x)$. Under the $\beta$-similarity this Hessian satisfies~\eqref{as:inexact_hessian} with $\delta=\beta$~\citep{zhang2015disco,agafonov2021accelerated}. To perform a step, the server collects the gradients from the workers and uses the full gradient with its local Hessian. In our experiments, $H(x)$ is computed from a fixed uniformly sampled subset of $\lceil n/10\rceil$ observations, where $n$ is the number of training samples. We compare only parameter-free methods that allow inexact Hessians.
\\[2pt]
\textbf{Results.}
Under exact Hessians, our Adaptive regularized Newton converges faster on all four problems, with ACRN close behind. Among the accelerated methods, our A-NPE methods (Algorithms~\ref{alg:ahpe_damped} and~\ref{alg:adaptive_ahpe}) perform best, and the optimal Algorithm~\ref{alg:ahpe_damped} usually outperforms the suboptimal Algorithm~\ref{alg:adaptive_ahpe}. As in the logistic regression experiments of~\citet{carmon2022optimal}, the non-accelerated methods outperform the accelerated ones. However, on toy problems satisfying~\eqref{as:lh}, optimal acceleration outperforms non-accelerated methods and suboptimal accelerated ones when starting far from the solution (see below). With inexact Hessians, our three methods outperform Extra-Newton, and Algorithm~\ref{alg:basic_pd_ms} is again the best. Experimental details and additional results are provided in Appendix~\ref{app:extra_experiments}.

\begin{figure}[t]
  \centering
  \includegraphics[width=\textwidth]{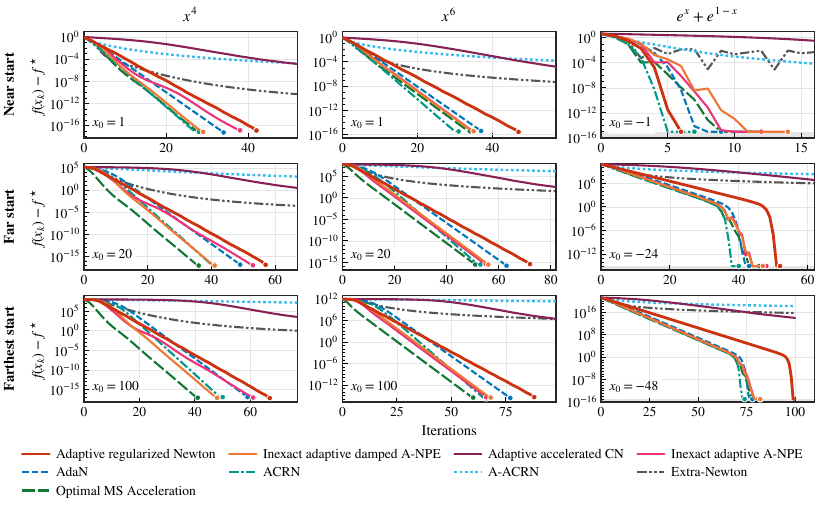}
  \caption{One-dimensional functions with exact Hessians, from the nearest start (top row) to the farthest (bottom row). Each panel ends ten iterations after Algorithm~\ref{alg:basic_pd_ms} stops.}
  \label{fig:toy1d}
\end{figure}

\textbf{One-dimensional functions satisfying~\eqref{as:lh}.}
We also run all the methods with exact Hessians on three one-dimensional functions from~\citet{gorbunov2025methods} and~\citet{tyurin2025near} (Figure~\ref{fig:toy1d}): $x^4$, $x^6$ and $e^x+e^{1-x}$. Up to scaling and a shift, $x^4$ and $e^x+e^{1-x}$ are Examples~\ref{ex:quartic_lh_not_gh} and~\ref{ex:cosh_lh_not_gh}, which satisfy~\eqref{as:lh} but not~\eqref{as:gh}. We start $x^4$ and $x^6$ from $x_0\in\{1,20,100\}$ and $e^x+e^{1-x}$ from $x_0\in\{-1,-24,-48\}$, compute derivatives in closed form, and set all parameters from $L(x_0)=M_0+M_1|f'(x_0)|$ by the rules above, where~\eqref{as:td} holds with $(M_0,M_1)=(16,2)$, $(48,12)$ and $(0,1)$, respectively. Every run has $100$ iterations.

Six methods converge from every start: Algorithms~\ref{alg:basic_pd_ms}, \ref{alg:ahpe_damped} and~\ref{alg:adaptive_ahpe}, AdaN, ACRN and Optimal MS Acceleration~\citep{carmon2022optimal}. On $x^4$ and $x^6$ they converge linearly, as expected at a degenerate minimizer. On $e^x+e^{1-x}$ they converge in a few iterations once close to $x^\ast$, after a far phase whose length grows linearly with $|x_0-x^\ast|$. Unlike in Figure~\ref{fig:exp}, acceleration pays off here: from far starts, the accelerated Monteiro--Svaiter-type methods, Algorithms~\ref{alg:ahpe_damped} and~\ref{alg:adaptive_ahpe} and Optimal MS Acceleration, converge faster than the non-accelerated Algorithm~\ref{alg:basic_pd_ms}.

\bibliography{agafonov}
\bibliographystyle{iclr2027_conference}

\newpage
\appendix
\tableofcontents
\addtocontents{toc}{\protect\setcounter{tocdepth}{2}}

\newpage
\section{Taylor Remainder Lemmas}
\label{app:taylor_remainders}
The following two lemmas record the Taylor remainder bounds used throughout
our analysis.

\begin{lemma}
\label{lem:global_taylor_remainders}
Suppose that \eqref{as:gh} holds. Then, for all $x,s\in\R^d$,
\begin{gather*}
  \norm{\grad f(x+s)-\grad f(x)-\nabla^2 f(x)s}
  \leq \tfrac{L(x)}{2}\norm{s}^2,\\
  \left|
    f(x+s)-f(x)-\la\grad f(x),s\ra
    -\tfrac{1}{2}\la\nabla^2 f(x)s,s\ra
  \right|
  \leq \tfrac{L(x)}{6}\norm{s}^3.
\end{gather*}
\end{lemma}

\begin{lemma}
\label{lem:local_taylor_remainders}
Suppose that \eqref{as:lh} holds. Then, for all $x,s\in\R^d$ satisfying
$\norm{s}\leq\rho$,
\begin{gather*}
  \norm{\grad f(x+s)-\grad f(x)-\nabla^2 f(x)s}
  \leq \tfrac{L(x)}{2}\norm{s}^2,\\
  \left|
    f(x+s)-f(x)-\la\grad f(x),s\ra
    -\tfrac{1}{2}\la\nabla^2 f(x)s,s\ra
  \right|
  \leq \tfrac{L(x)}{6}\norm{s}^3.
\end{gather*}
\end{lemma}

\begin{proof}
We follow~\citet[Lemma~1.2.4]{nesterov2018lectures}. For $t\in[0,1]$,~\eqref{as:gh}, or~\eqref{as:lh} with $\norm{ts}\leq\rho$, gives $\norm{\nabla^2 f(x+ts)-\nabla^2 f(x)}\leq tL(x)\norm{s}$. Hence,
\begin{gather*}
  \norm{\grad f(x+s)-\grad f(x)-\nabla^2 f(x)s}
  =\norm{\int_0^1\ls\nabla^2 f(x+ts)-\nabla^2 f(x)\rs s\mathrm{d}t}\\
  \leq\int_0^1 tL(x)\norm{s}^2\mathrm{d}t
  =\frac{L(x)}{2}\norm{s}^2.
\end{gather*}
Since $\norm{ts}\leq\norm{s}$, the same bound holds for the step $ts$:
\begin{equation*}
  \norm{\grad f(x+ts)-\grad f(x)-t\nabla^2 f(x)s}
  \leq\frac{L(x)}{2}t^2\norm{s}^2,
  \qquad t\in[0,1].
\end{equation*}
Using $f(x+s)-f(x)=\int_0^1\la\grad f(x+ts),s\ra\mathrm{d}t$ and $\int_0^1t\mathrm{d}t=\frac{1}{2}$, we obtain
\begin{gather*}
  \left|
    f(x+s)-f(x)-\la\grad f(x),s\ra
    -\frac{1}{2}\la\nabla^2 f(x)s,s\ra
  \right|\\
  =\left|
    \int_0^1\la\grad f(x+ts)-\grad f(x)-t\nabla^2 f(x)s,s\ra\mathrm{d}t
  \right|
  \leq\int_0^1\frac{L(x)}{2}t^2\norm{s}^3\mathrm{d}t
  =\frac{L(x)}{6}\norm{s}^3.
\end{gather*}
\end{proof}

\section{Proofs for the Relations between Assumptions and Examples}
\label{app:assumption_relations}

\begin{proof}[Proof of Example~\ref{ex:gh_not_h}]
Direct differentiation gives
\begin{gather*}
  f'(x)=2x+\int_0^x\sin(t^2)\mathrm{d}t,\qquad
  f''(x)=2+\sin(x^2),\qquad
  f^{(3)}(x)=2x\cos(x^2).
\end{gather*}
Thus, $f''(x)\geq1$, so $f$ is $1$-strongly convex. Moreover,
$f'(0)=0$ and
\begin{equation*}
  |f'(x)|
  =\left|\int_0^x f''(t)\mathrm{d}t\right|
  \geq |x|.
\end{equation*}
If $|y-x|\leq1$, then
\begin{gather*}
  |f''(y)-f''(x)|
  =|\sin(y^2)-\sin(x^2)|
  \leq |y^2-x^2|
  \leq \ls2|x|+1\rs|y-x|
  \leq \ls2+2|f'(x)|\rs|y-x|.
\end{gather*}
If $|y-x|>1$, then
\begin{equation*}
  |f''(y)-f''(x)|
  \leq2
  \leq2|y-x|
  \leq\ls2+2|f'(x)|\rs|y-x|.
\end{equation*}
Hence,~\eqref{as:gh} holds with $M_0=M_1=2$. On the other hand, for
$x_n=\sqrt{2\pi n}$,
\begin{equation*}
  |f^{(3)}(x_n)|=2\sqrt{2\pi n}\longrightarrow+\infty.
\end{equation*}
Condition~\eqref{as:h} would imply $|f^{(3)}(x)|\leq L$ by taking the
directional limit in~\eqref{as:h}. Therefore,~\eqref{as:h} does not hold.
\end{proof}

\begin{proof}[Proof of Example~\ref{ex:cosh_lh_not_gh}]
We have
\begin{gather*}
  f'(x)=\sinh x,\qquad
  f''(x)=\cosh x,\qquad
  f^{(3)}(x)=\sinh x.
\end{gather*}
Since $f''(x)\geq1$, the function is $1$-strongly convex. Let
$h=y-x$ and $|h|\leq\rho$. For every $t\in[0,1]$,
\begin{gather*}
  |\sinh(x+th)|
  \leq |\sinh x|\cosh(th)+\cosh x|\sinh(th)|\\
  \leq |\sinh x|\cosh\rho+\cosh x\sinh\rho
  \leq \sinh\rho+e^\rho|\sinh x|,
\end{gather*}
where we used $\cosh x\leq1+|\sinh x|$. Therefore,
\begin{equation*}
  |f''(y)-f''(x)|
  \leq |h|\int_0^1|\sinh(x+th)|\mathrm{d}t 
  \leq\ls\sinh\rho+e^\rho|f'(x)|\rs|y-x|.
\end{equation*}
This proves~\eqref{as:lh}. If~\eqref{as:gh} held, then the choice $x=0$
would give
\begin{equation*}
  \cosh y-1\leq M_0|y|,
  \qquad y\in\R,
\end{equation*}
which is impossible as $y\to+\infty$.
\end{proof}

\begin{proof}[Proof of Example~\ref{ex:exp_lh_not_gh}]
We have
\begin{gather*}
  f'(x)=e^x-1,\qquad
  f''(x)=f^{(3)}(x)=e^x.
\end{gather*}
Let $h=y-x$ and $|h|\leq\rho$. The mean value theorem gives
\begin{equation*}
  |f''(y)-f''(x)|
  =e^x|e^h-1|
  \leq e^\rho e^x|h|
  \leq e^\rho\ls1+|f'(x)|\rs|y-x|.
\end{equation*}
Hence,~\eqref{as:lh} holds with $M_0=M_1=e^\rho$. If~\eqref{as:gh}
held, then the choice $x=0$ would give
\begin{equation*}
  e^y-1\leq M_0|y|,
  \qquad y\in\R,
\end{equation*}
which is impossible as $y\to+\infty$.
\end{proof}

\begin{proof}[Proof of Example~\ref{ex:quartic_lh_not_gh}]
We have
\begin{gather*}
  \grad f(x)=\ls\mu+\norm{x}^2\rs x,\qquad
  \nabla^2f(x)=\ls\mu+\norm{x}^2\rs I+2xx^\top.
\end{gather*}
Hence, $\nabla^2f(x)\succeq\mu I$. Let $s=y-x$ and
$r=\norm{s}\leq\rho$. Using
\begin{equation*}
  yy^\top-xx^\top=xs^\top+sx^\top+ss^\top,
\end{equation*}
we obtain
\begin{gather*}
  \norm{\nabla^2f(y)-\nabla^2f(x)}
  \leq\left|\norm{y}^2-\norm{x}^2\right|
  +2\norm{yy^\top-xx^\top}\\
  \leq\ls2\norm{x}+r\rs r
  +2\ls2\norm{x}r+r^2\rs
  =\ls6\norm{x}+3r\rs r.
\end{gather*}
For every $t\geq0$,
\begin{equation*}
  4+2t^3-6t=2(t-1)^2(t+2)\geq0.
\end{equation*}
Since
$\norm{\grad f(x)}=\ls\mu+\norm{x}^2\rs\norm{x}\geq\norm{x}^3$,
we have
\begin{equation*}
  6\norm{x}\leq4+2\norm{\grad f(x)}.
\end{equation*}
Consequently,
\begin{equation*}
  \norm{\nabla^2f(y)-\nabla^2f(x)}
  \leq\ls3\rho+4+2\norm{\grad f(x)}\rs\norm{y-x},
\end{equation*}
which proves~\eqref{as:lh} with the stated constants. At $x=0$,
\begin{equation*}
  \norm{\nabla^2f(y)-\nabla^2f(0)}=3\norm{y}^2,
  \qquad \grad f(0)=0.
\end{equation*}
Thus,~\eqref{as:gh} would require $3\norm{y}^2\leq M_0\norm{y}$ for
all $y\in\R^d$, which is impossible.
\end{proof}

\begin{proof}[Proof of Example~\ref{ex:poisson_lh_not_gh}]
Let $\mu_i(x)\eqdef e^{\la a_i,x\ra}$. We have
\begin{gather*}
  \grad f(x)=\sum_{i=1}^n\ls\mu_i(x)-b_i\rs a_i,\qquad
  \nabla^2f(x)=\sum_{i=1}^n\mu_i(x)a_ia_i^\top\succeq0.
\end{gather*}
Since $[a_i]_1=1$ and $b_i\geq0$, the first coordinate of $\grad f(x)$ equals $\sum_{i=1}^n\mu_i(x)-\norm{b}_1$. Hence,
\begin{equation}
  \sum_{i=1}^n\mu_i(x)\leq\norm{b}_1+\norm{\grad f(x)}.
  \label{eq:poisson_intensity}
\end{equation}
Let $s=y-x$ and $\norm{s}\leq\rho$. Then $|\la a_i,s\ra|\leq c\norm{s}\leq c\rho$. Using $|e^u-1|\leq e^{|u|}|u|$ for $u\in\R$, we obtain
\begin{gather*}
  \norm{\nabla^2f(y)-\nabla^2f(x)}
  =\left\|\sum_{i=1}^n\mu_i(x)\ls e^{\la a_i,s\ra}-1\rs a_ia_i^\top\right\|
  \leq c^3e^{c\rho}\norm{s}\sum_{i=1}^n\mu_i(x)\\
  \stackrel{\eqref{eq:poisson_intensity}}{\leq}
  \ls\norm{b}_1c^3e^{c\rho}+c^3e^{c\rho}\norm{\grad f(x)}\rs\norm{y-x}.
\end{gather*}
This proves~\eqref{as:lh} with the stated constants. Let $e_1$ be the first coordinate vector. Since $\la a_i,e_1\ra=1$, we have $\nabla^2f(x+te_1)=e^t\nabla^2f(x)$ for $t\in\R$ and
\begin{equation*}
  \norm{\nabla^2f(x)}
  \geq\la\nabla^2f(x)e_1,e_1\ra
  =\sum_{i=1}^n\mu_i(x)>0.
\end{equation*}
If~\eqref{as:gh} held, then, for all $t>0$,
\begin{equation*}
  \ls e^t-1\rs\norm{\nabla^2f(x)}
  =\norm{\nabla^2f(x+te_1)-\nabla^2f(x)}
  \leq\ls M_0+M_1\norm{\grad f(x)}\rs t,
\end{equation*}
which is impossible as $t\to+\infty$.
\end{proof}

\begin{proof}[Proof of Example~\ref{ex:lp_lh_not_gh}]
Let $a_i^\top$ be the rows of $A$, $r\eqdef Ax-b$, and $\psi(u)\eqdef|u|^{p-2}u$ for $u\in\R$. We have
\begin{gather*}
  \grad f(x)=\sum_{i=1}^n\psi(r_i)a_i,\qquad
  \nabla^2f(x)=(p-1)\sum_{i=1}^n|r_i|^{p-2}a_ia_i^\top\succeq0.
\end{gather*}
First, we bound $\norm{r}_p$ by $\norm{\grad f(x)}$. Let $P$ be the orthogonal projector onto the range of $A^\top$. The gradient $\grad f(x)=A^\top\psi(r)$ belongs to this range, and $\norm{Az}\geq\sigma_+(A)\norm{z}$ for all $z$ in this range. Therefore, H\"older's inequality and $\norm{v}\leq\sqrt{n}\norm{v}_\infty\leq\sqrt{n}\norm{v}_p$ for $v\in\R^n$ give
\begin{gather*}
  \norm{r}_p^p-\norm{r}_p^{p-1}\norm{b}_p
  \leq\sum_{i=1}^n\psi(r_i)\ls r_i+b_i\rs
  =\la\grad f(x),x\ra
  =\la\grad f(x),Px\ra
  \leq\frac{\norm{\grad f(x)}\norm{APx}}{\sigma_+(A)}\\
  =\frac{\norm{\grad f(x)}\norm{r+b}}{\sigma_+(A)}
  \leq\frac{\sqrt{n}}{\sigma_+(A)}\norm{\grad f(x)}\ls\norm{r}_p+\norm{b}_p\rs.
\end{gather*}
If $\norm{r}_p>\max\lb1,2\norm{b}_p\rb$, then $\norm{r}_p-\norm{b}_p>\norm{r}_p/2$ and $\norm{r}_p+\norm{b}_p<3\norm{r}_p/2$. In this case, the previous display gives $\norm{r}_p^{p-1}<\frac{3\sqrt{n}}{\sigma_+(A)}\norm{\grad f(x)}=\frac{M_1}{M_0}\norm{\grad f(x)}$, and $\norm{r}_p>1$ implies
\begin{gather*}
  \ls\norm{r}_p+\norm{A}\rho\rs^{p-3}
  \leq(1+\norm{A}\rho)^{p-3}\norm{r}_p^{p-3}
  \leq(1+\norm{A}\rho)^{p-3}\norm{r}_p^{p-1}\\
  <\ls1+2\norm{b}_p+\norm{A}\rho\rs^{p-3}\frac{M_1}{M_0}\norm{\grad f(x)}.
\end{gather*}
Otherwise, $\norm{r}_p\leq1+2\norm{b}_p$ and $\ls\norm{r}_p+\norm{A}\rho\rs^{p-3}\leq\ls1+2\norm{b}_p+\norm{A}\rho\rs^{p-3}$. In both cases, the definition of $M_0$ gives
\begin{equation}
  (p-1)(p-2)\norm{A}^3\ls\norm{r}_p+\norm{A}\rho\rs^{p-3}
  \leq M_0+M_1\norm{\grad f(x)}.
  \label{eq:lp_residual}
\end{equation}
Now let $s=y-x$ and $\norm{s}\leq\rho$. Then $Ay-b=r+As$ and $|\la a_i,s\ra|\leq\norm{A}\norm{s}\leq\norm{A}\rho$. By the mean value theorem and $|r_i|\leq\norm{r}_p$,
\begin{gather*}
  \left||r_i+\la a_i,s\ra|^{p-2}-|r_i|^{p-2}\right|
  \leq(p-2)\ls|r_i|+|\la a_i,s\ra|\rs^{p-3}|\la a_i,s\ra|\\
  \leq(p-2)\ls\norm{r}_p+\norm{A}\rho\rs^{p-3}\norm{A}\norm{s}.
\end{gather*}
Since $\norm{\sum_{i=1}^nw_ia_ia_i^\top}\leq\max_i|w_i|\norm{A^\top A}$ for $w\in\R^n$, we obtain
\begin{gather*}
  \norm{\nabla^2f(y)-\nabla^2f(x)}
  =(p-1)\left\|\sum_{i=1}^n\ls|r_i+\la a_i,s\ra|^{p-2}-|r_i|^{p-2}\rs a_ia_i^\top\right\|\\
  \leq(p-1)(p-2)\norm{A}^3\ls\norm{r}_p+\norm{A}\rho\rs^{p-3}\norm{s}
  \stackrel{\eqref{eq:lp_residual}}{\leq}
  \ls M_0+M_1\norm{\grad f(x)}\rs\norm{y-x}.
\end{gather*}
This proves~\eqref{as:lh}. Since $A\neq0$, there is $j$ with $a_j\neq0$. For $t\geq|b_j|/\norm{a_j}^2$,
\begin{equation*}
  \norm{\nabla^2f(ta_j)}\norm{a_j}^2
  \geq\la\nabla^2f(ta_j)a_j,a_j\ra
  \geq(p-1)\ls t\norm{a_j}^2-|b_j|\rs^{p-2}\norm{a_j}^4.
\end{equation*}
If~\eqref{as:gh} held, then the choice $x=0$ and $y=ta_j$ would give
\begin{equation*}
  \norm{\nabla^2f(ta_j)}
  \leq\norm{\nabla^2f(0)}+\ls M_0+M_1\norm{\grad f(0)}\rs t\norm{a_j},
  \qquad t>0.
\end{equation*}
Since $p-2>1$, this contradicts the previous display as $t\to+\infty$.
\end{proof}

\begin{proof}[Pointwise consequence of the pairwise conditions]
For any $x\in\R^d$, $\norm{u}=1$, and all sufficiently small $|t|$,
either pairwise condition gives
\begin{gather*}
  \frac{\norm{\nabla^2 f(x+tu)-\nabla^2 f(x)}}{|t|}
  \leq M_0+M_1\norm{\grad f(x)}.
\end{gather*}
Taking $t\to0$ and then the supremum over $\norm{u}=1$
proves~\eqref{as:td}.
\end{proof}

The pointwise first-order condition used in
Proposition~\ref{prop:td_to_lh} is denoted by
\begin{equation*}
  \norm{\nabla^2 f(x)}
  \leq L_0+L_1\norm{\grad f(x)},
  \qquad x\in\R^d.
  \tag{FO}\label{as:fo}
\end{equation*}

\begin{proof}[Proof of Proposition~\ref{prop:td_to_lh}]
Let $s=y-x$, $r=\norm{s}\leq\rho$, and
$g(t)=\norm{\grad f(x+ts)}$. The function $g$ is absolutely continuous.
For almost every $t\in[0,1]$, condition~\eqref{as:fo} gives
\begin{gather*}
  g'(t)
  \leq r\norm{\nabla^2 f(x+ts)}
  \leq r\ls L_0+L_1g(t)\rs.
\end{gather*}
Gr\"onwall's inequality yields
\begin{gather*}
  g(t)
  \leq e^{L_1rt}\ls g(0)+L_0rt\rs
  \leq e^{L_1\rho}\ls\norm{\grad f(x)}+L_0\rho\rs.
\end{gather*}
Using~\eqref{as:td} along the segment from $x$ to $y$, we obtain
\begin{gather*}
  \norm{\nabla^2 f(y)-\nabla^2 f(x)}
  \leq r\int_0^1
  \ls\widehat M_0+\widehat M_1g(t)\rs\mathrm{d}t 
  \leq\ls
    \widehat M_0+\widehat M_1L_0\rho e^{L_1\rho}
    +\widehat M_1e^{L_1\rho}\norm{\grad f(x)}
  \rs r.
\end{gather*}
This is~\eqref{as:lh} with the stated constants.
\end{proof}

\section{Proof of Algorithm~\ref{alg:basic_pd_ms}}
\label{app:basic_method_proof}

We consider completed iterations before a zero-gradient return.

For $g\geq0$ and $k\geq0$, recall
\begin{equation}
\begin{aligned}
  \Theta_{\rm c}(g)&\eqdef\max\lb8\delta+\sqrt{48M_0g},\frac{4g}{\rho},12\sqrt{M_1}g\rb,\\
  U_k&\eqdef\max\lb\tau_{k-1},\Theta_{\rm c}(\norm{\grad f(x_k)})\rb,
  \qquad \tau_{-1}\eqdef2\eta_0.
\end{aligned}
  \label{eq:basic_pd_ms_threshold}
\end{equation}

\subsection{Step Estimates}

\begin{lemma}
\label{lem:basic_pd_ms_step}
Suppose that $f$ is convex. For every accepted step of
Algorithm~\ref{alg:basic_pd_ms}, set $s_k\eqdef x_{k+1}-x_k$. Then
\begin{gather}
  \frac{\tau_k}{2}\norm{s_k}
  \leq\norm{\grad f(x_{k+1})}
  \leq\frac{3\tau_k}{2}\norm{s_k},
  \label{eq:basic_pd_ms_output_gradient}\\
  f(x_k)-f(x_{k+1})
  \geq\max\lb\frac{\tau_k}{2}\norm{s_k}^2,
              \frac{2\norm{\grad f(x_{k+1})}^2}{3\tau_k}\rb,
  \label{eq:basic_pd_ms_decrease}\\
  \norm{\grad f(x_{k+1})}
  \leq\frac2{\sqrt3}\norm{\grad f(x_k)}.
  \label{eq:basic_pd_ms_gradient_growth}
\end{gather}
\end{lemma}

\begin{proof}
By the $\MS$ condition in~\eqref{eq:pd_ms_tests} 
\begin{gather*}
  \norm{\grad f(x_{k+1})}
  \geq\tau_k\norm{s_k}-\norm{\grad f(x_{k+1})+\tau_ks_k}
  \geq\frac{\tau_k}{2}\norm{s_k},\\
  \norm{\grad f(x_{k+1})}
  \leq\norm{\grad f(x_{k+1})+\tau_ks_k}+\tau_k\norm{s_k}
  \leq\frac{3\tau_k}{2}\norm{s_k}.
\end{gather*}
we prove ~\eqref{eq:basic_pd_ms_output_gradient}.
Squaring the $\MS$ condition and expanding the left-hand side gives
\begin{equation}
  \norm{\grad f(x_{k+1})+\tau_ks_k}^2
  =\norm{\grad f(x_{k+1})}^2
    +2\tau_k\la\grad f(x_{k+1}),s_k\ra
    +\tau_k^2\norm{s_k}^2
  \leq\frac{\tau_k^2}{4}\norm{s_k}^2.
  \label{eq:basic_pd_ms_squared_test}
\end{equation}
Next, by convexity,~\eqref{eq:basic_pd_ms_squared_test}, and~\eqref{eq:basic_pd_ms_output_gradient},
\begin{gather*}
  f(x_k)-f(x_{k+1})
  \geq-\la\grad f(x_{k+1}),s_k\ra
  \geq\frac{\norm{\grad f(x_{k+1})}^2}{2\tau_k}
    +\frac{3\tau_k}{8}\norm{s_k}^2\\
  \geq\frac{1}{2\tau_k}\ls\frac{\tau_k}{2}\norm{s_k}\rs^2
    +\frac{3\tau_k}{8}\norm{s_k}^2
  =\frac{\tau_k}{2}\norm{s_k}^2.
\end{gather*}
Using instead $\norm{s_k}\geq2\norm{\grad f(x_{k+1})}/(3\tau_k)$
from~\eqref{eq:basic_pd_ms_output_gradient}, we obtain
\begin{gather*}
  f(x_k)-f(x_{k+1})
  \geq\frac{\norm{\grad f(x_{k+1})}^2}{2\tau_k}
    +\frac{3\tau_k}{8}
       \ls\frac{2\norm{\grad f(x_{k+1})}}{3\tau_k}\rs^2
  =\frac{2\norm{\grad f(x_{k+1})}^2}{3\tau_k}.
\end{gather*}
This proves~\eqref{eq:basic_pd_ms_decrease}.
For~\eqref{eq:basic_pd_ms_gradient_growth}, convexity at both endpoints gives
\begin{equation*}
  \la\grad f(x_k),s_k\ra
  \leq f(x_{k+1})-f(x_k)
  \leq\la\grad f(x_{k+1}),s_k\ra.
\end{equation*}
Next, \eqref{eq:basic_pd_ms_squared_test} gives
\begin{equation*}
  \norm{\grad f(x_{k+1})}^2+\frac34\tau_k^2\norm{s_k}^2
  \leq-2\tau_k\la\grad f(x_{k+1}),s_k\ra
  \leq-2\tau_k\la\grad f(x_k),s_k\ra
  \leq2\tau_k\norm{\grad f(x_k)}\norm{s_k}.
\end{equation*}
Rearranging and completing the square gives
\begin{gather*}
  \norm{\grad f(x_{k+1})}^2
  \leq2 \tau_k \norm{\grad f(x_k)}\norm{s_k}
       -\frac34\tau_k^2\norm{s_k}^2\\
  =\frac43\norm{\grad f(x_k)}^2
   -\frac34\ls\tau_k\norm{s_k}-\frac43\norm{\grad f(x_k)}\rs^2
  \leq\frac43\norm{\grad f(x_k)}^2.
\end{gather*}
\end{proof}

\begin{lemma}
\label{lem:pd_ms_backtracking}
Suppose that $f$ is convex,~\eqref{as:inexact_hessian} holds, and
either~\eqref{as:gh} or~\eqref{as:lh} holds.
Set
\begin{equation}
  \label{eq:backtrack_upper_bound}
  \overline\tau(x)\eqdef
  \begin{cases}
    \max\lb4\delta,2\sqrt{L(x)\norm{\grad f(x)}}\rb,
    & \text{under~\eqref{as:gh}},\\
    \max\lb4\delta,2\sqrt{L(x)\norm{\grad f(x)}},
      \frac{2\norm{\grad f(x)}}{\rho}\rb,
    & \text{under~\eqref{as:lh}}.
  \end{cases}
\end{equation}
Every trial of $\MSBacktrack$ with
$\tau\geq\overline\tau(x)$ satisfies both tests in~\eqref{eq:pd_ms_tests}.
Consequently, initialized with $\tau=\eta>0$, the algorithm terminates
after at most
\begin{equation*}
  1+\left\lceil\log_2\max\lb1,\frac{\overline\tau(x)}{\eta}\rb\right\rceil
\end{equation*}
trials, and its accepted parameter satisfies
\begin{equation*}
  \eta\leq\tau\leq\max\lb\eta,2\overline\tau(x)\rb.
\end{equation*}
\end{lemma}

\begin{proof}
Fix a trial parameter $\tau\geq\overline\tau(x)$.
By convexity and~\eqref{as:inexact_hessian}, $H(x)\succeq-\delta I$.
Since $\tau\geq4\delta$,
\begin{equation*}
  H(x)+\tau I\succeq(\tau-\delta)I
  \succeq\frac{3\tau}{4}I\succ0.
\end{equation*}
Thus $\PD(x,\tau)$ holds, and the linear system
$(H(x)+\tau I)s=-\grad f(x)$ has a unique solution satisfying
\begin{equation*}
  \norm{s}
  =\norm{(H(x)+\tau I)^{-1}\grad f(x)}
  \leq\frac{\norm{\grad f(x)}}{\tau-\delta}
  \leq\frac{2\norm{\grad f(x)}}{\tau}.
\end{equation*}
Under~\eqref{as:lh}, the definition of $\overline\tau(x)$ also gives
$\norm{s}\leq\rho$.
The gradient remainder bound in Lemma~\ref{lem:global_taylor_remainders}
or Lemma~\ref{lem:local_taylor_remainders}, together
with~\eqref{as:inexact_hessian} and the linear system, yields
\begin{gather*}
  \norm{\grad f(x+s)+\tau s}
  =\norm{\grad f(x+s)-\grad f(x)-H(x)s}\\
  \leq\frac{L(x)}2\norm{s}^2+\delta\norm{s}
  \leq\ls\frac{L(x)\norm{\grad f(x)}}{\tau}+\delta\rs\norm{s}\\
  \leq\ls\frac\tau4+\frac\tau4\rs\norm{s}
  =\frac\tau2\norm{s}.
\end{gather*}
Hence $\MS(x,\tau,s)$ also holds.
The trial parameters are $2^j\eta$, $j=0,1,\ldots$, so the first index
with $2^j\eta\geq\overline\tau(x)$ gives the stated trial bound.
If no trial is rejected, the accepted parameter is $\tau=\eta$.
Otherwise, $\tau/2$ was rejected, so $\tau/2<\overline\tau(x)$ and
$\tau<2\overline\tau(x)$.
\end{proof}

\begin{lemma}
\label{lem:ahpe_damped_movement}
Suppose that $f$ is convex,~\eqref{as:inexact_hessian} holds, and
either~\eqref{as:gh} or~\eqref{as:lh} holds. Set
\begin{equation}
  r_0\eqdef
  \begin{cases}
    \frac1{6\sqrt{M_1}},
    & \text{under~\eqref{as:gh}},\\
    \min\lb\frac\rho2,\frac1{6\sqrt{M_1}}\rb,
    & \text{under~\eqref{as:lh}},
  \end{cases}
  \label{eq:ahpe_damped_movement_radius}
\end{equation}
where $1/\sqrt{M_1}=+\infty$ when $M_1=0$.
Suppose that $\MSBacktrack$, initialized at $x$ with
$\tau=\eta>0$, returns $(y,\tau)$. Set $s\eqdef y-x$.
If at least one trial is rejected (equivalently $\tau>\eta$), then
\begin{equation}
  \norm{s}>r_0
  \quad\text{or}\quad
  \tau\leq8\delta+24M_0\norm{s}.
  \label{eq:ahpe_damped_movement}
\end{equation}
If, in addition, $\tau\geq16\delta$, then
\begin{equation}
  \norm{s}\geq\min\lb r_0,\frac{\tau}{48M_0}\rb.
  \label{eq:ahpe_damped_movement_lower_bound}
\end{equation}
\end{lemma}

\begin{proof}
We prove~\eqref{eq:ahpe_damped_movement} by considering three cases.

First, if $\tau<4\delta$, then $\tau\leq8\delta+24M_0\norm{s}$,
so the second inequality holds.

Second, if $\tau\geq4\delta$ and $\norm{s}>r_0$, the first inequality holds.

It remains to consider the case $\tau\geq4\delta$ and $\norm{s}\leq r_0$.
Since $\tau>\eta$, the trial with parameter $\tau/2$ was rejected.
Convexity and~\eqref{as:inexact_hessian} give
\begin{equation}
  H(x)+\frac{\tau}{2}I
  \succeq\ls\frac{\tau}{2}-\delta\rs I
  \succeq\frac{\tau}{4}I\succ0,
  \label{eq:ahpe_damped_rejected_matrix_bound}
\end{equation}
so this trial failed the $\MS$ test. Its step $s_{1/2}$ and the accepted
step $s$ satisfy
\begin{gather*}
  \ls H(x)+\frac{\tau}{2}I\rs s_{1/2}=-\grad f(x),\\
  \ls H(x)+\tau I\rs s=-\grad f(x).
\end{gather*}
Subtracting these equations gives
\begin{equation}
  s_{1/2}-s
  =\frac{\tau}{2}\ls H(x)+\frac{\tau}{2}I\rs^{-1}s,
  \qquad
  \norm{s_{1/2}-s}\leq2\norm{s},
  \label{eq:ahpe_damped_trial_difference}
\end{equation}
where the last inequality follows from~\eqref{eq:ahpe_damped_rejected_matrix_bound}.
For every $t\in[0,1]$,~\eqref{eq:ahpe_damped_trial_difference} gives
\begin{gather*}
  \norm{x+ts_{1/2}-y}
  =\norm{t(s_{1/2}-s)-(1-t)s}\\
  \leq t\norm{s_{1/2}-s}+(1-t)\norm{s}
  \leq(1+t)\norm{s}\leq2\norm{s}.
\end{gather*}
Under~\eqref{as:lh},~\eqref{eq:ahpe_damped_movement_radius} gives
\begin{equation*}
  \norm{x-y}=\norm{s}\leq r_0\leq\frac\rho2, \qquad
  \norm{x+ts_{1/2}-y}\leq2\norm{s}\leq2r_0\leq\rho.
\end{equation*}
Thus, applying~\eqref{as:gh} or~\eqref{as:lh} at $y$ gives
\begin{gather}
  \norm{\nabla^2f(x+ts_{1/2})-\nabla^2f(x)}
  \leq\norm{\nabla^2f(x+ts_{1/2})-\nabla^2f(y)}
  +\norm{\nabla^2f(x)-\nabla^2f(y)} \notag\\
  \leq L(y)\ls\norm{x+ts_{1/2}-y}+\norm{x-y}\rs
  \leq3L(y)\norm{s}.
\label{eq:ahpe_damped_hessian_variation}
\end{gather}
Since the trial with parameter $\tau/2$ fails the $\MS$ test, we have
\begin{equation}
  \frac{\tau}{4}\norm{s_{1/2}}
  <\norm{\grad f(x+s_{1/2})+\frac{\tau}{2}s_{1/2}}.
  \label{eq:ahpe_damped_failed_ms_test}
\end{equation}
The linear system for $s_{1/2}$ gives
\begin{equation*}
  \frac{\tau}{2}s_{1/2}=-\grad f(x)-H(x)s_{1/2}.
\end{equation*}
Therefore,
\begin{gather}
  \grad f(x+s_{1/2})+\frac{\tau}{2}s_{1/2}
  =\grad f(x+s_{1/2})-\grad f(x)-H(x)s_{1/2} \notag \\
  =\int_0^1\ls\nabla^2f(x+ts_{1/2})-H(x)\rs s_{1/2}dt.
\label{eq:ahpe_damped_residual_integral}
\end{gather}
Therefore, by~\eqref{eq:ahpe_damped_failed_ms_test},~\eqref{eq:ahpe_damped_residual_integral},~\eqref{eq:ahpe_damped_hessian_variation}, and~\eqref{as:inexact_hessian},
\begin{equation*}
  \frac{\tau}{4}\norm{s_{1/2}}
  <\norm{\int_0^1\ls\nabla^2f(x+ts_{1/2})-H(x)\rs s_{1/2}dt}
  \leq\ls\delta+3L(y)\norm{s}\rs\norm{s_{1/2}}.
\end{equation*}
Since this trial fails $\MS$, its step $s_{1/2}$ is nonzero.
The accepted $\MS$ test in~\eqref{eq:pd_ms_tests} gives
\begin{equation*}
  \norm{\grad f(y)}
  \leq\norm{\grad f(y)+\tau s}+\tau\norm{s}
  \leq\frac{3\tau}{2}\norm{s}.
\end{equation*}
Dividing by $\norm{s_{1/2}}$ and using this bound, we obtain
\begin{gather*}
  \frac{\tau}{4}
  <\delta+3L(y)\norm{s}
  =\delta+3M_0\norm{s}+3M_1\norm{\grad f(y)}\norm{s}\\
  \leq\delta+3M_0\norm{s}+\frac92M_1\tau\norm{s}^2
  \leq\delta+3M_0\norm{s}+\frac{\tau}{8}.
\end{gather*}
The last inequality follows from $\norm{s}\leq r_0$ and~\eqref{eq:ahpe_damped_movement_radius}.
Indeed, if $M_1>0$, then
\begin{equation*}
  M_1\norm{s}^2\leq M_1r_0^2
  \leq M_1\ls\frac1{6\sqrt{M_1}}\rs^2=\frac1{36}.
\end{equation*}
If $M_1=0$, then $M_1\norm{s}^2=0$. Thus, in both cases,
\begin{equation*}
  \frac92M_1\tau\norm{s}^2
  \leq \frac{\tau}{8}.
\end{equation*}
Rearranging proves~\eqref{eq:ahpe_damped_movement}.
Finally, if $\tau\geq16\delta$,~\eqref{eq:ahpe_damped_movement} gives
\begin{equation*}
  \norm{s}
  \geq\min\lb r_0,\frac{\tau-8\delta}{24M_0}\rb
  \geq\min\lb r_0,\frac{\tau}{48M_0}\rb,
\end{equation*}
which proves~\eqref{eq:ahpe_damped_movement_lower_bound}.
\end{proof}

By~\eqref{eq:basic_pd_ms_threshold} and~\eqref{eq:ahpe_damped_movement_radius},
\begin{equation*}
  \Theta_{\rm c}(g)=\max\lb8\delta+\sqrt{48M_0g},\frac{2g}{r_0}\rb,
  \qquad g\geq0,
\end{equation*}
where $1/(+\infty)=0$.

\begin{lemma}
\label{lem:basic_pd_ms_rejection}
Suppose that $f$ is convex,~\eqref{as:inexact_hessian} holds, and
either~\eqref{as:gh} or~\eqref{as:lh} holds.
If iteration $k$ of Algorithm~\ref{alg:basic_pd_ms} rejects at least one trial,
then
\begin{equation}
  \tau_k\leq\Theta_{\rm c}(\norm{\grad f(x_{k+1})}).
  \label{eq:basic_pd_ms_rejection}
\end{equation}
\end{lemma}

\begin{proof}
Set $s_k\eqdef x_{k+1}-x_k$ and apply Lemma~\ref{lem:ahpe_damped_movement}
with $x=x_k$, $y=x_{k+1}$, $\eta=\eta_k$, and $\tau=\tau_k$.
If $\norm{s_k}>r_0$, then~\eqref{eq:basic_pd_ms_output_gradient}
and~\eqref{eq:basic_pd_ms_threshold} give
\begin{equation*}
  \tau_k<\frac{2\norm{\grad f(x_{k+1})}}{r_0}
  \leq\Theta_{\rm c}(\norm{\grad f(x_{k+1})}).
\end{equation*}
If $\norm{s_k}\leq r_0$,~\eqref{eq:ahpe_damped_movement}
and~\eqref{eq:basic_pd_ms_output_gradient} imply
\begin{gather*}
  \tau_k\leq8\delta+24M_0\norm{s_k}
  \leq8\delta+\frac{48M_0\norm{\grad f(x_{k+1})}}{\tau_k},\\
  \tau_k\leq4\delta+\sqrt{16\delta^2+48M_0\norm{\grad f(x_{k+1})}}
  \leq8\delta+\sqrt{48M_0\norm{\grad f(x_{k+1})}},
\end{gather*}
which proves~\eqref{eq:basic_pd_ms_rejection}.
\end{proof}

\subsection{Iteration Count}

We call a completed iteration $k$ of Algorithm~\ref{alg:basic_pd_ms} controlled if
\begin{equation}
  \tau_k\leq\Theta_{\rm c}(\norm{\grad f(x_{k+1})}).
  \label{eq:basic_pd_ms_controlled}
\end{equation}
By Lemma~\ref{lem:basic_pd_ms_rejection}, every iteration with a rejection
is controlled. An iteration accepted at its first trial can belong to
either type. 

\begin{lemma}
\label{lem:basic_pd_ms_tracking}
Under the assumptions of Lemma~\ref{lem:basic_pd_ms_rejection}, suppose that
$m$ iterations of Algorithm~\ref{alg:basic_pd_ms} are completed and $f(x_m)>f^\ast$.
If $m_{\rm c}$ of these iterations are controlled, then
\begin{equation}
  m\leq\frac54m_{\rm c}
   +\log_2\frac{U_0}{\Theta_{\rm c}(\norm{\grad f(x_m)})}.
  \label{eq:basic_pd_ms_tracking_count}
\end{equation}
\end{lemma}

\begin{proof}
By~\eqref{eq:basic_pd_ms_threshold},
$\Theta_{\rm c}(\norm{\grad f(x_k)})\leq U_k$ for every $k\geq0$.
The function $\Theta_{\rm c}$ is nondecreasing, and~\eqref{eq:basic_pd_ms_threshold} gives
$\Theta_{\rm c}(ag)\leq a\Theta_{\rm c}(g)$ for $g\geq0$ and $a\geq1$.

If iteration $k$ is controlled, then $\tau_k\leq\Theta_{\rm c}(\norm{\grad f(x_{k+1})})$
by~\eqref{eq:basic_pd_ms_controlled}. The definition of $U_{k+1}$
and~\eqref{eq:basic_pd_ms_gradient_growth} give
\begin{equation*}
  U_{k+1}=\Theta_{\rm c}(\norm{\grad f(x_{k+1})})
  \leq\frac2{\sqrt3}\Theta_{\rm c}(\norm{\grad f(x_k)})
  \leq\frac2{\sqrt3}U_k.
\end{equation*}
If iteration $k$ is not controlled, then
$\tau_k>\Theta_{\rm c}(\norm{\grad f(x_{k+1})})$, so $U_{k+1}=\tau_k$ by definition.
By Lemma~\ref{lem:basic_pd_ms_rejection}, this iteration has no rejected trials,
so $\tau_k=\eta_k$.
For $k\geq1$, line~\ref{line:basic_pd_ms_eta_update} of
Algorithm~\ref{alg:basic_pd_ms}, applied at iteration $k-1$, gives
$\eta_k\leq\tau_{k-1}/2$.
Hence,
\begin{equation*}
  U_{k+1}=\tau_k=\eta_k\leq\frac{\tau_{k-1}}2\leq\frac{U_k}{2}.
\end{equation*}
For $k=0$, the same bound holds because
$U_1=\tau_0=\eta_0\leq U_0/2$. Applying these bounds over the $m_{\rm c}$ controlled iterations and
$m-m_{\rm c}$ remaining iterations gives
\begin{equation*}
  \Theta_{\rm c}(\norm{\grad f(x_m)})\leq U_m
  \leq U_0\ls\frac2{\sqrt3}\rs^{m_{\rm c}}2^{-(m-m_{\rm c})}.
\end{equation*}
Taking logarithms and rearranging, we obtain
\begin{equation*}
  m\leq\ls1+\log_2\frac2{\sqrt3}\rs m_{\rm c}
       +\log_2\frac{U_0}{\Theta_{\rm c}(\norm{\grad f(x_m)})}.
\end{equation*}
Since $(2/\sqrt3)^4=16/9<2$, this
proves~\eqref{eq:basic_pd_ms_tracking_count}.
\end{proof}

To count the controlled iterations, assume $D<+\infty$.
By~\eqref{eq:basic_pd_ms_decrease}, every accepted point remains in
$\mathcal{L}(x_0)$, which also contains $x^\ast$. Hence, by convexity,
\begin{equation}
\begin{aligned}
  \norm{x_k-x^\ast}&\leq\norm{x_k-x_0}+\norm{x_0-x^\ast}\leq2D,\\
  f(x_k)-f^\ast
  &\leq\la\grad f(x_k),x_k-x^\ast\ra
  \leq2D\norm{\grad f(x_k)}.
\end{aligned}
  \label{eq:basic_pd_ms_localization}
\end{equation}
By~\eqref{eq:basic_pd_ms_scales} and~\eqref{eq:ahpe_damped_movement_radius},
$B=\max\lb1,D/r_0\rb$, where $D/(+\infty)=0$.
For $t>0$, define
\begin{equation}
  \Phi(t)\eqdef\max\lb\frac{6D}{r_0},
    12\sqrt{6}\sqrt{\frac{M_0D^3}{t}}+\frac{48\delta D^2}{t}\rb.
  \label{eq:basic_pd_ms_count_scale}
\end{equation}
The function $\Phi$ is positive and nonincreasing, and $t(1+\Phi(t))$ is nondecreasing.

\begin{lemma}
\label{lem:basic_pd_ms_controlled_count}
For the prefix and count of Lemma~\ref{lem:basic_pd_ms_tracking}, assume
in addition that $0<D<+\infty$. Then
\begin{equation}
  m_{\rm c}\leq\int_{f(x_m)-f^\ast}^{f(x_0)-f^\ast}\frac{1+\Phi(t)}t dt.
  \label{eq:basic_pd_ms_controlled_count}
\end{equation}
\end{lemma}
\begin{proof}
For a controlled iteration, we have
\begin{gather*}
  f(x_k)-f(x_{k+1})
  \stackrel{\eqref{eq:basic_pd_ms_decrease}}{\geq}
  \frac{2\norm{\grad f(x_{k+1})}^2}{3\tau_k}
  \stackrel{\eqref{eq:basic_pd_ms_controlled}}{\geq}
  \frac{2\norm{\grad f(x_{k+1})}^2}
       {3\Theta_{\rm c}(\norm{\grad f(x_{k+1})})}\\
  \stackrel{\eqref{eq:basic_pd_ms_threshold}}{=}
  \min\lb
      \frac{2\norm{\grad f(x_{k+1})}^2}
           {3\ls8\delta+\sqrt{48M_0\norm{\grad f(x_{k+1})}}\rs},
      \frac{r_0\norm{\grad f(x_{k+1})}}3
     \rb.
\end{gather*}
By convexity and Cauchy-Schwarz,
\begin{equation*}
  f(x_{k+1})-f^\ast
  \leq\la\grad f(x_{k+1}),x_{k+1}-x^\ast\ra
  \leq\norm{\grad f(x_{k+1})}\norm{x_{k+1}-x^\ast}
  \stackrel{\eqref{eq:basic_pd_ms_localization}}{\leq}
  2D\norm{\grad f(x_{k+1})}.
\end{equation*}
Substituting $\norm{\grad f(x_{k+1})}\geq\frac{f(x_{k+1})-f^\ast}{2D}$ gives
\begin{gather*}
  f(x_k)-f(x_{k+1})
  \geq\min\lb
      \frac{(f(x_{k+1})-f^\ast)^2}
           {48\delta D^2+12\sqrt{6}\sqrt{M_0D^3(f(x_{k+1})-f^\ast)}},
      \frac{r_0(f(x_{k+1})-f^\ast)}{6D}
     \rb\\
  =\frac{f(x_{k+1})-f^\ast}{\Phi(f(x_{k+1})-f^\ast)}.
\end{gather*}
By~\eqref{eq:basic_pd_ms_decrease},
$f(x_{k+1})-f^\ast\leq f(x_k)-f^\ast$.
For every $t\in[f(x_{k+1})-f^\ast,f(x_k)-f^\ast]$,
the monotonicity of $t(1+\Phi(t))$ gives
\begin{equation*}
  \frac{1+\Phi(t)}t
  \geq\frac{(f(x_{k+1})-f^\ast)(1+\Phi(f(x_{k+1})-f^\ast))}{t^2}.
\end{equation*}
Thus, we obtain
\begin{gather*}
  \int_{f(x_{k+1})-f^\ast}^{f(x_k)-f^\ast}\frac{1+\Phi(t)}t dt 
  \geq(f(x_{k+1})-f^\ast)(1+\Phi(f(x_{k+1})-f^\ast))
       \int_{f(x_{k+1})-f^\ast}^{f(x_k)-f^\ast}\frac{dt}{t^2}\\
  =(f(x_{k+1})-f^\ast)(1+\Phi(f(x_{k+1})-f^\ast))
       \ls\frac1{f(x_{k+1})-f^\ast}-\frac1{f(x_k)-f^\ast}\rs\\
  =(1+\Phi(f(x_{k+1})-f^\ast))
        \frac{f(x_k)-f(x_{k+1})}{f(x_k)-f^\ast}
  \geq\frac{1+\Phi(f(x_{k+1})-f^\ast)}{\Phi(f(x_{k+1})-f^\ast)+1}
  =1.
\end{gather*}
Summing over the controlled iterations
proves~\eqref{eq:basic_pd_ms_controlled_count}.
\end{proof}

\subsection{Proof of Theorem~\ref{thm:basic_pd_ms_convergence}}

\begin{theorem}
\label{thm:basic_pd_ms_convergence_appendix}
Suppose that $f$ is convex,~\eqref{as:inexact_hessian} holds, and
either~\eqref{as:gh} or~\eqref{as:lh} holds.
Define
\begin{equation}
  B\eqdef\max\lb1,\frac{2D}{\rho},6\sqrt{M_1}D\rb,
  \qquad
  \e_{\tr}\eqdef\max\lb\frac{M_0D^3}{B^2},\frac{\delta D^2}{B}\rb,
  \label{eq:app_basic_pd_ms_scales}
\end{equation}
where $\rho=+\infty$ under~\eqref{as:gh}.
For every $\e>0$, Algorithm~\ref{alg:basic_pd_ms} generates an iterate $x_N$
satisfying $f(x_N)-f^\ast\leq\e$ within
\begin{equation}
\begin{aligned}
  N\leq 2+\log_2\max\lb1,\frac{2D\norm{\grad f(x_0)}}{\e_{\tr}},
    \frac{\eta_0D^2}{\sqrt6 B\e_{\tr}}\rb
  +100B\log_+\frac{f(x_0)-f^\ast}{\max\{\e,\e_{\tr}\}}
    +80\sqrt{\frac{M_0D^3}{\e}}+64\frac{\delta D^2}{\e}.
\end{aligned}
  \label{eq:basic_pd_ms_explicit_complexity}
\end{equation}
iterations. 

\end{theorem}

\begin{proof}

Each search terminates by Lemma~\ref{lem:pd_ms_backtracking}.
Consider a prefix $x_0,\ldots,x_m$ with
$f(x_m)-f^\ast>\e$.
By~\eqref{eq:basic_pd_ms_localization},
$\norm{\grad f(x_m)}>\e/(2D)$.
Combining~\eqref{eq:basic_pd_ms_tracking_count}
and~\eqref{eq:basic_pd_ms_controlled_count} gives
\begin{equation}
\begin{aligned}
  m&\leq\frac54\int_{f(x_m)-f^\ast}^{f(x_0)-f^\ast}\frac{1+\Phi(t)}t dt
       +\log_2\frac{U_0}{\Theta_{\rm c}(\norm{\grad f(x_m)})}\\
  &\leq\frac54\int_\e^{f(x_0)-f^\ast}\frac{1+\Phi(t)}t dt
       +\log_{2,+}\frac{U_0}{\Theta_{\rm c}(\e/(2D))}.
\end{aligned}
  \label{eq:basic_pd_ms_prefix_count}
\end{equation}
Let $N$ be the first index with
$f(x_N)-f^\ast\leq\e$.
If $N\geq1$, then $f(x_{N-1})-f^\ast>\e$. Therefore,
\begin{equation}
  N\leq1+\frac54\int_\e^{f(x_0)-f^\ast}\frac{1+\Phi(t)}t dt
       +\log_{2,+}\frac{U_0}{\Theta_{\rm c}(\e/(2D))}.
  \label{eq:basic_pd_ms_integral_count}
\end{equation}
We begin with bounding integral term in~\eqref{eq:basic_pd_ms_integral_count}. 

By~\eqref{eq:app_basic_pd_ms_scales} and~\eqref{eq:basic_pd_ms_count_scale},
$1+\Phi(t)\leq79B$ for $t\geq\e_{\tr}$.
For $0<t\leq\e_{\tr}$, we have
\begin{gather*}
  \max\lb\sqrt{\frac{M_0D^3}{t}},\frac{\delta D^2}{t}\rb
  \geq B\geq\max\lb1,\frac{D}{r_0}\rb,
  1+\Phi(t)
  \leq(12\sqrt{6}+1)\sqrt{\frac{M_0D^3}{t}}
       +49\frac{\delta D^2}{t}.
\end{gather*}
If $\e_{\tr}\leq\e<f(x_0)-f^\ast$, the whole integration interval
lies above $\e_{\tr}$, so
\begin{equation*}
  \int_\e^{f(x_0)-f^\ast}\frac{1+\Phi(t)}t dt
  \leq79B\int_\e^{f(x_0)-f^\ast}\frac{dt}t
  =79B\log\frac{f(x_0)-f^\ast}{\e}.
\end{equation*}
If $0<\e<\e_{\tr}\leq f(x_0)-f^\ast$, split the integral as
\begin{gather*}
  \int_\e^{f(x_0)-f^\ast}\frac{1+\Phi(t)}t dt
  =\int_\e^{\e_{\tr}}\frac{1+\Phi(t)}t dt
   +\int_{\e_{\tr}}^{f(x_0)-f^\ast}\frac{1+\Phi(t)}t dt\\
  \leq(12\sqrt{6}+1)\sqrt{M_0D^3}
       \int_\e^{\e_{\tr}}\frac{dt}{t^{3/2}}
    +49\delta D^2\int_\e^{\e_{\tr}}\frac{dt}{t^2}
    +79B\int_{\e_{\tr}}^{f(x_0)-f^\ast}\frac{dt}t\\
  =(24\sqrt{6}+2)\sqrt{M_0D^3}
      \ls\frac1{\sqrt\e}-\frac1{\sqrt{\e_{\tr}}}\rs
    +49\delta D^2\ls\frac1\e-\frac1{\e_{\tr}}\rs
    +79B\log\frac{f(x_0)-f^\ast}{\e_{\tr}}\\
  \leq(24\sqrt{6}+2)\sqrt{\frac{M_0D^3}{\e}}
    +49\frac{\delta D^2}{\e}
    +79B\log\frac{f(x_0)-f^\ast}{\e_{\tr}}.
\end{gather*}
If $0<\e<f(x_0)-f^\ast<\e_{\tr}$, the whole interval lies below
$\e_{\tr}$, so
\begin{gather*}
  \int_\e^{f(x_0)-f^\ast}\frac{1+\Phi(t)}t dt
  \leq(12\sqrt{6}+1)\sqrt{M_0D^3}
       \int_\e^{f(x_0)-f^\ast}\frac{dt}{t^{3/2}}
    +49\delta D^2\int_\e^{f(x_0)-f^\ast}\frac{dt}{t^2}\\
  =(24\sqrt{6}+2)\sqrt{M_0D^3}
      \ls\frac1{\sqrt\e}-\frac1{\sqrt{f(x_0)-f^\ast}}\rs
    +49\delta D^2\ls\frac1\e-\frac1{f(x_0)-f^\ast}\rs\\
  \leq(24\sqrt{6}+2)\sqrt{\frac{M_0D^3}{\e}}
    +49\frac{\delta D^2}{\e}.
\end{gather*}
Here $\log_+\frac{f(x_0)-f^\ast}{\e_{\tr}}=0$ by definition.

Combining the three cases yields
\begin{equation}
  \int_\e^{f(x_0)-f^\ast}\frac{1+\Phi(t)}t dt
  \leq79B\log_+\frac{f(x_0)-f^\ast}{\max\{\e,\e_{\tr}\}}
  \quad+(24\sqrt{6}+2)\sqrt{\frac{M_0D^3}{\e}}
    +49\frac{\delta D^2}{\e}.
  \label{eq:basic_pd_ms_clipped_integral}
\end{equation}
Now, we need to bound the last term in~\eqref{eq:basic_pd_ms_integral_count}.

For $g_2\geq g_1>0$,~\eqref{eq:basic_pd_ms_threshold} gives
$\Theta_{\rm c}(g_2)\leq(g_2/g_1)\Theta_{\rm c}(g_1)$.
Together with monotonicity of $\Theta_{\rm c}$, this yields
\begin{equation*}
  \log_{2,+}\frac{U_0}{\Theta_{\rm c}(\e/(2D))}
  \leq\log_{2,+}\frac{U_0}{\Theta_{\rm c}(\e_{\tr}/(2D))}
    +\log_{2,+}\frac{\e_{\tr}}{\e}.
\end{equation*}
If $\e\geq\e_{\tr}$, then $\log_{2,+}(\e_{\tr}/\e)=0$.
For $0<\e<\e_{\tr}$, we have
\begin{gather*}
  \log_2\frac{\e_{\tr}}{\e}
  \leq3\sqrt{\frac{\e_{\tr}}{\e}}
  =\frac3B\sqrt{\frac{M_0D^3}{\e}},
  \qquad\text{if }\e_{\tr}=\frac{M_0D^3}{B^2},\\
  \log_2\frac{\e_{\tr}}{\e}
  \leq2\frac{\e_{\tr}}{\e}
  =\frac2B\frac{\delta D^2}{\e},
  \qquad\text{if }\e_{\tr}=\frac{\delta D^2}{B}.
\end{gather*}
Since $B\geq1$,
\begin{equation}
  \log_{2,+}\frac{U_0}{\Theta_{\rm c}(\e/(2D))}
  \leq\log_{2,+}\frac{U_0}{\Theta_{\rm c}(\e_{\tr}/(2D))}
    +3\sqrt{\frac{M_0D^3}{\e}}+2\frac{\delta D^2}{\e}.
  \label{eq:basic_pd_ms_clipped_tracking}
\end{equation}
Substituting~\eqref{eq:basic_pd_ms_clipped_integral}
and~\eqref{eq:basic_pd_ms_clipped_tracking}
into~\eqref{eq:basic_pd_ms_integral_count}, and using
\begin{equation*}
  \frac54\cdot79\leq100,
  \qquad
  \frac54(24\sqrt{6}+2)+3\leq80,
  \qquad
  \frac54\cdot49+2\leq64,
\end{equation*}
we obtain
\begin{equation*}
  N\leq 2+\log_{2,+}\frac{U_0}{\Theta_{\rm c}(\e_{\tr}/(2D))}
  +100B\log_+\frac{f(x_0)-f^\ast}{\max\{\e,\e_{\tr}\}}
    +80\sqrt{\frac{M_0D^3}{\e}}+64\frac{\delta D^2}{\e}.
\end{equation*}
Recall that $\Theta_{\rm c}$ is nondecreasing and $\Theta_{\rm c}(ag)\leq a\Theta_{\rm c}(g)$
for $a\geq1$.
Thus,
\begin{gather*}
  \Theta_{\rm c}(\norm{\grad f(x_0)})
  \leq\Theta_{\rm c}\ls
    \max\lb1,\frac{2D\norm{\grad f(x_0)}}{\e_{\tr}}\rb
    \frac{\e_{\tr}}{2D}\rs
  \leq\max\lb1,\frac{2D\norm{\grad f(x_0)}}{\e_{\tr}}\rb
    \Theta_{\rm c}\ls\frac{\e_{\tr}}{2D}\rs.
\end{gather*}
If $\e_{\tr}=M_0D^3/B^2$
\begin{equation*}
  \Theta_{\rm c}\ls\frac{\e_{\tr}}{2D}\rs
  \geq\sqrt{\frac{24M_0\e_{\tr}}D}
  =\frac{2\sqrt6 B\e_{\tr}}{D^2}.
\end{equation*}
Otherwise $\e_{\tr}=\delta D^2/B$
\begin{equation*}
  \Theta_{\rm c}\ls\frac{\e_{\tr}}{2D}\rs
  \geq8\delta=\frac{8B\e_{\tr}}{D^2}
  \geq\frac{2\sqrt6 B\e_{\tr}}{D^2}.
\end{equation*}
Since $U_0=\max\lb2\eta_0,\Theta_{\rm c}(\norm{\grad f(x_0)})\rb$, these estimates give
\begin{equation}
\begin{aligned}
  \log_{2,+}\frac{U_0}{\Theta_{\rm c}(\e_{\tr}/(2D))}
  &=\log_2\max\lb1,
    \frac{\Theta_{\rm c}(\norm{\grad f(x_0)})}{\Theta_{\rm c}(\e_{\tr}/(2D))},
    \frac{2\eta_0}{\Theta_{\rm c}(\e_{\tr}/(2D))}\rb\\
  &\leq\log_2\max\lb1,\frac{2D\norm{\grad f(x_0)}}{\e_{\tr}},
    \frac{\eta_0D^2}{\sqrt6 B\e_{\tr}}\rb.
\end{aligned}
  \label{eq:basic_pd_ms_initial_log_bound}
\end{equation}
Substituting~\eqref{eq:basic_pd_ms_initial_log_bound} into the bound for $N$
proves~\eqref{eq:basic_pd_ms_explicit_complexity}.

\end{proof}

\subsection{Far and Near Convergence Rates}

\begin{corollary}
\label{cor:basic_pd_ms_rates}
Under the assumptions of Theorem~\ref{thm:basic_pd_ms_convergence_appendix},
suppose that $0<D<+\infty$.
Define
\begin{equation*}
  \ell_0\eqdef\max\lb
    1,\frac{2D\norm{\grad f(x_0)}}{\e_{\tr}},
    \frac{\eta_0D^2}{\sqrt6 B\e_{\tr}}
  \rb.
\end{equation*}
For every $k\geq0$ with $f(x_k)-f^\ast>\e_{\tr}$,
\begin{equation}
  f(x_k)-f^\ast\leq\ls f(x_0)-f^\ast\rs
    \exp\ls-\frac{(k-\log_2\ell_0)_+}{100B}\rs,
  \label{eq:basic_pd_ms_far_rate_entry}
\end{equation}
where we write $(t)_+\eqdef\max\lb0,t\rb$.
If $f(x_k)-f^\ast\leq\e_{\tr}$, the subsequent iterates satisfy
\begin{equation}
  f(x_{k+n})-f^\ast
  \leq\frac{7396M_0D^3}{\bigl(B+(n-\log_2\ell_0)_+\bigr)^2}
    +\frac{140\delta D^2}{B+(n-\log_2\ell_0)_+},
  \qquad n=0,1,\ldots,
  \label{eq:basic_pd_ms_near_rate_entry}
\end{equation}
where $t_+\eqdef\max\lb t,0\rb$.
\end{corollary}

\begin{proof}
If $f(x_k)-f^\ast>\e_{\tr}$,~\eqref{eq:basic_pd_ms_prefix_count} gives
\begin{equation*}
  k\leq\frac54\int_{f(x_k)-f^\ast}^{f(x_0)-f^\ast}
    \frac{1+\Phi(t)}t dt
    +\log_2\frac{U_0}{\Theta_{\rm c}(\norm{\grad f(x_k)})}.
\end{equation*}
For $t\geq\e_{\tr}$, definitions~\eqref{eq:app_basic_pd_ms_scales}
and~\eqref{eq:basic_pd_ms_count_scale} give
\begin{gather*}
  \Phi(t)\leq\max\lb
    6B,12\sqrt6\sqrt{\frac{M_0D^3}{\e_{\tr}}}
      +48\frac{\delta D^2}{\e_{\tr}}
  \rb
  \leq\max\lb6B,(12\sqrt6+48)B\rb
  =(12\sqrt6+48)B.
\end{gather*}
Therefore, $1+\Phi(t)\leq79B$.
By~\eqref{eq:basic_pd_ms_localization}, monotonicity of $\Theta_{\rm c}$,
and~\eqref{eq:basic_pd_ms_initial_log_bound},
\begin{equation*}
  \log_2\frac{U_0}{\Theta_{\rm c}(\norm{\grad f(x_k)})}
  \leq\log_{2,+}\frac{U_0}{\Theta_{\rm c}(\e_{\tr}/(2D))}
  \leq\log_2\ell_0.
\end{equation*}
Substituting these bounds yields
\begin{gather*}
  k\leq\frac54\cdot79B
    \int_{f(x_k)-f^\ast}^{f(x_0)-f^\ast}\frac{dt}t+\log_2\ell_0\\
  =\frac{395}{4}B\log\frac{f(x_0)-f^\ast}{f(x_k)-f^\ast}+\log_2\ell_0
  \leq100B\log\frac{f(x_0)-f^\ast}{f(x_k)-f^\ast}+\log_2\ell_0.
\end{gather*}
Rearranging and using $f(x_k)\leq f(x_0)$ when $k\leq\log_2\ell_0$
proves~\eqref{eq:basic_pd_ms_far_rate_entry}.

Now suppose that $f(x_k)-f^\ast\leq\e_{\tr}$.
Let $n_{\tr}\leq k$ be the first index with $f(x_{n_{\tr}})-f^\ast\leq\e_{\tr}$.
For $n\geq0$ with $f(x_{n_{\tr}+n})>f^\ast$,
the first inequality in~\eqref{eq:basic_pd_ms_prefix_count},
applied from iteration $n_{\tr}$, gives
\begin{equation}
  n\leq\frac54\int_{f(x_{n_{\tr}+n})-f^\ast}^{f(x_{n_{\tr}})-f^\ast}
    \frac{1+\Phi(t)}t dt
    +\log_2\frac{U_{n_{\tr}}}{\Theta_{\rm c}(\norm{\grad f(x_{n_{\tr}+n})})}.
  \label{eq:basic_pd_ms_segment_count}
\end{equation}
For $0<t\leq\e_{\tr}$, the definitions of $\e_{\tr}$ and $B$ give
\begin{equation*}
  \max\lb\sqrt{\frac{M_0D^3}{t}},\frac{\delta D^2}{t}\rb
  \geq B=\max\lb1,\frac{D}{r_0}\rb.
\end{equation*}
Hence, by the definition of $\Phi$ in~\eqref{eq:basic_pd_ms_count_scale},
\begin{gather*}
  1+\Phi(t)\leq1+12\sqrt6\sqrt{\frac{M_0D^3}{t}}
    +48\frac{\delta D^2}{t}
  \leq(12\sqrt6+1)\sqrt{\frac{M_0D^3}{t}}
    +49\frac{\delta D^2}{t}.
\end{gather*}
Since $0<f(x_{n_{\tr}+n})-f^\ast\leq f(x_{n_{\tr}})-f^\ast\leq\e_{\tr}$,
\begin{gather*}
  \int_{f(x_{n_{\tr}+n})-f^\ast}^{f(x_{n_{\tr}})-f^\ast}\frac{1+\Phi(t)}t dt
  \leq(12\sqrt6+1)\sqrt{M_0D^3}
    \int_{f(x_{n_{\tr}+n})-f^\ast}^{f(x_{n_{\tr}})-f^\ast}\frac{dt}{t^{3/2}}
    +49\delta D^2\int_{f(x_{n_{\tr}+n})-f^\ast}^{f(x_{n_{\tr}})-f^\ast}\frac{dt}{t^2}\\
  =(24\sqrt6+2)\sqrt{M_0D^3}
    \ls\frac1{\sqrt{f(x_{n_{\tr}+n})-f^\ast}}-\frac1{\sqrt{f(x_{n_{\tr}})-f^\ast}}\rs\\
  +49\delta D^2\ls\frac1{f(x_{n_{\tr}+n})-f^\ast}-\frac1{f(x_{n_{\tr}})-f^\ast}\rs\\
  \leq(24\sqrt6+2)\sqrt{\frac{M_0D^3}{f(x_{n_{\tr}+n})-f^\ast}}
    +49\frac{\delta D^2}{f(x_{n_{\tr}+n})-f^\ast}.
\end{gather*}
To bound the logarithm in~\eqref{eq:basic_pd_ms_segment_count}, first express
$U_{n_{\tr}}$ through $\Theta_{\rm c}$.
If there is no controlled iteration before $n_{\tr}$, the proof of
Lemma~\ref{lem:basic_pd_ms_tracking} gives $U_{n_{\tr}}\leq U_0$.
Otherwise, let $i<n_{\tr}$ be the last controlled iteration. Then
\begin{equation*}
  U_{n_{\tr}}\leq U_{i+1}=\Theta_{\rm c}(\norm{\grad f(x_{i+1})}).
\end{equation*}
Next, by 
$f(x_0)-f^\ast>\e_{\tr}
=\max\lb M_0D^3/B^2,\delta D^2/B\rb$
\begin{equation*}
  \delta\leq\frac{B(f(x_0)-f^\ast)}{D^2},
  \qquad
  M_0\leq\frac{B^2(f(x_0)-f^\ast)}{D^3}.
\end{equation*}
By the definition of $B$,
\begin{equation*}
  \max\lb\frac4\rho,12\sqrt{M_1}\rb
  =\frac2D\max\lb\frac{2D}{\rho},6\sqrt{M_1}D\rb
  \leq\frac{2B}{D}.
\end{equation*}
If $\norm{\grad f(x_{i+1})}<(f(x_0)-f^\ast)/D$, then $B\geq1$ gives $\norm{\grad f(x_{i+1})}<B\frac{f(x_0)-f^\ast}{D}$.

Otherwise, since iteration $i$ is controlled, by~\eqref{eq:basic_pd_ms_decrease},~\eqref{eq:basic_pd_ms_controlled},
and the definition of $\Theta_{\rm c}$
\begin{gather*}
  \norm{\grad f(x_{i+1})}
  \leq\frac32(f(x_0)-f^\ast)
  \max\lb
    \frac{8\delta}{\norm{\grad f(x_{i+1})}}
      +\sqrt{\frac{48M_0}{\norm{\grad f(x_{i+1})}}},
    \frac4\rho,12\sqrt{M_1}
  \rb\\
  \leq\frac{3(8+\sqrt{48})B(f(x_0)-f^\ast)}{2D}
  \leq\frac{23B(f(x_0)-f^\ast)}{D}.
\end{gather*}
Thus, 
\begin{equation*}
  \norm{\grad f(x_{i+1})}
  \leq23B\frac{f(x_0)-f^\ast}{D}
  \leq23B\norm{\grad f(x_0)}.
\end{equation*}
Substituting into $\Theta_{\rm c}$ and using its monotonicity and
$\Theta_{\rm c}(ag)\leq a\Theta_{\rm c}(g)$ for $a\geq1$, we obtain
\begin{gather*}
  U_{n_{\tr}}\leq\Theta_{\rm c}(\norm{\grad f(x_{i+1})})
  \leq\Theta_{\rm c}(23B\norm{\grad f(x_0)}) 
  \leq23B\Theta_{\rm c}(\norm{\grad f(x_0)})
  \leq23B U_0.
\end{gather*}
Using this bound, $\norm{\grad f(x_{n_{\tr}+n})}\geq(f(x_{n_{\tr}+n})-f^\ast)/(2D)$,
and~\eqref{eq:basic_pd_ms_clipped_tracking}
with~\eqref{eq:basic_pd_ms_initial_log_bound}, we obtain
\begin{gather*}
  \log_2\frac{U_{n_{\tr}}}{\Theta_{\rm c}(\norm{\grad f(x_{n_{\tr}+n})})}
  \leq\log_2(23B)
    +\log_{2,+}\frac{U_0}{\Theta_{\rm c}((f(x_{n_{\tr}+n})-f^\ast)/(2D))}\\
  \leq\log_2\ell_0+\log_2(23B)
    +3\sqrt{\frac{M_0D^3}{f(x_{n_{\tr}+n})-f^\ast}}
    +2\frac{\delta D^2}{f(x_{n_{\tr}+n})-f^\ast}.
\end{gather*}
Substitution into~\eqref{eq:basic_pd_ms_segment_count} yields
\begin{equation}
  n\leq\log_2\ell_0+\log_2(23B)
    +80\sqrt{\frac{M_0D^3}{f(x_{n_{\tr}+n})-f^\ast}}
    +64\frac{\delta D^2}{f(x_{n_{\tr}+n})-f^\ast}.
  \label{eq:basic_pd_ms_near_segment_count}
\end{equation}
The definition of $\e_{\tr}$ and $f(x_{n_{\tr}+n})-f^\ast\leq\e_{\tr}$ imply
\begin{equation*}
  B\leq\max\lb
    \sqrt{\frac{M_0D^3}{f(x_{n_{\tr}+n})-f^\ast}},
    \frac{\delta D^2}{f(x_{n_{\tr}+n})-f^\ast}\rb,
  \qquad \log_2(23B)\leq5B.
\end{equation*}
Inequality~\eqref{eq:basic_pd_ms_near_segment_count} therefore gives
\begin{gather*}
  B+(n-\log_2\ell_0)_+
  \leq6B+80\sqrt{\frac{M_0D^3}{f(x_{n_{\tr}+n})-f^\ast}}
    +64\frac{\delta D^2}{f(x_{n_{\tr}+n})-f^\ast}\\
  \leq86\sqrt{\frac{M_0D^3}{f(x_{n_{\tr}+n})-f^\ast}}
    +70\frac{\delta D^2}{f(x_{n_{\tr}+n})-f^\ast}.
\end{gather*}
Multiplying by $f(x_{n_{\tr}+n})-f^\ast$ and using $2ab\leq a^2+b^2$, we obtain
\begin{gather*}
  \ls B+(n-\log_2\ell_0)_+\rs\ls f(x_{n_{\tr}+n})-f^\ast\rs
  \leq86\sqrt{M_0D^3\ls f(x_{n_{\tr}+n})-f^\ast\rs}+70\delta D^2\\
  \leq\frac{\ls B+(n-\log_2\ell_0)_+\rs\ls f(x_{n_{\tr}+n})-f^\ast\rs}{2}
    +\frac{3698M_0D^3}{B+(n-\log_2\ell_0)_+}+70\delta D^2,\\
  f(x_{n_{\tr}+n})-f^\ast
  \leq\frac{7396M_0D^3}{\ls B+(n-\log_2\ell_0)_+\rs^2}
    +\frac{140\delta D^2}{B+(n-\log_2\ell_0)_+}.
\end{gather*}
Since $f(x_{k+n})\leq f(x_{n_{\tr}+n})$,
this proves~\eqref{eq:basic_pd_ms_near_rate_entry}.
\end{proof}

\subsection{Oracle complexity}
\label{sec:pd_ms_oracle}

\begin{corollary}
\label{cor:basic_pd_ms_oracle_complexity}
Let $N\eqdef\min\lb k\geq0:f(x_k)-f^\ast\leq\e\rb$ and $g_0 = \|\nabla f(x_0)\|$.
The total number of rejected trials of Algorithm~\ref{alg:basic_pd_ms} before obtaining $x_N$,
including failed $\PD$ tests, is at most
\begin{equation}
  \frac54(N-1)+4
    +\log_{2,+}\frac{Bg_0}{D\eta_0}
    +\log_{2,+}\frac{\e_{\tr}}{\e}.
  \label{eq:basic_pd_ms_total_backtracking}
\end{equation}
Thus, the total number of gradient evaluations and rejected trials is
bounded by
\begin{equation}
  \cO\bigg(1+\log\ell_0
    +\log_+\frac{Bg_0}{D\eta_0}
    +B\log_+\frac{f(x_0)-f^\ast}{\max\lb\e,\e_{\tr}\rb}
    +\sqrt{\frac{M_0D^3}{\e}}+\frac{\delta D^2}{\e}\bigg).
  \label{eq:basic_pd_ms_total_complexity}
\end{equation}
\end{corollary}

\begin{proof}
  Let $b_j$ be the number of rejected trials on iteration $j$. Then, $b_0 = \log_2 \tfrac{\tau_0}{\eta_0}$ and $\tau_j = 2^{b_j}\eta_j$. Let $g_j = \|\nabla f(x_j)\|$ Therefore, for $j \geq 1$ by line~\ref{line:basic_pd_ms_eta_update} of Algorithm~\ref{alg:basic_pd_ms}
  \begin{equation*}
    b_j = \log_2 \frac{\tau_j}{\eta_j} = \log_2 \frac{2 \tau_j}{\tau_{j-1} \min \lb 1, g_j/g_{j-1}\rb} = 1 + \log_2 \frac{\tau_j}{\tau_{j-1}} + \log_2 \max\lb 1, \frac{g_{j-1}}{g_{j}}\rb.
  \end{equation*}
  Next, let us sum the number of rejected trials to obtain $x_N$
  \begin{equation*}
    \sum_{j=0}^{N-1} b_j = \log_2 \frac{\tau_0}{\eta_0} + (N - 1) + \log_2 \frac{\tau_{N-1}}{\tau_{0}} + \sum_{j=1}^{N-1} \log_{2, +} \frac{g_{j-1}}{g_{j}}.
  \end{equation*}
  Let us bound the last sum 
  \begin{equation*}
    \sum_{j=1}^{N-1} \log_{2, +} \frac{g_{j-1}}{g_{j}} = \sum_{j=1}^{N-1} \ls  \log_{2, +} \frac{g_{j}}{g_{j-1}} + \log_2 \frac{g_{j-1}}{g_j} \rs  = \sum_{j=1}^{N-1} \log_{2, +} \frac{g_{j}}{g_{j-1}} + \log_2 \frac{g_0}{g_{N-1}}
  \end{equation*}
  Therefore, 
  \begin{equation}
    \label{eq:basic_pd_ms_rejections_boundary}
    \sum_{j=0}^{N-1} b_j = N - 1 + \log_2 \frac{g_0 \tau_{N-1}}{\eta_0 g_{N-1}} + \sum_{j=1}^{N-1} \log_{2, +} \frac{g_{j}}{g_{j-1}}.
  \end{equation}
  By~\eqref{eq:basic_pd_ms_gradient_growth}, we have $g_j \leq (2/\sqrt{3}) g_{j-1}$. Therefore,
  \begin{equation}
    \label{eq:basic_pd_ms_gradient_log_sum}
     \sum_{j=1}^{N-1} \log_{2, +} \frac{g_{j}}{g_{j-1}} \leq (N-1) \log_{2, +} \frac{2}{\sqrt{3}} \leq \frac{N-1}{4}.
  \end{equation}
  Now, it remains to bound $\log_{2, + } (\tau_{N-1}/g_{N-1})$. By Lemma~\ref{lem:pd_ms_backtracking}, every accepted $\tau$ satisfies 
  \begin{equation*}
    \tau_j \leq \max \lb \eta_j, 8 \delta, 4 \sqrt{L(x_j)g_j}, \frac{4 g_j}{\rho}\rb.
  \end{equation*}
  Therefore,
  \begin{equation}
    \label{eq:app_tau_g_bound}
    \frac{\tau_j}{g_j} \leq \max \lb \frac{\eta_j}{g_j}, \frac{8 \delta}{g_j}, 4 \sqrt{\frac{M_0}{g_j} + M_1}, \frac{4}{\rho}\rb.
  \end{equation}
By the definitions~\eqref{eq:app_basic_pd_ms_scales} of $B$ and $\e_{\tr}$
\begin{equation*}
  \delta\leq\frac{B\e_{\tr}}{D^2},
  \qquad M_0\leq\frac{B^2\e_{\tr}}{D^3},
  \qquad M_1\leq\frac{B^2}{36D^2},
  \qquad \frac1\rho\leq\frac{B}{2D}.
\end{equation*}
Next, by definition, $N$ is the first iteration, such that $f(x_N) - f^\ast \leq \e$. Then for $j<N$ we have $g_j \geq \e / (2D)$. Therefore, 
\begin{gather*}
  \frac{8\delta}{g_j}
  \leq\frac{16\delta D}{\e}
  \leq\frac{16B}{D}\frac{\e_{\tr}}{\e},
\qquad   \frac4\rho\leq\frac{2B}{D}\leq\frac{16B}{D}\\
  4\sqrt{\frac{M_0}{g_j}+M_1}
  \leq4\sqrt{\frac{2M_0D}{\e}+M_1}
  \leq\frac{4B}{D}\sqrt{2\frac{\e_{\tr}}{\e}+\frac1{36}}
  \leq\frac{16B}{D}\max\lb1,\frac{\e_{\tr}}{\e}\rb,
\end{gather*}
Moreover, for $j\geq1$, the update for $\eta_j$ gives
\begin{equation*}
  \frac{\eta_j}{g_j}
  =\frac{\tau_{j-1}}{2g_{j-1}}
    \min\lb1,\frac{g_{j-1}}{g_j}\rb
  \leq\frac{\tau_{j-1}}{2g_{j-1}}.
\end{equation*}
Substituting these estimates
into~\eqref{eq:app_tau_g_bound} gives for $1\leq j<N$
\begin{equation*}
  \frac{\tau_0}{g_0}
  \leq\max\lb\frac{\eta_0}{g_0},\frac{16B}{D}\max\lb1,\frac{\e_{\tr}}{\e}\rb\rb,
  \quad
  \frac{\tau_j}{g_j}
  \leq\max\lb\frac{\tau_{j-1}}{2g_{j-1}},\frac{16B}{D}\max\lb1,\frac{\e_{\tr}}{\e}\rb\rb. 
\end{equation*}
Induction therefore gives
\begin{equation}
  \label{eq:basic_pd_ms_oracle_normalized_bound}
  \frac{\tau_j}{g_j}
  \leq\max\lb\frac{\eta_0}{g_0},\frac{16B}{D}\max\lb1,\frac{\e_{\tr}}{\e}\rb\rb,
  \qquad 0\leq j<N.
\end{equation}
Finally, by combining~\eqref{eq:basic_pd_ms_rejections_boundary},~\eqref{eq:basic_pd_ms_gradient_log_sum}, and~\eqref{eq:basic_pd_ms_oracle_normalized_bound}, we obtain
\begin{gather*}
   \sum_{j=0}^{N-1} b_j \leq\frac54(N-1)
    +\log_2\max\lb1,\frac{16Bg_0}{D\eta_0}\max\lb1,\frac{\e_{\tr}}{\e}\rb\rb\\
  \leq\frac54(N-1)+4
    +\log_{2,+}\frac{Bg_0}{D\eta_0}
    +\log_{2,+}\frac{\e_{\tr}}\e.
\end{gather*}
This proves~\eqref{eq:basic_pd_ms_total_backtracking}.

Each trial requires at most one
new gradient evaluation. Hence, the total number of gradient
evaluations is at most $1+N+\sum_{j=0}^{N-1}b_j$. Next, by the definition of $\e_{\tr}$, we have $\max\lb
    \sqrt{\frac{M_0D^3}{\e_{\tr}}},
    \frac{\delta D^2}{\e_{\tr}}
\rb=B\geq1$. Using $\log_2 t\leq\min\lb2\sqrt t,t\rb$ for $t\geq1$, we obtain $\log_{2,+}\frac{\e_{\tr}}\e \leq2\sqrt{\frac{M_0D^3}{\e}}+\frac{\delta D^2}{\e}$. Combining these bounds with~\eqref{eq:basic_pd_ms_complexity}
proves~\eqref{eq:basic_pd_ms_total_complexity}.

\end{proof}

\subsection{Strongly Convex Analysis}
\label{app:basic_pd_ms_strong}

Let $f$ be strongly convex
$\mu$-strongly convex for $\mu>0$:
\begin{equation}
  f(y)\geq f(x)+\la\grad f(x),y-x\ra
    +\frac\mu2\norm{y-x}^2,
  \qquad x,y\in\R^d.
  \label{as:basic_pd_ms_strong}
\end{equation}
This implies
\begin{equation}
  \frac\mu2\norm{x-x^\ast}^2
  \leq f(x)-f^\ast
  \leq\frac{\norm{\grad f(x)}^2}{2\mu}.
  \label{eq:basic_pd_ms_strong_gap}
\end{equation}
In particular, strong convexity with monotonicity (Lemma~\ref{lem:basic_pd_ms_step}) imply
\begin{equation}
  \norm{x_k-x^\ast}
  \leq\sqrt{\frac{2(f(x_k)-f^\ast)}\mu}
  \leq\sqrt{\frac{2(f(x_0)-f^\ast)}\mu}.
  \label{eq:basic_pd_ms_strong_localization}
\end{equation}
We also denote $\delta_k \eqdef \|\nabla^2 f(x_k) - H(x_k)\| \leq \delta$ for all $h \geq 0$.

\subsubsection{One Step Estimates}

\begin{lemma}
\label{lem:basic_pd_ms_strong_step}
Under~\eqref{as:basic_pd_ms_strong}, every completed iteration of
Algorithm~\ref{alg:basic_pd_ms} satisfies
\begin{gather}
  f(x_{k+1})-f^\ast
  \leq\ls\frac{3\tau_k}{3\tau_k+2\mu}\rs^2
       (f(x_k)-f^\ast),
  \label{eq:basic_pd_ms_strong_step}\\
  \norm{s_k}\leq
  \frac{\norm{\grad f(x_k)}}{\mu+\tau_k/2}
  \leq\frac{\norm{\grad f(x_k)}}\mu.
  \label{eq:basic_pd_ms_strong_length}
\end{gather}
\end{lemma}

\begin{proof}
By~\eqref{eq:basic_pd_ms_squared_test},
\begin{equation*}
  -\la\grad f(x_{k+1}),s_k\ra
  \geq\frac{\norm{\grad f(x_{k+1})}^2}{2\tau_k}
    +\frac{3\tau_k}{8}\norm{s_k}^2.
\end{equation*}
Apply~\eqref{as:basic_pd_ms_strong} at $x_{k+1}$, then
use~\eqref{eq:basic_pd_ms_output_gradient}
and~\eqref{eq:basic_pd_ms_strong_gap}:
\begin{gather*}
  f(x_k)-f(x_{k+1})
  \geq-\la\grad f(x_{k+1}),s_k\ra+\frac\mu2\norm{s_k}^2
  \geq\frac{\norm{\grad f(x_{k+1})}^2}{2\tau_k}
    +\ls\frac{3\tau_k}{8}+\frac\mu2\rs\norm{s_k}^2\\
  \geq\ls\frac1{2\tau_k}
      +\frac4{9\tau_k^2}\ls\frac{3\tau_k}{8}+\frac\mu2\rs\rs
      \norm{\grad f(x_{k+1})}^2
  =\ls\frac2{3\tau_k}+\frac{2\mu}{9\tau_k^2}\rs
      \norm{\grad f(x_{k+1})}^2\\
  \geq\ls\frac{4\mu}{3\tau_k}+\frac{4\mu^2}{9\tau_k^2}\rs
      (f(x_{k+1})-f^\ast).
\end{gather*}
The last inequality uses~\eqref{eq:basic_pd_ms_strong_gap} at $x_{k+1}$:
\begin{equation*}
  \norm{\grad f(x_{k+1})}^2
  \geq2\mu(f(x_{k+1})-f^\ast).
\end{equation*}
Adding $f(x_{k+1})-f^\ast$ to the bound for $f(x_k)-f(x_{k+1})$
and completing the square gives
\begin{gather*}
  f(x_k)-f^\ast
  \geq\ls1+\frac{4\mu}{3\tau_k}+\frac{4\mu^2}{9\tau_k^2}\rs
      (f(x_{k+1})-f^\ast)\\
  =\ls1+\frac{2\mu}{3\tau_k}\rs^2
      (f(x_{k+1})-f^\ast),
\end{gather*}
which is~\eqref{eq:basic_pd_ms_strong_step}.
For~\eqref{eq:basic_pd_ms_strong_length}, recall that $s_k=x_{k+1}-x_k$.
Strong convexity gives
\begin{equation*}
  \mu\norm{s_k}^2
  \leq\la\grad f(x_{k+1})-\grad f(x_k),s_k\ra.
\end{equation*}
The accepted $\MS$ test~\eqref{eq:pd_ms_tests} gives
\begin{equation*}
  \norm{\grad f(x_{k+1})+\tau_k s_k}
  \leq\frac{\tau_k}{2}\norm{s_k}.
\end{equation*}
Hence, by Cauchy-Schwarz,
\begin{gather*}
  \la\grad f(x_{k+1}),s_k\ra
  =\la\grad f(x_{k+1})+\tau_k s_k,s_k\ra-\tau_k\norm{s_k}^2\\
  \leq\norm{\grad f(x_{k+1})+\tau_k s_k}\norm{s_k}
    -\tau_k\norm{s_k}^2
  \leq\frac{\tau_k}{2}\norm{s_k}^2-\tau_k\norm{s_k}^2
  =-\frac{\tau_k}{2}\norm{s_k}^2.
\end{gather*}
Combining the two bounds and applying Cauchy-Schwarz once more gives
\begin{gather*}
  \mu\norm{s_k}^2
  \leq-\frac{\tau_k}{2}\norm{s_k}^2-\la\grad f(x_k),s_k\ra,\\
  \ls\mu+\frac{\tau_k}{2}\rs\norm{s_k}^2
  \leq-\la\grad f(x_k),s_k\ra
  \leq\norm{\grad f(x_k)}\norm{s_k}.
\end{gather*}
Dividing by $\norm{s_k}>0$ and using $\tau_k>0$ gives
\begin{equation*}
  \norm{s_k}\leq\frac{\norm{\grad f(x_k)}}{\mu+\tau_k/2}
  \leq\frac{\norm{\grad f(x_k)}}\mu.
\end{equation*}
\end{proof}

\subsubsection{Global Linear Convergence with Inexact Hessians}

We recall definition of $r_0$ from~\eqref{eq:ahpe_damped_movement_radius}:
\begin{equation*}
  r_0=
  \begin{cases}
    \frac1{6\sqrt{M_1}},
    & \text{under~\eqref{as:gh}},\\
    \min\lb\frac\rho2,\frac1{6\sqrt{M_1}}\rb,
    & \text{under~\eqref{as:lh}}.
  \end{cases}
\end{equation*}

\begin{theorem}
\label{thm:basic_pd_ms_strong_global}
Suppose that~\eqref{as:basic_pd_ms_strong}
and~\eqref{as:inexact_hessian} hold, and either~\eqref{as:gh}
or~\eqref{as:lh} holds. Define
\begin{equation}
  \overline\tau\eqdef\max\lb
    \eta_0,\ 8\delta+48M_0\sqrt{\frac{2(f(x_0)-f^\ast)}\mu},
    \frac{2(f(x_0)-f^\ast)}{r_0^2}\rb.
  \label{eq:basic_pd_ms_strong_scale}
\end{equation}
Every completed iteration satisfies $\tau_k\leq\overline\tau$, and
\begin{equation}
  f(x_k)-f^\ast
  \leq\ls\frac{3\overline\tau}{3\overline\tau+2\mu}\rs^{2k}
       (f(x_0)-f^\ast),\qquad k\geq0.
  \label{eq:basic_pd_ms_strong_global_rate}
\end{equation}
\end{theorem}

\begin{proof}
All searches terminate by Lemma~\ref{lem:pd_ms_backtracking}.
We prove $\tau_k\leq\overline\tau$ by induction.

Suppose that on iteration $k$ at least one trial is rejected, so $\tau_k>\eta_k$. If $\norm{s_k}\leq r_0$, then~\eqref{eq:ahpe_damped_movement}
and~\eqref{eq:basic_pd_ms_strong_localization} give
\begin{gather*}
  \tau_k\leq8\delta_k+24M_0\norm{s_k}
  \leq8\delta+24M_0\ls\norm{x_k-x^\ast}
                         +\norm{x_{k+1}-x^\ast}\rs\\
  \leq8\delta+48M_0\sqrt{\frac{2(f(x_0)-f^\ast)}\mu}
  \leq\overline\tau.
\end{gather*}
If $\norm{s_k}>r_0$, then~\eqref{eq:basic_pd_ms_decrease} gives
\begin{equation*}
  \frac{\tau_k}{2}r_0^2
  <\frac{\tau_k}{2}\norm{s_k}^2
  \leq f(x_k)-f(x_{k+1})\leq f(x_0)-f^\ast,
\end{equation*}
so again $\tau_k\leq\overline\tau$. 

Now, let's assume that in iteration $k$ there is no rejection.
For $k=0$, the accepted value is $\tau_0=\eta_0\leq\overline\tau$.
For $k\geq1$, line~\ref{line:basic_pd_ms_eta_update} gives
\begin{equation*}
  \tau_k=\eta_k\leq\frac{\tau_{k-1}}2\leq\overline\tau.
\end{equation*}
Thus the induction closes. Since $0<\tau_k\leq\overline\tau$,
\begin{gather*}
  \ls\frac{3\tau_k}{3\tau_k+2\mu}\rs^2
  =\ls1-\frac{2\mu}{3\tau_k+2\mu}\rs^2 
  \leq\ls1-\frac{2\mu}{3\overline\tau+2\mu}\rs^2
  =\ls\frac{3\overline\tau}{3\overline\tau+2\mu}\rs^2.
\end{gather*}
Substitution into~\eqref{eq:basic_pd_ms_strong_step} and iteration
prove~\eqref{eq:basic_pd_ms_strong_global_rate}.
\end{proof}

\begin{corollary}
\label{cor:basic_pd_ms_strong_gradient_complexity}
Under the assumptions of Theorem~\ref{thm:basic_pd_ms_strong_global},
let $0<\e<f(x_0)-f^\ast$ and
$N\eqdef\min\lb k\geq0:f(x_k)-f^\ast\leq\e\rb$.
The total number of gradient evaluations before obtaining $x_N$ is at most
\begin{equation}
  1+\frac94N
  +\log_{2,+}\ls\frac{\norm{\grad f(x_0)}}{\eta_0}
    \max\lb\frac{24M_0}{\mu}+\frac{8\delta}{\sqrt{2\mu\e}},
      \frac2{r_0}\rb\rs.
  \label{eq:basic_pd_ms_strong_gradient_cost}
\end{equation}
Consequently, it is bounded by
\begin{equation*}
  1+\frac94\left\lceil
    \frac{\log((f(x_0)-f^\ast)/\e)}
         {2\log\ls1+\frac{2\mu}{3\overline\tau}\rs}
  \right\rceil
  +\log_{2,+}\ls\frac{\norm{\grad f(x_0)}}{\eta_0}
    \max\lb\frac{24M_0}{\mu}+\frac{8\delta}{\sqrt{2\mu\e}},
      \frac2{r_0}\rb\rs.
\end{equation*}
For $\delta=0$, the logarithmic term in~\eqref{eq:basic_pd_ms_strong_gradient_cost}
is independent of $\e$, and the bound holds for any $N\geq1$ completed
iterations.
\end{corollary}

\begin{proof}
Let $g_k\eqdef\norm{\grad f(x_k)}$ and let $b_k$ count the rejected trials
on iteration $k$.
By~\eqref{eq:basic_pd_ms_rejections_boundary} and~\eqref{eq:basic_pd_ms_gradient_log_sum},
\begin{equation*}
  \sum_{k=0}^{N-1}b_k
  \leq\frac54(N-1)+\log_2\frac{g_0\tau_{N-1}}{\eta_0g_{N-1}}.
\end{equation*}
We prove by induction that
\begin{equation}
  \frac{\tau_k}{g_k}
  \leq\max\lb\frac{2^{-k}\eta_0}{g_0},
    \frac{24M_0}{\mu}+\frac{8\delta}{g_k},
    \frac2{r_0}\rb.
  \label{eq:basic_pd_ms_strong_gradient_ratio}
\end{equation}
Suppose first that at least one trial is rejected.
If $\norm{s_k}\leq r_0$, then~\eqref{eq:ahpe_damped_movement}
and~\eqref{eq:basic_pd_ms_strong_length} give
\begin{equation*}
  \tau_k\leq8\delta+24M_0\norm{s_k}
  \leq8\delta+\frac{24M_0}{\mu}g_k.
\end{equation*}
If $\norm{s_k}>r_0$, then~\eqref{eq:basic_pd_ms_strong_length} gives
\begin{equation*}
  \tau_k\leq\frac{2g_k}{\norm{s_k}}
  <\frac{2g_k}{r_0}.
\end{equation*}
Both estimates imply~\eqref{eq:basic_pd_ms_strong_gradient_ratio}
for every $k\geq0$ with a rejection.

For $k=0$, the estimates above apply if a trial is rejected.
Otherwise, $\tau_0=\eta_0$, and
\begin{equation*}
  \frac{\tau_0}{g_0}=\frac{\eta_0}{g_0}
  \leq\max\lb\frac{\eta_0}{g_0},
    \frac{24M_0}{\mu}+\frac{8\delta}{g_0},\frac2{r_0}\rb.
\end{equation*}
Thus~\eqref{eq:basic_pd_ms_strong_gradient_ratio} holds for $k=0$.
Now let $k\geq1$ and assume it holds for $k-1$.
If a trial is rejected, the estimates above apply again. If there is no
rejection, line~\ref{line:basic_pd_ms_eta_update} gives
\begin{gather*}
  \tau_k=\frac{\tau_{k-1}}2\min\lb1,\frac{g_k}{g_{k-1}}\rb
  \leq\max\lb2^{-k}\eta_0\frac{g_k}{g_0},
    4\delta+\frac{12M_0}{\mu}g_k,
    \frac{g_k}{r_0}\rb\\
  \leq\max\lb2^{-k}\eta_0\frac{g_k}{g_0},
    8\delta+\frac{24M_0}{\mu}g_k,
    \frac{2g_k}{r_0}\rb.
\end{gather*}
This proves~\eqref{eq:basic_pd_ms_strong_gradient_ratio} by induction.
By the definition of $N$, we have $f(x_{N-1})-f^\ast>\e$.
By~\eqref{eq:basic_pd_ms_strong_gap}, $g_{N-1}>\sqrt{2\mu\e}$.
Substituting this into~\eqref{eq:basic_pd_ms_strong_gradient_ratio}, we obtain
\begin{gather*}
  \frac{g_0\tau_{N-1}}{\eta_0g_{N-1}}
  \leq\max\lb2^{-(N-1)},
    \frac{g_0}{\eta_0}\ls\frac{24M_0}{\mu}+\frac{8\delta}{g_{N-1}}\rs,
    \frac{2g_0}{\eta_0r_0}\rb\\
  \leq\max\lb1,
    \frac{g_0}{\eta_0}\ls\frac{24M_0}{\mu}+\frac{8\delta}{\sqrt{2\mu\e}}\rs,
    \frac{2g_0}{\eta_0r_0}\rb.
\end{gather*}
Each trial passing $\PD$ requires at most one new gradient evaluation.
The accepted gradient is reused at the next iteration. Hence the total
number of gradient evaluations is at most
\begin{equation*}
  1+N+\sum_{k=0}^{N-1}b_k,
\end{equation*}
which, together with the preceding bounds, proves~\eqref{eq:basic_pd_ms_strong_gradient_cost}.

\end{proof}

\subsubsection{Proof of Theorem~\ref{thm:basic_pd_ms_strong_main}}

We refine the global bound from Theorem~\ref{thm:basic_pd_ms_strong_global} to show superlinear convergence by retaining the individual Hessian errors
$\delta_k$ instead of replacing them by $\delta$.

\begin{theorem}
\label{thm:basic_pd_ms_refined_global}
Under the assumptions of Theorem~\ref{thm:basic_pd_ms_strong_global},
let $\overline\tau$ be defined by~\eqref{eq:basic_pd_ms_strong_scale} and set
\begin{equation*}
  \begin{aligned}
    \overline\tau_k \eqdef\min\lb\overline\tau,
      8\max_{0\leq j\leq k}2^{j-k}\delta_j
        +\overline\tau
          \ls\frac{3\overline\tau}{3\overline\tau+2\mu}\rs^k\rb, \qquad
    \zeta_k \eqdef
      \ls\frac{3\overline\tau_k}{3\overline\tau_k+2\mu}\rs^2.
  \end{aligned}
\end{equation*}
Every completed iteration satisfies
\begin{gather}
  \tau_k\leq\overline\tau_k,
  \nonumber\\
  f(x_{k+1})-f^\ast\leq\zeta_k(f(x_k)-f^\ast),
  \qquad 0<\zeta_k<1.
  \label{eq:basic_pd_ms_refined_step}
\end{gather}
\end{theorem}

\begin{proof}
For every $i\geq1$, by~\eqref{eq:basic_pd_ms_output_gradient}
and~\eqref{eq:basic_pd_ms_strong_length},
\begin{equation*}
  \norm{\grad f(x_i)}
  \leq\frac{3\tau_{i-1}}2\norm{s_{i-1}}
  \leq\frac{3\tau_{i-1}}{\tau_{i-1}+2\mu}
    \norm{\grad f(x_{i-1})}.
\end{equation*}
Dividing by $\norm{\grad f(x_{i-1})}>0$ and substituting into the update
in line~\ref{line:basic_pd_ms_eta_update} gives
\begin{equation}
    \frac{\eta_i}{\tau_{i-1}}
    =\frac12\min\lb1,\frac{\norm{\grad f(x_i)}}
                              {\norm{\grad f(x_{i-1})}}\rb
    \leq\min\lb\frac12,
      \frac{3\tau_{i-1}}{2\tau_{i-1}+4\mu}\rb
    \leq\min\lb\frac12,
      \frac{3\overline\tau}{2\overline\tau+4\mu}\rb
    \leq\frac{3\overline\tau}{3\overline\tau+2\mu},
  \label{eq:basic_pd_ms_refined_first_trial}
\end{equation}
where we used $\tau_{i-1}\leq\overline\tau$ from
Theorem~\ref{thm:basic_pd_ms_strong_global}.
Fix a completed iteration $k$. 

If no trial was rejected during iterations $0,\ldots,k$, then
$\tau_i=\eta_i$ for $i=0,\ldots,k$. By~\eqref{eq:basic_pd_ms_refined_first_trial},
\begin{gather*}
  \tau_k=\eta_0\prod_{i=1}^k\frac{\eta_i}{\tau_{i-1}}
  \leq\eta_0\ls\frac{3\overline\tau}{3\overline\tau+2\mu}\rs^k
  \leq\overline\tau\ls\frac{3\overline\tau}{3\overline\tau+2\mu}\rs^k.
\end{gather*}
Otherwise, let $j$ be the last iteration among $0,\ldots,k$
on which a trial is rejected.
By~\eqref{eq:basic_pd_ms_strong_gap}
and~\eqref{eq:basic_pd_ms_strong_global_rate},
\begin{gather*}
  \norm{s_j}\leq\norm{x_j-x^\ast}+\norm{x_{j+1}-x^\ast}\\
  \leq\sqrt{\frac{2(f(x_j)-f^\ast)}\mu}
    +\sqrt{\frac{2(f(x_{j+1})-f^\ast)}\mu}
  \leq2\sqrt{\frac{2(f(x_0)-f^\ast)}\mu}
    \ls\frac{3\overline\tau}{3\overline\tau+2\mu}\rs^j.
\end{gather*}
Since a trial is rejected on iteration $j$, the two cases in the proof of
Theorem~\ref{thm:basic_pd_ms_strong_global} give the first bound below.
Substituting the step estimate above
and~\eqref{eq:basic_pd_ms_strong_global_rate}, we obtain
\begin{gather*}
  \tau_j\leq\max\lb8\delta_j+24M_0\norm{s_j},
    \frac{2(f(x_j)-f^\ast)}{r_0^2}\rb\\
  \leq\max\lb
    8\delta_j+48M_0\sqrt{\frac{2(f(x_0)-f^\ast)}\mu}
      \ls\frac{3\overline\tau}{3\overline\tau+2\mu}\rs^j,
    \frac{2(f(x_0)-f^\ast)}{r_0^2}
      \ls\frac{3\overline\tau}{3\overline\tau+2\mu}\rs^{2j}\rb.
\end{gather*}
By~\eqref{eq:basic_pd_ms_strong_scale} we obtain
\begin{equation}
  \tau_j\leq8\delta_j
    +\overline\tau\ls\frac{3\overline\tau}{3\overline\tau+2\mu}\rs^j.
  \label{eq:basic_pd_ms_refined_rejected}
\end{equation}
By the choice of $j$, no trial is rejected on iterations $j+1,\ldots,k$,
so $\tau_i=\eta_i$ for $i=j+1,\ldots,k$.
Applying~\eqref{eq:basic_pd_ms_refined_first_trial}, we obtain
\begin{gather*}
  \prod_{i=j+1}^k\frac{\eta_i}{\tau_{i-1}}
  \leq\prod_{i=j+1}^k\frac12=2^{j-k},\quad
  \prod_{i=j+1}^k\frac{\eta_i}{\tau_{i-1}}
  \leq\prod_{i=j+1}^k\frac{3\overline\tau}{3\overline\tau+2\mu}
  =\ls\frac{3\overline\tau}{3\overline\tau+2\mu}\rs^{k-j}.
\end{gather*}
Substituting~\eqref{eq:basic_pd_ms_refined_rejected}, we obtain
\begin{gather*}
  \tau_k=\tau_j\prod_{i=j+1}^k\frac{\eta_i}{\tau_{i-1}}
  \leq\ls8\delta_j+\overline\tau
    \ls\frac{3\overline\tau}{3\overline\tau+2\mu}\rs^j\rs
    \prod_{i=j+1}^k\frac{\eta_i}{\tau_{i-1}}
  \leq8\delta_j2^{j-k}+\overline\tau
    \ls\frac{3\overline\tau}{3\overline\tau+2\mu}\rs^j
    \ls\frac{3\overline\tau}{3\overline\tau+2\mu}\rs^{k-j}\\
  =8\delta_j2^{j-k}+\overline\tau
    \ls\frac{3\overline\tau}{3\overline\tau+2\mu}\rs^k
  \leq8\max_{0\leq i\leq k}2^{i-k}\delta_i+\overline\tau
    \ls\frac{3\overline\tau}{3\overline\tau+2\mu}\rs^k.
\end{gather*}
Thus, in both cases,
\begin{equation}
  \tau_k\leq8\max_{0\leq j\leq k}2^{j-k}\delta_j
    +\overline\tau\ls\frac{3\overline\tau}{3\overline\tau+2\mu}\rs^k.
  \label{eq:basic_pd_ms_refined_trial_bound}
\end{equation}
Combining~\eqref{eq:basic_pd_ms_refined_trial_bound}
with $\tau_k\leq\overline\tau$ gives
$\tau_k\leq\overline\tau_k$.

Since $0<\tau_k\leq\overline\tau_k$,
\begin{equation*}
  \ls\frac{3\tau_k}{3\tau_k+2\mu}\rs^2
  =\ls1-\frac{2\mu}{3\tau_k+2\mu}\rs^2
  \leq\ls1-\frac{2\mu}{3\overline\tau_k+2\mu}\rs^2
  =\zeta_k.
\end{equation*}
Substitution into~\eqref{eq:basic_pd_ms_strong_step}
proves~\eqref{eq:basic_pd_ms_refined_step}.
\end{proof}

\subsubsection{Global Superlinear Convergence under Relative Hessian Accuracy}

We consider Hessian approximations satisfying
\begin{equation}
  \delta_k=\norm{H(x_k)-\nabla^2 f(x_k)}
  \leq\nu\norm{\grad f(x_k)},\qquad k\geq0,
  \label{as:basic_pd_ms_relative_appendix}
\end{equation}
for some $\nu>0$.

\begin{corollary}
\label{cor:basic_pd_ms_relative_complexity}
Under the assumptions of Theorem~\ref{thm:basic_pd_ms_strong_global},
suppose additionally that~\eqref{as:basic_pd_ms_relative_appendix} holds.
The method either returns a minimizer or converges superlinearly
in function value:
\begin{equation*}
  \frac{f(x_{k+1})-f^\ast}{f(x_k)-f^\ast}
  \leq\ls\frac{3\tau_k}{3\tau_k+2\mu}\rs^2
  \longrightarrow0.
\end{equation*}
The first $N\geq1$ completed iterations require at most $N$ Hessian
evaluations and
\begin{equation}
  1+\frac94N
  +\log_{2,+}\ls
    \frac{\norm{\grad f(x_0)}}{\eta_0}
    \max\lb\frac{24M_0}{\mu}+8\nu,\frac2{r_0}\rb
  \rs
  \label{eq:basic_pd_ms_relative_gradient_cost}
\end{equation}
gradient evaluations.
\end{corollary}

\begin{proof}
Let $g_k\eqdef\norm{\grad f(x_k)}$.
In the proof of~\eqref{eq:basic_pd_ms_strong_gradient_ratio}, retain
the current error $\delta_k$ on iterations with a rejection and use
$\delta_k/g_k\leq\nu$ from~\eqref{as:basic_pd_ms_relative_appendix}.
The same bounds for the two cases give
\begin{gather*}
  \frac{\tau_0}{g_0}
  \leq\max\lb\frac{\eta_0}{g_0},
    \frac{24M_0}{\mu}+8\nu,\frac2{r_0}\rb, \qquad
  \frac{\tau_k}{g_k}
  \leq\max\lb\frac{\tau_{k-1}}{2g_{k-1}},
    \frac{24M_0}{\mu}+8\nu,\frac2{r_0}\rb,\quad k\geq1.
\end{gather*}
Iterating this inequality gives
\begin{equation}
  \frac{\tau_k}{g_k}
  \leq\max\lb\frac{2^{-k}\eta_0}{g_0},
    \frac{24M_0}{\mu}+8\nu,\frac2{r_0}\rb.
  \label{eq:basic_pd_ms_relative_ratio}
\end{equation}
By~\eqref{eq:basic_pd_ms_decrease}
and~\eqref{eq:basic_pd_ms_strong_global_rate}, for $k\geq1$,
\begin{gather*}
  g_k^2
  \leq\frac{3\tau_{k-1}}2\ls f(x_{k-1})-f(x_k)\rs
  \leq\frac{3\overline\tau}2\ls f(x_{k-1})-f^\ast\rs 
  \leq\frac{3\overline\tau(f(x_0)-f^\ast)}2
    \ls\frac{3\overline\tau}{3\overline\tau+2\mu}\rs^{2(k-1)}.
\end{gather*}
Substituting into~\eqref{eq:basic_pd_ms_relative_ratio} gives
\begin{equation*}
  \tau_k
  \leq\sqrt{\frac{3\overline\tau(f(x_0)-f^\ast)}2}
    \ls\frac{3\overline\tau}{3\overline\tau+2\mu}\rs^{k-1}
    \max\lb\frac{2^{-k}\eta_0}{g_0},
      \frac{24M_0}{\mu}+8\nu,\frac2{r_0}\rb
  \longrightarrow0.
\end{equation*}
Thus,~\eqref{eq:basic_pd_ms_strong_step} proves superlinear convergence.

Let $b_k$ count the rejected trials on iteration $k$. By~\eqref{eq:basic_pd_ms_rejections_boundary}
and~\eqref{eq:basic_pd_ms_gradient_log_sum},
\begin{equation*}
  \sum_{k=0}^{N-1}b_k
  \leq\frac54(N-1)
    +\log_2\frac{g_0\tau_{N-1}}{\eta_0g_{N-1}}.
\end{equation*}
By~\eqref{eq:basic_pd_ms_relative_ratio},
\begin{gather*}
  \frac{g_0\tau_{N-1}}{\eta_0g_{N-1}}
  \leq\max\lb2^{-(N-1)},
    \frac{g_0}{\eta_0}\ls\frac{24M_0}{\mu}+8\nu\rs,
    \frac{2g_0}{\eta_0r_0}\rb
  \leq\max\lb1,
    \frac{g_0}{\eta_0}\ls\frac{24M_0}{\mu}+8\nu\rs,
    \frac{2g_0}{\eta_0r_0}\rb.
\end{gather*}
Hence
\begin{equation*}
  \sum_{k=0}^{N-1}b_k
  \leq\frac54(N-1)
    +\log_{2,+}\ls\frac{g_0}{\eta_0}
      \max\lb\frac{24M_0}{\mu}+8\nu,\frac2{r_0}\rb\rs.
\end{equation*}
The gradient count is at most $1+N+\sum_{k=0}^{N-1}b_k$,
which proves~\eqref{eq:basic_pd_ms_relative_gradient_cost}.
\end{proof}

\subsubsection{Local Quadratic Convergence}

\begin{theorem}
\label{thm:basic_pd_ms_quadratic}
Under the assumptions of Corollary~\ref{cor:basic_pd_ms_relative_complexity},
set
\begin{equation*}
  C\eqdef\frac{9\ls\norm{\nabla^2 f(x^\ast)}+\mu\rs}{2\mu^2}
    \ls8\nu+\frac{24M_0}{\mu}\rs^2.
\end{equation*}
Unless a minimizer is returned earlier, let $K_0\geq1$ satisfy
\begin{equation}
  f(x_{K_0})-f^\ast
  \leq\min\lb\frac{\mu r_0^2}{8},\frac1{2C}\rb
  \label{eq:basic_pd_ms_quadratic_entry}
\end{equation}
and define
\begin{equation}
  K\eqdef K_0+\left\lceil\log_{2,+}
    \frac{\tau_{K_0-1}}
    {2\ls8\nu+\frac{24M_0}{\mu}\rs
      \sqrt{2\mu(f(x_{K_0-1})-f^\ast)}}
    \right\rceil.
  \label{eq:basic_pd_ms_quadratic_entry_index}
\end{equation}
Then, for every $k\geq K$,
\begin{equation}
  f(x_{k+1})-f^\ast\leq C\ls f(x_k)-f^\ast\rs^2.
  \label{eq:basic_pd_ms_quadratic_gap}
\end{equation}
The gradient norms and the distances to $x^\ast$ also converge
quadratically. After a finite initial cost independent of $\e$,
attaining $f(x)-f^\ast\leq\e$ requires
$\cO(\log\log(1/\e))$ additional iterations, Hessian evaluations,
and gradient evaluations.
\end{theorem}

\begin{proof}
By Theorem~\ref{thm:basic_pd_ms_strong_global}, an index $K_0$
satisfying~\eqref{eq:basic_pd_ms_quadratic_entry} exists unless
a minimizer is returned earlier.
Since $\norm{\nabla^2 f(x^\ast)}\geq\mu$, the definition of $C$ gives
\begin{equation*}
  C\geq\frac9\mu\ls\frac{24M_0}{\mu}\rs^2
  \geq\frac{4M_0^2}{\mu^3},
  \qquad
  \frac1{2C}\leq\frac{\mu^3}{8M_0^2}.
\end{equation*}
Together with monotonicity and~\eqref{eq:basic_pd_ms_strong_gap},
condition~\eqref{eq:basic_pd_ms_quadratic_entry} gives, for $k\geq K_0$,
\begin{gather*}
  \norm{x_k-x^\ast}
  \leq\sqrt{\frac{2(f(x_k)-f^\ast)}{\mu}}
  \leq\frac12\min\lb r_0,\frac{\mu}{M_0}\rb,\\
  \norm{s_k}
  \leq\norm{x_k-x^\ast}+\norm{x_{k+1}-x^\ast}
  \leq2\sqrt{\frac{2(f(x_{K_0})-f^\ast)}{\mu}}
  \leq r_0.
\end{gather*}
Let $g_k\eqdef\norm{\grad f(x_k)}$.
Since $\norm{s_k}\leq r_0$, the proof
of~\eqref{eq:basic_pd_ms_relative_ratio} applies without the term
$2/r_0$. Starting at $K_0$, it gives, for $j\geq0$,
\begin{gather*}
  \frac{\tau_{K_0+j}}{g_{K_0+j}}
  \leq\max\lb
    \frac{2^{-j-1}\tau_{K_0-1}}{g_{K_0-1}},
    8\nu+\frac{24M_0}{\mu}\rb
  \leq\max\lb
    \frac{2^{-j-1}\tau_{K_0-1}}
      {\sqrt{2\mu(f(x_{K_0-1})-f^\ast)}},
    8\nu+\frac{24M_0}{\mu}\rb.
\end{gather*}
By~\eqref{eq:basic_pd_ms_quadratic_entry_index}, for every $k\geq K$,
\begin{equation}
  \tau_k\leq\ls8\nu+\frac{24M_0}{\mu}\rs g_k.
  \label{eq:basic_pd_ms_quadratic_shift}
\end{equation}
Combining~\eqref{eq:basic_pd_ms_output_gradient}
and~\eqref{eq:basic_pd_ms_strong_length}, we obtain
\begin{equation}
  g_{k+1}
  \leq\frac{3\tau_k}{2}\norm{s_k}
  \leq\frac{3\tau_k}{2\mu}g_k
  \leq\frac{3}{2\mu}\ls8\nu+\frac{24M_0}{\mu}\rs g_k^2.
  \label{eq:basic_pd_ms_quadratic_gradient}
\end{equation}
Thus, the gradient norms converge quadratically.
Next, we relate the gradient norm to the function gap.
For $\norm{z-x^\ast}\leq\min\lb r_0,\mu/M_0\rb$,
either~\eqref{as:gh} or~\eqref{as:lh} gives
\begin{equation*}
  \norm{\nabla^2 f(z)}
  \leq\norm{\nabla^2 f(x^\ast)}+M_0\norm{z-x^\ast}
  \leq\norm{\nabla^2 f(x^\ast)}+\mu,
\end{equation*}
since $\grad f(x^\ast)=0$ and $r_0\leq\rho/2$ under~\eqref{as:lh}.
Thus, $f$ has $(\norm{\nabla^2 f(x^\ast)}+\mu)$-Lipschitz gradient on this ball.
Hence, for $k\geq K_0$,
\begin{gather*}
  g_k\leq\ls\norm{\nabla^2 f(x^\ast)}+\mu\rs\norm{x_k-x^\ast},\\
  \left \| x_k-\frac{\grad f(x_k)}{\norm{\nabla^2 f(x^\ast)}+\mu}-x^\ast \right \|
  \leq\norm{x_k-x^\ast}+\frac{g_k}{\norm{\nabla^2 f(x^\ast)}+\mu}
  \leq2\norm{x_k-x^\ast}\leq\min\lb r_0,\frac{\mu}{M_0}\rb.
\end{gather*}
By Lipschitznes of the gradient, for any $y$ in the ball,
\begin{equation*}
  f(y)\leq f(x_k)+\la\grad f(x_k),y-x_k\ra
  +\frac{\norm{\nabla^2 f(x^\ast)}+\mu}{2}\norm{y-x_k}^2.
\end{equation*}
Substituting the point above, we obtain
\begin{gather*}
  f^\ast\leq f\ls x_k-\frac{\grad f(x_k)}{\norm{\nabla^2 f(x^\ast)}+\mu}\rs
  \leq f(x_k)-\frac{g_k^2}{\norm{\nabla^2 f(x^\ast)}+\mu}
  +\frac{\ls\norm{\nabla^2 f(x^\ast)}+\mu\rs g_k^2}
    {2\ls\norm{\nabla^2 f(x^\ast)}+\mu\rs^2}\\
  =f(x_k)-\frac{g_k^2}{2\ls\norm{\nabla^2 f(x^\ast)}+\mu\rs}.
\end{gather*}
Rearranging gives
\begin{equation*}
  g_k^2\leq2\ls\norm{\nabla^2 f(x^\ast)}+\mu\rs(f(x_k)-f^\ast).
\end{equation*}
Substituting this and~\eqref{eq:basic_pd_ms_quadratic_shift}
into~\eqref{eq:basic_pd_ms_strong_step}, for $k\geq K$ we get
\begin{gather*}
  f(x_{k+1})-f^\ast
  \leq\ls\frac{3\tau_k}{3\tau_k+2\mu}\rs^2
    (f(x_k)-f^\ast)\\
  \leq\frac{9}{4\mu^2}\ls8\nu+\frac{24M_0}{\mu}\rs^2
    g_k^2(f(x_k)-f^\ast)
  \leq C(f(x_k)-f^\ast)^2.
\end{gather*}
This proves~\eqref{eq:basic_pd_ms_quadratic_gap}.
For the distances, strong convexity
and~\eqref{eq:basic_pd_ms_quadratic_gradient} give
\begin{equation*}
  \norm{x_{k+1}-x^\ast}
  \leq\frac{g_{k+1}}{\mu}
  \leq\frac{3\ls\norm{\nabla^2 f(x^\ast)}+\mu\rs^2}{2\mu^2}
    \ls8\nu+\frac{24M_0}{\mu}\rs\norm{x_k-x^\ast}^2.
\end{equation*}
By~\eqref{eq:basic_pd_ms_quadratic_entry},
$C(f(x_K)-f^\ast)\leq1/2$.
Iterating~\eqref{eq:basic_pd_ms_quadratic_gap} gives
\begin{equation*}
  f(x_{K+n})-f^\ast
  \leq\frac1C\ls C(f(x_K)-f^\ast)\rs^{2^n}
  \leq\frac{2^{-2^n}}C,\qquad n\geq0.
\end{equation*}
Thus, $f(x_{K+n})-f^\ast\leq\e$ after
$n=\cO(\log\log(1/\e))$ additional iterations.
By Corollary~\ref{cor:basic_pd_ms_relative_complexity}, this requires
$\cO(K+\log\log(1/\e))$ Hessian and gradient evaluations in total.
\end{proof}

\section{Nesterov-type Acceleration}
\label{sec:nesterov_acceleration}

We present an accelerated cubic method for the global
condition~\eqref{as:gh}. Fix a positive radius estimate $\bar R\geq R$ and an initial scale $\Gamma_0>0$ satisfying 
\begin{equation}
  \label{eq:inital_cnd}
  \Gamma_0\bar R^3\geq 3\ls f(x_0)-f^\ast\rs.
\end{equation}
We discuss how to choose parameters satisfying~\eqref{eq:inital_cnd}
in Section~\ref{sec:fully_adaptive}.

We present the resulting method as Algorithm~\ref{alg:global_accelerated_cubic}. 
The acceleration mechanism aggregates the linear models
$\ell_w(x)\eqdef f(w)+\la\grad f(w),x-w\ra$. At each iteration, $\alpha_k$
is defined by the scalar equation in line~\ref{line:alpha_k} of
Algorithm~\ref{alg:global_accelerated_cubic}. This definition is well posed
because the mapping $\alpha\mapsto\alpha^3/(1-\alpha)$ is continuous and
strictly increasing from $0$ to $+\infty$ on $(0,1)$. Hence, for
$\Gamma_k>0$ and $\M_k>0$, the equation has a unique solution, which can be
found by bisection.

\begin{algorithm}
  \caption{Accelerated cubic regularization under global generalized Hessian smoothness}
  \label{alg:global_accelerated_cubic}
  \begin{algorithmic}[1]
    \STATE \textbf{Input:} $x_0\in\R^d$, $M_0>0$,
      $M_1\geq0$, $\bar R\geq R$ with $\bar R>0$, and $\Gamma_0>0$ satisfying
      $\Gamma_0\bar R^3\geq3\ls f(x_0)-f^\ast\rs$
    \STATE Set $A_0\eqdef1$, $w_0\eqdef x_0$, and $z_0\eqdef x_0$, $\psi_0(x)\eqdef f(x_0)+\tfrac{\Gamma_0}{3}\norm{x-x_0}^3$.
    \FOR{$k \geq 0$}
      \STATE Set
        \begin{equation}
          \label{eq:regularization_parameter}
          \M_k\eqdef
          \max\lb
            2M_0,
            12\sqrt{2}M_1^{3/2}\bar R^3\Gamma_k
          \rb.
        \end{equation}
      \STATE Find the unique $\alpha_k\in(0,1)$ satisfying
        $\tfrac{\alpha_k^3}{1-\alpha_k}
        =\tfrac{\Gamma_k}{72\M_k}$.\label{line:alpha_k}
      \STATE Set 
        \begin{equation}
          \label{eq:a_params}
          A_{k+1}\eqdef\tfrac{A_k}{1-\alpha_k}, \quad a_{k+1}\eqdef\alpha_kA_{k+1}, \quad \Gamma_{k+1}\eqdef(1-\alpha_k)\Gamma_k.
        \end{equation}
      \STATE Set $v_k\eqdef(1-\alpha_k)w_k+\alpha_kz_k$.
      \STATE Compute
        \begin{equation}
          \label{eq:cubic_step}
          s_k\eqdef\argmin_{s\in\R^d}\lb
            \la\grad f(v_k),s\ra
            +\tfrac{1}{2}\la\nabla^2 f(v_k)s,s\ra
            +\tfrac{15\M_k}{6}\norm{s}^3
          \rb,\quad w_{k+1}\eqdef v_k+s_k.
        \end{equation}
      \STATE Set $\psi_{k+1}(x)\eqdef
        \psi_k(x)+a_{k+1}\ell_{w_{k+1}}(x)$ and compute $z_{k+1}\eqdef\argmin_{x\in\R^d}\psi_{k+1}(x)$.
    \ENDFOR
  \end{algorithmic}
\end{algorithm}

\begin{theorem}
\label{thm:global_convergence}
Assume that $f$ is convex,~\eqref{as:gh} and~\eqref{eq:inital_cnd} hold. For every
$\e>0$, Algorithm~\ref{alg:global_accelerated_cubic} returns an iterate
$w_{N}$ satisfying $f(w_{N})-f^\ast\leq\e$ for
\begin{equation}
  \label{eq:gh_conv_rate}
  N = \cO \ls
    \ls\tfrac{M_0\bar R^3}{\e}\rs^{1/3} + \ls1+\sqrt{M_1}\bar R\rs
    \log\max\lb
      1,
      \min\lb
        \tfrac{M_1^{3/2}\bar R^3\Gamma_0}{M_0},
        \tfrac{\Gamma_0\bar R^3}{\e}
      \rb
    \rb
  \rs.
\end{equation}
\end{theorem}

\textbf{Proof sketch.}
We use the standard accelerated-cubic estimating sequence technique~\cite{nesterov2008accelerating}. The main obstruction is that under~\eqref{as:gh} the regularization coefficient $\M_k$ must dominate the generalized smoothness parameter $L(w_{k+1})$ at the \emph{output} $w_{k+1}$ of the cubic step~\eqref{eq:cubic_step}, although this output is not known when the coefficient is selected. Inspired by~\cite{tyurin2025near}, we provide an explicit way of choosing regularization parameter $\M_k$~\eqref{eq:regularization_parameter}, such that this condition is satisfied. The complete proof is given in Appendix~\ref{app:global_method_proof}.

Let us consider one step of Algorithm~\ref{alg:global_accelerated_cubic}. For some $\M>0$, let $A_+(\M)$, $w(\M)$, and $\psi_+(\cdot;\M)$ denote the
quantities produced by one iteration, and set
$\psi_+^\ast(\M)=\min_{x\in\R^d}\psi_+(x;\M)$. Assume first that
$\M\geq L(w(\M))$. Then the standard accelerated-cubic estimating-sequence argument gives
\begin{equation}
  \label{eq:sketch_strengthened_estimate}
  A_+(\M)f(w(\M))
  +\tfrac{A_+(\M)}{9}
  \tfrac{\norm{\grad f(w(\M))}^{3/2}}{\sqrt{\M}}
  \leq\psi_+^\ast(\M)
  \leq A_+(\M)f^\ast+\tfrac{2}{3}\Gamma_0\bar R^3.
\end{equation}
Thus, conditioned on proper choice of $\M \geq L(w(\M))$, the usual estimating-sequence invariant closes. It remains to choose $\M$ before $w(\M)$ is known. Consider the virtual family of steps and
set $\mathcal U_k=\lb \M>0:\M<L(w(\M))\rb$. If this set is nonempty, let
$\M_k^\dagger=\sup\mathcal U_k$. Since all sufficiently large $\M$ are admissible, $\M_k^\dagger<+\infty$,
and continuity gives
$\M_k^\dagger=L(w(\M_k^\dagger))$. Let $g^\dagger=\grad f(w(\M_k^\dagger))$. Combining the two sides of
\eqref{eq:sketch_strengthened_estimate} at this boundary yields
\begin{equation*}
  \tfrac{\norm{g^\dagger}^{3/2}}{\sqrt{\M_k^\dagger}}
  \leq6\Gamma_k\bar R^3,
  \qquad
  \M_k^\dagger
  \leq\max\lb
    2M_0,
    12\sqrt{2}M_1^{3/2}\bar R^3\Gamma_k
  \rb
  =\M_k.
\end{equation*}
Hence the explicit $\M_k$ used by the algorithm is admissible. This is why we retain the positive remainder
in~\eqref{eq:sketch_strengthened_estimate}. Although it can be
omitted in the standard convergence argument, here it controls the
gradient at the unknown output and closes the explicit choice of $\M$. For the actual step, admissibility gives
\begin{equation*}
  A_{k+1} (f(w_{k+1})-f^\ast)
  \leq\psi_{k+1}^\ast-A_{k+1}f^\ast
  \leq\tfrac{2}{3}\Gamma_0\bar R^3.
\end{equation*}
Then, recurrence for $\Gamma_k$ then gives the rate~\eqref{eq:gh_conv_rate}.

\subsection{Proof of Nesterov-type Acceleration under~\eqref{as:gh}}
\label{app:global_method_proof}

\subsubsection{Auxiliary Lemmas}

\begin{lemma}[{\citet[Lemma~2]{nesterov2008accelerating}}]
\label{lem:cubic_conjugate}
For any $a\geq0$, $g\in\R^d$, and $\beta>0$,
\begin{equation}
  \min_{h\in\R^d}
  \lb
    \frac{\beta}{3!}\norm{h}^3+a\la g,h\ra
  \rb
  =-\frac{2\sqrt{2}}{3}
  \frac{a^{3/2}\norm{g}^{3/2}}{\sqrt{\beta}}.
  \label{eq:cubic_conjugate}
\end{equation}
\end{lemma}

\begin{lemma}[{\citet[Lemma~7]{agafonov2023inexact}}]
\label{lem:cubic_prox}
Let $h:\R^d\to\R$ be convex, $x_0\in\R^d$, and $\theta>0$. Define
\begin{equation*}
  \varphi(x)
  \eqdef
  h(x)
  +\frac{\theta}{3}\norm{x-x_0}^3,
  \qquad
  \bar x\eqdef\argmin_{x\in\R^d}\varphi(x).
\end{equation*}
Then, for every $u\in\R^d$,
\begin{equation*}
  \varphi(\bar x+u)
  \geq
  \varphi(\bar x)
  +\frac{\theta}{6}\norm{u}^3.
\end{equation*}
\end{lemma}

\begin{lemma}
\label{lem:equivalent_weights}
The weight updates in Algorithm~\ref{alg:global_accelerated_cubic} admit the
following equivalent representation:
\begin{gather}
  a_{k+1} \eqdef A_{k+1}-A_k,~~
  \Gamma_k=\frac{\Gamma_0}{A_k},~~
  v_k=\frac{A_kw_k+a_{k+1}z_k}{A_{k+1}},~~
  a_{k+1}^3=\frac{\Gamma_0A_{k+1}^2}{72\M_k}.
  \label{eq:equivalent_weights}
\end{gather}
\end{lemma}

\begin{proof}
Since $A_0=1$, the updates of $A_k$ and $\Gamma_k$~\eqref{eq:a_params} give
\begin{gather*}
  A_{k+1}\Gamma_{k+1}
  =\frac{A_k}{1-\alpha_k}(1-\alpha_k)\Gamma_k
  =A_k\Gamma_k
  =\Gamma_0,
\end{gather*}
and hence $\Gamma_k=\Gamma_0/A_k$. Moreover,
\begin{gather*}
  A_{k+1}-A_k
  =A_{k+1}-(1-\alpha_k)A_{k+1}
  =\alpha_kA_{k+1}
  =a_{k+1},\\
  A_{k+1}v_k
  =(1-\alpha_k)A_{k+1}w_k+\alpha_kA_{k+1}z_k
  =A_kw_k+a_{k+1}z_k.
\end{gather*}
Finally,
\begin{gather*}
  \frac{\alpha_k^3}{1-\alpha_k}
  =\frac{\Gamma_k}{72\M_k}
  \quad\Longleftrightarrow\quad
  \alpha_k^3A_{k+1}^3
  =\frac{\Gamma_kA_kA_{k+1}^2}{72\M_k}
  \quad\Longleftrightarrow\quad
  a_{k+1}^3
  =\frac{\Gamma_0A_{k+1}^2}{72\M_k}.
\end{gather*}
\end{proof}

\subsubsection{Cubic Step Progress}

\begin{lemma}
\label{lem:backward_progress}
Fix $v\in\R^d$ and $\M>0$, and define
\begin{equation}
  s\eqdef\argmin_{h\in\R^d}
  \lb
    \la\grad f(v),h\ra
    +\frac{1}{2}\la\nabla^2 f(v)h,h\ra
    +\frac{15\M}{6}\norm{h}^3
  \rb.
  \label{eq:local_model}
\end{equation}
Set $w\eqdef v+s$. Suppose that either~\eqref{as:gh} holds, or~\eqref{as:lh} holds and $\norm{s}\leq\rho$. If
\begin{equation*}
  \M\geq L(w)=M_0+M_1\norm{\grad f(w)},
\end{equation*}
then
\begin{gather}
  \norm{\grad f(w)}
  \leq 9\M\norm{s}^2,
  \label{eq:q_step}\\
  \la\grad f(w),v-w\ra
  \geq 6\M\norm{s}^3,
  \label{eq:inner_step}\\
  \la\grad f(w),v-w\ra
  \geq\frac{2}{9}
  \frac{\norm{\grad f(w)}^{3/2}}{\sqrt{\M}}.
  \label{eq:progress_final}
\end{gather}
\end{lemma}

\begin{proof}
By optimality condition for the cubic subproblem~\eqref{eq:local_model}
\begin{equation}
  \label{eq:opt_cnd}
  \grad f(v)+\nabla^2 f(v)s+\frac{15\M}{2}\norm{s}s=0.
\end{equation}
Under~\eqref{as:gh}, Lemma~\ref{lem:global_taylor_remainders} applies at $w$ with $v-w=-s$. Under~\eqref{as:lh}, Lemma~\ref{lem:local_taylor_remainders} applies in the same way, since $\norm{s}\leq\rho$. Together with the corresponding smoothness condition, this gives
\begin{gather*}
  \norm{\grad f(v)-\grad f(w)+\nabla^2 f(w)s}
  \leq\frac{L(w)}{2}\norm{s}^2
  \leq\frac{\M}{2}\norm{s}^2,\\
  \norm{\nabla^2 f(v)-\nabla^2 f(w)}
  \leq L(w)\norm{s}
  \leq \M\norm{s}.
\end{gather*}
By the triangle inequality
\begin{gather*}
  \norm{\grad f(v)-\grad f(w)+\nabla^2 f(v)s}
  \\
  \leq \norm{\grad f(v)-\grad f(w)+\nabla^2 f(w)s}
  +\norm{\nabla^2 f(v)-\nabla^2 f(w)}\norm{s} \leq\frac{3\M}{2}\norm{s}^2.
\end{gather*}
Combining this estimate with~\eqref{eq:opt_cnd}
\begin{gather*}
  \norm{\grad f(w)}
  =
  \norm{\grad f(w) - \ls \grad f(v)+\nabla^2 f(v)s+\tfrac{15\M}{2}\norm{s}s\rs}
  \\
  \leq
  \norm{\grad f(v)-\grad f(w)+\nabla^2 f(v)s}
  +\frac{15\M}{2}\norm{s}^2
  \leq9\M\norm{s}^2,
\end{gather*}
which proves \eqref{eq:q_step}. 

Next, 
\begin{gather*}
  \la\grad f(w),-s\ra
  \stackrel{\eqref{eq:opt_cnd}}{=}
  \la\grad f(v)-\grad f(w)+\nabla^2 f(v)s,s\ra
  +\frac{15\M}{2}\norm{s}^3,\\
  \la\grad f(w),-s\ra
  \geq-\frac{3\M}{2}\norm{s}^3+\frac{15\M}{2}\norm{s}^3
  =6\M\norm{s}^3,
\end{gather*}
which proves \eqref{eq:inner_step}. Finally, \eqref{eq:q_step} implies
\begin{equation*}
  \norm{s}\geq
  \sqrt{\frac{\norm{\grad f(w)}}{9\M}}.
\end{equation*}
Substituting this bound into \eqref{eq:inner_step} proves
\eqref{eq:progress_final}.
\end{proof}

\subsubsection{Estimating Sequence One Step Analysis}

\begin{lemma}
\label{lem:one_safe_step}
Suppose that, at the beginning of iteration $k$,
\begin{equation}
  \psi_k^\ast\eqdef\min_{x\in\R^d}\psi_k(x)
  \geq A_kf(w_k).
  \label{eq:current_invariant}
\end{equation}
Then
\begin{equation}
  \psi_{k+1}^\ast
  \geq A_{k+1}f(w_{k+1})
  +A_{k+1}\la\grad f(w_{k+1}),v_k-w_{k+1}\ra
  -\frac{A_{k+1}}{9}
  \frac{\norm{\grad f(w_{k+1})}^{3/2}}{\sqrt{\M_k}}.
  \label{eq:raw_one_step}
\end{equation}
If, in addition, the hypotheses of Lemma~\ref{lem:backward_progress} hold for $v=v_k$ and $\M=\M_k$, then
\begin{equation}
  \psi_{k+1}^\ast
  \geq A_{k+1}f(w_{k+1})
  +\frac{A_{k+1}}{9}
  \frac{\norm{\grad f(w_{k+1})}^{3/2}}{\sqrt{\M_k}}.
  \label{eq:strengthened_lower}
\end{equation}
Consequently,
\begin{equation}
  \psi_{k+1}^\ast\geq A_{k+1}f(w_{k+1}).
  \label{eq:invariant_update}
\end{equation}
\end{lemma}

\begin{proof}
By the definition of $\psi_k$,
\begin{equation*}
  \psi_k(x)
  =f(x_0)+\sum_{i=0}^{k-1}a_{i+1}\ell_{w_{i+1}}(x)
  +\frac{\Gamma_0}{3}\norm{x-x_0}^3.
\end{equation*}
The first two terms are affine and hence convex. Therefore,
Lemma~\ref{lem:cubic_prox} with $\theta=\Gamma_0$ gives for any $h \in \R^d$
\begin{equation}
  \psi_k(z_k+h)
  \geq\psi_k^\ast+\frac{\Gamma_0}{6}\norm{h}^3.
  \label{eq:psi_lower_bound}
\end{equation}
Next, using assumption of this Lemma~\eqref{eq:current_invariant} and Lemma~\ref{lem:cubic_conjugate},
\begin{gather}
  \psi_{k+1}^\ast
  =\min_{x\in\R^d}\psi_{k+1}(x)
  =\min_{h\in\R^d}\lb
    \psi_k(z_k+h)+a_{k+1}\ell_{w_{k+1}}(z_k+h)
  \rb
  \notag\\
  =a_{k+1}f(w_{k+1})
  +\min_{h\in\R^d}\lb
    \psi_k(z_k+h)
    +a_{k+1}\la\grad f(w_{k+1}),z_k+h-w_{k+1}\ra
  \rb
  \notag\\
  \stackrel{\eqref{eq:psi_lower_bound}}{\geq}
  \psi_k^\ast+a_{k+1}f(w_{k+1})
  +a_{k+1}\la\grad f(w_{k+1}),z_k-w_{k+1}\ra
  \notag\\
  +\min_{h\in\R^d}\lb
    \frac{\Gamma_0}{6}\norm{h}^3
    +a_{k+1}\la\grad f(w_{k+1}),h\ra
  \rb
  \notag\\
  \stackrel{\eqref{eq:current_invariant}}{\geq}
  A_kf(w_k)+a_{k+1}f(w_{k+1})
  +a_{k+1}\la\grad f(w_{k+1}),z_k-w_{k+1}\ra
  \notag\\
  +\min_{h\in\R^d}\lb
    \frac{\Gamma_0}{6}\norm{h}^3
    +a_{k+1}\la\grad f(w_{k+1}),h\ra
  \rb
  \notag\\
  \stackrel{\eqref{eq:cubic_conjugate}}{=}
  A_kf(w_k)+a_{k+1}f(w_{k+1})
  +a_{k+1}\la\grad f(w_{k+1}),z_k-w_{k+1}\ra
  -\tfrac{2\sqrt{2}}{3}
  \tfrac{a_{k+1}^{3/2}\norm{\grad f(w_{k+1})}^{3/2}}
  {\sqrt{\Gamma_0}}.
  \label{eq:one_step_lower}
\end{gather}
Relations~\eqref{eq:equivalent_weights} imply
\begin{gather}
  A_{k+1}v_k=A_kw_k+a_{k+1}z_k,
  \qquad
  A_{k+1}w_{k+1}=A_kw_{k+1}+a_{k+1}w_{k+1},
  \notag\\
  a_{k+1}(z_k-w_{k+1})
  =A_{k+1}(v_k-w_{k+1})-A_k(w_k-w_{k+1}).
  \label{eq:interpolation_identity}
\end{gather}
Consequently,
\begin{gather}
  A_kf(w_k)+a_{k+1}f(w_{k+1})
  +a_{k+1}\la\grad f(w_{k+1}),z_k-w_{k+1}\ra
  \notag\\
  \stackrel{\eqref{eq:interpolation_identity}}{=}
  A_{k+1}f(w_{k+1})
  +A_{k+1}\la\grad f(w_{k+1}),v_k-w_{k+1}\ra
  \notag\\
  \phantom{\stackrel{\eqref{eq:interpolation_identity}}{=}{}}
  +A_k\ls
    f(w_k)-f(w_{k+1})
    -\la\grad f(w_{k+1}),w_k-w_{k+1}\ra
  \rs
  \notag\\
  \geq A_{k+1}f(w_{k+1})
  +A_{k+1}\la\grad f(w_{k+1}),v_k-w_{k+1}\ra,
  \label{eq:convex_combination_bound}
\end{gather}
where the inequality follows from convexity of $f$. By~\eqref{eq:equivalent_weights},
\begin{equation*}
  \frac{2\sqrt{2}}{3}
  \frac{a_{k+1}^{3/2}}{\sqrt{\Gamma_0}}
  =\frac{A_{k+1}}{9\sqrt{\M_k}}.
\end{equation*}
Substituting this identity and~\eqref{eq:convex_combination_bound} into~\eqref{eq:one_step_lower} proves~\eqref{eq:raw_one_step}. Under the hypotheses of Lemma~\ref{lem:backward_progress}, substituting~\eqref{eq:progress_final} into~\eqref{eq:raw_one_step} gives
\begin{equation*}
  \psi_{k+1}^\ast
  \geq A_{k+1}f(w_{k+1})
  +\ls
    \frac{2A_{k+1}}{9\sqrt{\M_k}}
    -\frac{A_{k+1}}{9\sqrt{\M_k}}
  \rs
  \norm{\grad f(w_{k+1})}^{3/2},
\end{equation*}
which proves~\eqref{eq:strengthened_lower}. Dropping the nonnegative second term
proves~\eqref{eq:invariant_update}.
\end{proof}

\begin{lemma}
\label{lem:upper_estimate}
For every $k\ge0$,
\begin{equation}
  \label{eq:upper_estimate}
  \psi_k(x^\ast)
  \leq A_kf^\ast+\Delta_0+\frac{\Gamma_0R^3}{3},
\end{equation}
where $\Delta_0\eqdef f(x_0)-f^\ast$.
Moreover, under~\eqref{eq:inital_cnd},
\begin{equation}
  \label{eq:upper_estimate_exact}
  \Delta_0+\frac{\Gamma_0R^3}{3}\leq\frac{2}{3}\Gamma_0\bar R^3.
\end{equation}
\end{lemma}

\begin{proof}
By convexity, $\ell_w(x^\ast)\le f^\ast$ for every $w$. Since $A_k=1+\sum_{i=0}^{k-1}a_{i+1}$,
\begin{gather*}
 \psi_k(x^\ast)
 =f(x_0)+\frac{\Gamma_0}{3}R^3
 +\sum_{i=0}^{k-1}a_{i+1}\ell_{w_{i+1}}(x^\ast)\\
 \leq f(x_0)+\frac{\Gamma_0}{3}R^3+(A_k-1)f^\ast
 =A_kf^\ast+\Delta_0+\frac{\Gamma_0R^3}{3}.
\end{gather*}
Finally,~\eqref{eq:inital_cnd} gives
$\Delta_0\leq\Gamma_0\bar R^3/3$, and $R\leq\bar R$ gives
$\Gamma_0R^3/3\leq\Gamma_0\bar R^3/3$, proving the statement.
\end{proof}

\subsubsection{Explicit Choice of Regularization Parameter}

Fix an iteration $k$ and suppose that~\eqref{eq:current_invariant} holds.
For $\M>0$, let $\alpha(\M)\in(0,1)$ solve
\begin{equation}
  \label{eq:alpha_M_equality}
  \frac{\alpha(\M)^3}{1-\alpha(\M)}
  =\frac{\Gamma_k}{72\M}.
\end{equation}
Define
\begin{gather}
  A_+(\M)\eqdef\frac{A_k}{1-\alpha(\M)},
  \qquad
  a(\M)\eqdef\alpha(\M)A_+(\M), \notag \\
  v(\M)\eqdef\ls1-\alpha(\M)\rs w_k+\alpha(\M)z_k,
  \label{eq:interpolation_point_M}
\end{gather}
and let $w(\M)$ be the output of the cubic subproblem~\eqref{eq:local_model}
with base point $v(\M)$ and regularization parameter $\M$. Finally, set
\begin{equation*}
  \psi_+(x;\M)\eqdef\psi_k(x)+a(\M)\ell_{w(\M)}(x),
  \qquad
  \psi_+^\ast(\M)\eqdef\min_{x\in\R^d}\psi_+(x;\M).
\end{equation*}
We call $\M$ admissible if $\M\geq L(w(\M))$. The proof of
Lemma~\ref{lem:one_safe_step} directly applies to these quantities,
since the weight identities in Lemma~\ref{lem:equivalent_weights} remain valid.
Thus, every admissible $\M$ satisfies
\begin{equation}
  \psi_+^\ast(\M)
  \geq A_+(\M)f(w(\M))
  +\frac{A_+(\M)}{9}
  \frac{\norm{\grad f(w(\M))}^{3/2}}{\sqrt{\M}}.
  \label{eq:virtual_safe_step}
\end{equation}
It remains to verify admissibility for the explicit choice in
Algorithm~\ref{alg:global_accelerated_cubic}.

\begin{lemma}
\label{lem:virtual_step_continuity}
The mappings $\M\mapsto\norm{w(\M)-v(\M)}$ and $\M\mapsto\M-L(w(\M))$
are continuous on $(0,+\infty)$, and
\begin{equation}
  \norm{w(\M)-v(\M)}\to0,
  \qquad
  \M-L(w(\M))\to+\infty
  \quad\text{as }\M\to+\infty.
  \label{eq:M_L_limit}
\end{equation}
\end{lemma}

\begin{proof}
The mapping $\alpha\mapsto\alpha^3/(1-\alpha)$ is continuous and strictly
increasing on $(0,1)$. Hence $\alpha(\M)$ and $v(\M)$ are continuous.

Let $s(\M)\eqdef w(\M)-v(\M)$. For any sequence $\M_j\to \M>0$, the optimality
condition of the cubic subproblem and $\nabla^2f(v(\M_j))\succeq0$ give
\begin{equation}
  \label{eq:s_bound}
  \norm{s(\M_j)}^2
  \leq\frac{2\norm{\grad f(v(\M_j))}}{15\M_j}.
\end{equation}
By~\eqref{eq:interpolation_point_M} $\alpha(\M_j) \to \alpha(\M)$ as $\M_j \to \M$, therefore $v(\M_j) \to v(\M)$. Since $f$ is twice continuously differentiable $\nabla f(v(\M_j)) \to \nabla f(v(\M))$.
Thus $\{\nabla f(v(\M_j))\}$ is bounded and $\{s(\M_j)\}$ is bounded as well. 

Next, we show that $s(\M_j) \to s(\M)$. By Bolzano-Weierstrass theorem bounded sequence $\{s(\M_j)\}$ has a convergent subsequence $\{s(\M_{j_l})\} \to \hat{s}$. Taking the limit in the optimality condition of cubic subproblem shows that $\hat{s}$  satisfies this optimality condition corresponding to $\M$. Since cubic subproblem is strictly convex, it has unique minimizer and $\hat{s} = s(\M)$. The same argument applies to any convergent subsequence of $\{s(\M_j)\}$. Thus, $s(\M_j) \to s(\M)$ since ${s(\M_j)}$ is bounded sequence with only one accumulation point. Therefore $w(\M_j)\to w(\M)$ and thus $w(\M)$ is continuous. Since $L(\cdot)$ is continuous, the mappings $\norm{s(\M)}$ and $\M-L(w(\M))$ are continuous as well.

The estimate~\eqref{eq:s_bound} gives
$s(\M)\to0$ as $\M\to+\infty$, while
$\alpha(\M)\to0$ by~\eqref{eq:alpha_M_equality} and $v(\M)\to w_k$ by~\eqref{eq:interpolation_point_M}. Hence $w(\M)\to w_k$, so $L(w(\M))$ remains bounded. This gives us~\eqref{eq:M_L_limit}.
\end{proof}

\begin{lemma}
\label{thm:explicit_regularization}
Suppose that~\eqref{as:gh} and~\eqref{eq:current_invariant} hold. Then the
parameter
\begin{equation*}
  \M_k=\max\lb
    2M_0,
    12\sqrt{2}M_1^{3/2}\bar R^3\Gamma_k
  \rb
\end{equation*}
used in Algorithm~\ref{alg:global_accelerated_cubic} satisfies
\begin{equation}
  \M_k\geq L(w_{k+1})
  =M_0+M_1\norm{\grad f(w_{k+1})}.
  \label{eq:explicit_regularization_validity}
\end{equation}
Consequently, estimates~\eqref{eq:strengthened_lower}
and~\eqref{eq:invariant_update} hold for the actual iteration.
\end{lemma}

\begin{proof}
Define
\begin{equation*}
  \mathcal U_k\eqdef\lb \M>0:\M<L(w(\M))\rb.
\end{equation*}
If $\mathcal U_k=\emptyset$, the statement is immediate. Otherwise, let $\M_k^\dagger\eqdef\sup\mathcal U_k$. Lemma~\ref{lem:virtual_step_continuity} shows that $\M_k^\dagger<+\infty$.
By continuity,
\begin{equation}
  \M_k^\dagger=L(w(\M_k^\dagger)),
  \label{eq:boundary_admissibility}
\end{equation}
and every $\M\geq \M_k^\dagger$ is admissible.

Set
\begin{equation*}
  w^\dagger\eqdef w(\M_k^\dagger),
  \qquad
  g^\dagger\eqdef\grad f(w^\dagger),
  \qquad
  A_+^\dagger\eqdef A_+(\M_k^\dagger).
\end{equation*}
Applying~\eqref{eq:virtual_safe_step} at the admissible boundary and evaluating
$\psi_+(\cdot;\M_k^\dagger)$ at $x^\ast$, we obtain
\begin{gather*}
  A_+^\dagger f(w^\dagger)
  +\frac{A_+^\dagger}{9}
  \frac{\norm{g^\dagger}^{3/2}}{\sqrt{\M_k^\dagger}}
  \leq\psi_+^\ast(\M_k^\dagger)
  \leq\psi_+(x^\ast;\M_k^\dagger)
  \leq A_+^\dagger f^\ast
  +\Delta_0+\frac{\Gamma_0R^3}{3}.
\end{gather*}
Since $f(w^\dagger)\geq f^\ast$, $A_+^\dagger\geq A_k$, and
$\Gamma_k=\Gamma_0/A_k$, equation~\eqref{eq:upper_estimate_exact} gives
\begin{equation}
  \frac{\norm{g^\dagger}^{3/2}}{\sqrt{\M_k^\dagger}}
  \leq6\Gamma_k\bar R^3.
  \label{eq:boundary_gradient_bound}
\end{equation}
For every $g\in\R^d$,
\begin{equation}
  M_0+M_1\norm{g}
  \leq\max\lb
    2M_0,
    2\sqrt{2}M_1^{3/2}
    \frac{\norm{g}^{3/2}}{\sqrt{M_0+M_1\norm{g}}}
  \rb.
  \label{eq:L_max_bound}
\end{equation}
Indeed, the first term applies when $M_1\norm{g}\leq M_0$. Otherwise,
$M_0+M_1\norm{g}\leq2M_1\norm{g}$ and
\begin{equation*}
  \frac{\norm{g}^{3/2}}{\sqrt{M_0+M_1\norm{g}}}
  \geq\frac{\norm{g}}{\sqrt{2M_1}}.
\end{equation*}
Applying this bound with $g=g^\dagger$ and using
\eqref{eq:boundary_admissibility} and~\eqref{eq:boundary_gradient_bound},
we obtain
\begin{equation*}
  \M_k^\dagger
  \leq\max\lb
    2M_0,
    2\sqrt{2}M_1^{3/2}
    \frac{\norm{g^\dagger}^{3/2}}{\sqrt{\M_k^\dagger}}
  \rb 
  \leq\max\lb
    2M_0,
    12\sqrt{2}M_1^{3/2}\bar R^3\Gamma_k
  \rb
  =\M_k.
\end{equation*}
Therefore $\M_k$ is admissible and~\eqref{eq:explicit_regularization_validity}
holds. The last statement follows from Lemma~\ref{lem:one_safe_step}.
\end{proof}

\subsubsection{Estimating Sequence Bounds}

\begin{lemma}
\label{lem:global_invariant}
Suppose that~\eqref{as:gh} holds. For every $k\geq0$, the parameter $\M_k$ in
Algorithm~\ref{alg:global_accelerated_cubic} satisfies
\eqref{eq:explicit_regularization_validity}, and
\begin{equation}
  \psi_k^\ast\geq A_kf(w_k).
  \label{eq:global_invariant}
\end{equation}
Moreover,
\begin{equation}
  \psi_{k+1}^\ast
  \geq A_{k+1}f(w_{k+1})
  +\frac{A_{k+1}}{9}
  \frac{\norm{\grad f(w_{k+1})}^{3/2}}{\sqrt{\M_k}}.
  \label{eq:global_strengthened}
\end{equation}
\end{lemma}

\begin{proof}
For $k=0$,
\begin{equation*}
  \psi_0^\ast
  =\min_{x\in\R^d}\lb
    f(x_0)+\frac{\Gamma_0}{3}\norm{x-x_0}^3
  \rb
  =f(x_0)
  =A_0f(w_0).
\end{equation*}
Suppose that~\eqref{eq:global_invariant} holds at iteration $k$.
Theorem~\ref{thm:explicit_regularization} gives
\eqref{eq:explicit_regularization_validity}. Hence,
Lemma~\ref{lem:one_safe_step} gives~\eqref{eq:global_strengthened} and
\begin{equation*}
  \psi_{k+1}^\ast\geq A_{k+1}f(w_{k+1}),
\end{equation*}
which closes the induction.
\end{proof}

\begin{lemma}
\label{lem:global_error_bounds}
Suppose that~\eqref{as:gh} holds. For every $k\geq0$,
\begin{equation}
  f(w_k)-f^\ast
  \leq\frac{\Delta_0+\Gamma_0R^3/3}{A_k}
  \leq\frac{2}{3}\Gamma_k\bar R^3.
  \label{eq:global_error_bound}
\end{equation}
Moreover, every cubic output satisfies
\begin{equation}
  f(w_{k+1})-f^\ast
  +\frac{1}{9}
  \frac{\norm{\grad f(w_{k+1})}^{3/2}}{\sqrt{\M_k}}
  \leq\frac{\Delta_0+\Gamma_0R^3/3}{A_{k+1}}
  \leq\frac{2}{3}\Gamma_{k+1}\bar R^3.
  \label{eq:global_output_bound}
\end{equation}
\end{lemma}

\begin{proof}
Combining~\eqref{eq:global_invariant} and~\eqref{eq:upper_estimate},
\begin{gather*}
  A_kf(w_k)
  \stackrel{\eqref{eq:global_invariant}}{\leq}\psi_k^\ast
  \leq\psi_k(x^\ast)
  \stackrel{\eqref{eq:upper_estimate}}{\leq}
  A_kf^\ast+\Delta_0+\frac{\Gamma_0R^3}{3}.
\end{gather*}
Therefore,
\begin{gather*}
  f(w_k)-f^\ast
  \leq\frac{\Delta_0+\Gamma_0R^3/3}{A_k}
  \stackrel{\eqref{eq:upper_estimate_exact}}{\leq}
  \frac{2\Gamma_0\bar R^3}{3A_k}
  \stackrel{\eqref{eq:equivalent_weights}}{=}
  \frac{2}{3}\Gamma_k\bar R^3.
\end{gather*}
Similarly,
\begin{equation*}
  A_{k+1}f(w_{k+1})
  +\frac{A_{k+1}}{9}
  \frac{\norm{\grad f(w_{k+1})}^{3/2}}{\sqrt{\M_k}}
  \stackrel{\eqref{eq:global_strengthened}}{\leq}\psi_{k+1}^\ast
  \leq\psi_{k+1}(x^\ast)
  \stackrel{\eqref{eq:upper_estimate}}{\leq}
  A_{k+1}f^\ast+\Delta_0+\frac{\Gamma_0R^3}{3}.
\end{equation*}
Dividing by $A_{k+1}$ and using~\eqref{eq:upper_estimate_exact} and
\eqref{eq:equivalent_weights} proves~\eqref{eq:global_output_bound}.
\end{proof}

\subsubsection{Convergence Rates}

Since $\alpha_k\in(0,1)$ and $\Gamma_0>0$,
\begin{equation*}
  0<\Gamma_{k+1}=(1-\alpha_k)\Gamma_k<\Gamma_k.
\end{equation*}
Thus, $\{\Gamma_k\}_{k\geq0}$ is strictly decreasing. Define the transition
level
\begin{equation}
  \Gamma_{\tr}\eqdef
  \begin{cases}
    \dfrac{M_0}{6\sqrt{2}M_1^{3/2}\bar R^3}, & M_1>0,\\
    +\infty, & M_1=0.
  \end{cases}
  \label{eq:gamma_transition}
\end{equation}

\begin{lemma}[Far regime]
\label{lem:far_regime}
Suppose that $M_1>0$ and $\Gamma_k>\Gamma_{\tr}$. Then
\begin{equation}
  \alpha_k
  \geq\frac{1}{12\sqrt[6]{2}\ls1+\sqrt{M_1}\bar R\rs},
  \label{eq:far_weight_bound}
\end{equation}
and
\begin{equation}
  \Gamma_{k+1}
  \leq\exp\ls
    -\frac{1}{12\sqrt[6]{2}\ls1+\sqrt{M_1}\bar R\rs}
  \rs\Gamma_k.
  \label{eq:far_contraction}
\end{equation}
\end{lemma}

\begin{proof}
In this regime,
\begin{equation*}
  \M_k=12\sqrt{2}M_1^{3/2}\bar R^3\Gamma_k,
\end{equation*}
and hence
\begin{equation*}
  \frac{\alpha_k^3}{1-\alpha_k}
  =\frac{1}{864\sqrt{2}M_1^{3/2}\bar R^3}.
\end{equation*}
If $\alpha_k>1/2$, then~\eqref{eq:far_weight_bound} holds immediately.
If $\alpha_k\leq1/2$, then
\begin{equation*}
  \alpha_k^3
  =\tfrac{1-\alpha_k}{864\sqrt{2}M_1^{3/2}\bar R^3}
  \geq\tfrac{1}{1728\sqrt{2}M_1^{3/2}\bar R^3}, \quad
  \alpha_k
  \geq\tfrac{1}{12\sqrt[6]{2}\sqrt{M_1}\bar R}.
\end{equation*}
Combining the two cases proves~\eqref{eq:far_weight_bound}. Finally,
$1-\alpha_k\leq\exp(-\alpha_k)$ proves~\eqref{eq:far_contraction}.
\end{proof}

\begin{lemma}[Near regime]
\label{lem:near_regime}
Suppose that~\eqref{as:gh} holds and
$\Gamma_K\leq\Gamma_{\tr}$. Then $\M_j=2M_0$ for every
$j\geq K$, and
\begin{equation}
  A_{K+n}^{1/3}
  \geq A_K^{1/3}
  +\frac{n}{3}\ls\frac{\Gamma_0}{144M_0}\rs^{1/3},
  \qquad n\geq1.
  \label{eq:near_growth}
\end{equation}
Consequently,
\begin{equation}
  \Gamma_{K+n}\leq\frac{3888M_0}{n^3},
  \qquad n\geq1,
  \label{eq:near_gamma_bound}
\end{equation}
and
\begin{equation}
  f(w_{K+n})-f^\ast
  \leq\frac{2592M_0\bar R^3}{n^3},
  \qquad n\geq1.
  \label{eq:near_function_bound}
\end{equation}
\end{lemma}

\begin{proof}
The monotonicity of $\Gamma_k$ and~\eqref{eq:gamma_transition} give
$\M_j=2M_0$ for every $j\geq K$. Therefore,
\eqref{eq:equivalent_weights} gives
\begin{equation*}
  a_{j+1}^3=\frac{\Gamma_0A_{j+1}^2}{144M_0}.
\end{equation*}
Using $A_{j+1}-A_j=a_{j+1}$,
\begin{equation*}
  A_{j+1}^{1/3}-A_j^{1/3}
  =\frac{a_{j+1}}
  {A_{j+1}^{2/3}+A_{j+1}^{1/3}A_j^{1/3}+A_j^{2/3}}
  \geq\frac{a_{j+1}}{3A_{j+1}^{2/3}}
  =\frac{1}{3}\ls\frac{\Gamma_0}{144M_0}\rs^{1/3}.
\end{equation*}
Summing over $j=K,\ldots,K+n-1$ proves~\eqref{eq:near_growth}. Hence,
\begin{equation*}
  A_{K+n}
  \geq\frac{n^3\Gamma_0}{3888M_0}, \quad
  \Gamma_{K+n}
  \stackrel{\eqref{eq:equivalent_weights}}{=}
  \frac{\Gamma_0}{A_{K+n}}
  \leq\frac{3888M_0}{n^3},
\end{equation*}
which proves~\eqref{eq:near_gamma_bound}. Combining
\eqref{eq:global_error_bound} and~\eqref{eq:near_gamma_bound} proves
\eqref{eq:near_function_bound}.
\end{proof}

\begin{theorem}
\label{thm:explicit_global_complexity}
Suppose that $f$ is convex,~\eqref{as:gh} and~\eqref{eq:inital_cnd} hold. For every
$\e>0$, Algorithm~\ref{alg:global_accelerated_cubic} returns an iterate
$w_{N}$ satisfying
\begin{equation*}
  f(w_{N})-f^\ast\leq\e
\end{equation*}
after at most
\begin{equation}
\begin{aligned}
  N \leq \left\lceil
    12\sqrt[6]{2}\ls1+\sqrt{M_1}\bar R\rs
    \log\max\lb
      1,
      \min\lb
        \tfrac{6\sqrt{2}M_1^{3/2}\bar R^3\Gamma_0}{M_0},
        \tfrac{2\Gamma_0\bar R^3}{3\e}
      \rb
    \rb
  \right\rceil + 
  \left\lceil
    \ls\tfrac{2592M_0\bar R^3}{\e}\rs^{1/3}
  \right\rceil
\end{aligned}
\label{eq:global_complexity}
\end{equation}
iterations.
\end{theorem}

\begin{proof}
Define
\begin{equation*}
  \Gamma_\e\eqdef\frac{3\e}{2\bar R^3}.
\end{equation*}
By~\eqref{eq:global_error_bound}, if $\Gamma_k\leq\Gamma_\e$, then
\begin{equation*}
  f(w_k)-f^\ast
  \leq\frac{2}{3}\Gamma_k\bar R^3
  \leq\frac{2}{3}\Gamma_\e\bar R^3
  =\e.
\end{equation*}
Suppose first that $M_1>0$, and set
\begin{equation*}
\begin{aligned}
  N_{\rm far}\eqdef
  \left\lceil
    12\sqrt[6]{2}\ls1+\sqrt{M_1}\bar R\rs
    \log\max\lb
      1,
      \frac{\Gamma_0}{\max\lb\Gamma_{\tr},\Gamma_\e\rb}
    \rb
  \right\rceil.
\end{aligned}
\end{equation*}
Suppose that
$\Gamma_k>\max\lb\Gamma_{\tr},\Gamma_\e\rb$ for every
$k=0,\ldots,N_{\rm far}-1$. Then~\eqref{eq:far_contraction} applies at
each of these iterations, and
\begin{equation*}
  \Gamma_{N_{\rm far}}
  \leq\exp\ls
    -\frac{N_{\rm far}}
    {12\sqrt[6]{2}\ls1+\sqrt{M_1}\bar R\rs}
  \rs\Gamma_0
  \leq\max\lb\Gamma_{\tr},\Gamma_\e\rb.
\end{equation*}
Therefore, after at most $N_{\rm far}$ iterations, either
$\Gamma_k\leq\Gamma_\e$, or $\Gamma_k\leq\Gamma_{\tr}$ and the
near regime begins. Moreover,
\begin{gather*}
  \frac{\Gamma_0}{\Gamma_{\tr}}
  =\frac{6\sqrt{2}M_1^{3/2}\bar R^3\Gamma_0}{M_0},
  \qquad
  \frac{\Gamma_0}{\Gamma_\e}
  =\frac{2\Gamma_0\bar R^3}{3\e},\\
  \frac{\Gamma_0}
  {\max\lb\Gamma_{\tr},\Gamma_\e\rb}
  =\min\lb
    \frac{6\sqrt{2}M_1^{3/2}\bar R^3\Gamma_0}{M_0},
    \frac{2\Gamma_0\bar R^3}{3\e}
  \rb.
\end{gather*}
If the near regime is needed, Lemma~\ref{lem:near_regime} shows that
\begin{equation*}
  N_{\rm near}
  \eqdef\left\lceil
    \ls\frac{2592M_0\bar R^3}{\e}\rs^{1/3}
  \right\rceil
\end{equation*}
additional iterations suffice. Hence,
$N \leq N_{\rm far}+N_{\rm near}$, which
proves~\eqref{eq:global_complexity}.

If $M_1=0$, then $\Gamma_{\tr}=+\infty$ and the near regime starts at
$K=0$. The first term in~\eqref{eq:global_complexity} is zero, while
Lemma~\ref{lem:near_regime} gives the second term.
\end{proof}

\subsection{Adaptive Method under~\eqref{as:lh} with Exact Hessian}
\label{sec:fully_adaptive}

Algorithm~\ref{alg:global_accelerated_cubic} uses the problem parameters in its
regularization coefficient and requires an initial scale
satisfying~\eqref{eq:inital_cnd}. We now give one implementation that does not use
$M_0$, $M_1$, $R$, $\rho$, or $f^\ast$. 

\textbf{Initialization step.} Our first goal is to construct a computable initial scale $\Gamma_0$ satisfying 
\begin{equation}
  \label{eq:inital_cnd2}
  \Gamma_0 R^3\geq 3\ls f(x_0)-f^\ast\rs.
\end{equation}
without knowing $R$ and $f^\ast$. 
We construct them by solving the cubic subproblem~\eqref{eq:cubic_initial_step} with backtracking regularization parameter $\tilde{\M}$ until condition~\eqref{eq:cubic_initial_test} holds.
The whole warm-up backtracking procedure is listed in lines~\ref{line:adaptive_initial_start}-\ref{line:adaptive_initial_end} of Algorithm~\ref{alg:fully_adaptive_cubic}. The following lemma shows that the backtracking procedure terminates after a finite number of iterations and provides an exact way of satisfying initial condition~\eqref{eq:inital_cnd2}.

\begin{lemma}
\label{lem:cubic_initialization}
Suppose that either~\eqref{as:gh} or~\eqref{as:lh} holds. If
$\nabla f(x_0)=0$, Algorithm~\ref{alg:fully_adaptive_cubic} returns the
optimal point $x_0$. Otherwise, the initialization
backtracking (lines~\ref{line:adaptive_initial_start}-\ref{line:adaptive_initial_end} of Algorithm~\ref{alg:fully_adaptive_cubic}) satisfies~\eqref{eq:cubic_initial_test}  after at most $\left\lceil
    \log_2\max\lb
      1, \tfrac{T_{0}}{\tilde{\M}_0} \rb \right\rceil$ rejected trials, where 
\begin{equation}
  T_{0} = 
  \begin{cases}
    M_0+2M_1\norm{\nabla f(x_0)},
    & \text{under~\eqref{as:gh}},\\
    \max\lb
      M_0+2M_1\norm{\nabla f(x_0)},
      \tfrac{2\norm{\nabla f(x_0)}}{15\rho^2}
    \rb,
    & \text{under~\eqref{as:lh}}.
  \end{cases}
  \label{eq:initial_safe_threshold}
\end{equation}
For the accepted step,
$\Gamma_0 \eqdef \tfrac{3 \|\nabla f(x_0)\|^3}{(f(x_0) - f(\tilde{x}))^2}$
is well defined and satisfies initial condition~\eqref{eq:inital_cnd2}.
\end{lemma}

\begin{algorithm}[!htb]
  \caption{Fully adaptive accelerated cubic regularization}
  \label{alg:fully_adaptive_cubic}
  \begin{algorithmic}[1]
    \STATE \textbf{Input:} $x_0\in\R^d$ and $\tilde{\M}_0>0$.
    \STATE If $\nabla f(x_0)=0$, return $x_0$.
    \STATE Set $\tilde{\M}\eqdef\tilde{\M}_0$.
    \STATE\label{line:adaptive_initial_start} Compute
      \begin{equation}
        \tilde{s}
        \eqdef
        \argmin_{s\in\R^d}
        \lb
          \la \nabla f(x_0), s\ra
          +\tfrac{1}{2}\la\nabla^2f(x_0)s,s\ra
          +\tfrac{15\tilde{\M}}{6}\norm{s}^3
        \rb,
        \quad
        \tilde{x}\eqdef x_0+ \tilde{s}.
        \label{eq:cubic_initial_step}
      \end{equation}
    \WHILE{condition
        \begin{equation}
          \la \nabla f(\tilde x),-\tilde s\ra
          \geq
          \tfrac{2}{9\sqrt{\tilde{\M}}}
          \norm{\nabla f(\tilde{x})}^{3/2}.
          \label{eq:cubic_initial_test}
        \end{equation}
      does not hold}
      \STATE Set $\tilde{\M}:=2\tilde{\M}$.
      \STATE Recompute $\tilde{s}$ and $\tilde{x}$ according
        to~\eqref{eq:cubic_initial_step}.\label{line:adaptive_initial_end}
    \ENDWHILE
    \STATE Set $\Gamma_0 \eqdef \tfrac{3 \|\nabla f(x_0)\|^3}{(f(x_0) - f(\tilde{x}))^2}$, $A_0\eqdef1$, $w_0\eqdef x_0$, $z_0\eqdef x_0$, and
      $\psi_0(x)\eqdef f(x_0)+\tfrac{\Gamma_0}{3}\norm{x-x_0}^3$.
    \STATE Set $\M_0\eqdef\Gamma_0$.
    \FOR{$k \geq 0$}
      \STATE Find the unique $\alpha_k\in(0,1)$ satisfying
        $\tfrac{\alpha_k^3}{1-\alpha_k}
        =\tfrac{\Gamma_k}{72\M_k}$.
      \STATE Set $A_{k+1}$, $a_{k+1}$, and $\Gamma_{k+1}$ according
        to~\eqref{eq:a_params}.
      \STATE Set $v_k\eqdef(1-\alpha_k)w_k+\alpha_kz_k$.
      \STATE Compute $s_k$ and $w_{k+1}$ according
        to~\eqref{eq:cubic_step} with regularization coefficient $\M_k$.
      \WHILE{condition
        \begin{equation}
          \label{eq:adaptive_certificate}
          \la \nabla f(w_{k+1}), v_k - w_{k+1}\ra
          \geq
          \tfrac{2}{9\sqrt{\M_k}}
          \norm{\nabla f(w_{k+1})}^{3/2}.
        \end{equation}
      does not hold}
        \STATE Set $\M_k:=2\M_k$.
        \STATE Recompute $\alpha_k$, $A_{k+1}$, $a_{k+1}$,
          $\Gamma_{k+1}$, and $v_k$ as above.
        \STATE Recompute $w_{k+1}$ according
          to~\eqref{eq:cubic_step} with regularization coefficient $\M_k$.
      \ENDWHILE
      \STATE Set $\psi_{k+1}(x)\eqdef
        \psi_k(x)+a_{k+1}\ell_{w_{k+1}}(x)$ and compute
        $z_{k+1}\eqdef\argmin_{x\in\R^d}\psi_{k+1}(x)$.
      \STATE Set
        $\M_{k+1}
        :=\tfrac{\Gamma_{k+1}}{\Gamma_k}\M_k$.
    \ENDFOR
\end{algorithmic}
\end{algorithm}

For the statement below, define
\begin{equation}
  B\eqdef
  \begin{cases}
    12\sqrt{2}M_1^{3/2}R^3,
    & \text{under~\eqref{as:gh}},\\
    \max\lb
      12\sqrt{2}M_1^{3/2}R^3,
      \dfrac{2R^3}{9\rho^3}
    \rb,
    & \text{under~\eqref{as:lh}}.
  \end{cases}
  \label{eq:adaptive_B}
\end{equation}

\begin{theorem}
\label{thm:fully_adaptive_convergence}
Suppose that $f$ is convex and either~\eqref{as:gh} or~\eqref{as:lh} holds.
All backtracking loops in Algorithm~\ref{alg:fully_adaptive_cubic} terminate.
For every $\e>0$, under~\eqref{as:gh}, the method returns an iterate $w_N$
satisfying $f(w_N)-f^\ast\leq\e$ for
\begin{equation}
  N
  =\cO\ls
    \ls\tfrac{M_0R^3}{\e}\rs^{1/3}
  \rs
  +\widetilde{\cO}\ls
    1+\sqrt{M_1}R
  \rs.
  \label{eq:adaptive_global_rate}
\end{equation}
Under~\eqref{as:lh}, the corresponding bound is
\begin{equation}
  N
  =\cO\ls
    \ls\tfrac{M_0R^3}{\e}\rs^{1/3}
  \rs
  +\widetilde{\cO}\ls
    1+\sqrt{M_1}R+\tfrac{R}{\rho}
  \rs.
  \label{eq:adaptive_local_rate}
\end{equation}
The number of rejected initialization trials is bounded as in
Lemma~\ref{lem:cubic_initialization}. Before returning $w_N$, the main
backtracking loops reject at most
$
  \left\lceil
    \log_2\max\lb
      1,
      B,
      \tfrac{4M_0R^3}{3\e}
    \rb
  \right\rceil
$
trials in total, where $B$ is defined in~\eqref{eq:adaptive_B}.
\end{theorem}

\subsubsection{Proof of Lemma~\ref{lem:cubic_initialization}}

\begin{proof}
Suppose that $\nabla f(x_0)\neq0$ and fix a trial value $\tilde{\M}>0$.
The optimality condition for~\eqref{eq:cubic_initial_step} is
\begin{equation}
  \nabla f(x_0)+\nabla^2f(x_0)\tilde s
  +\frac{15\tilde{\M}}{2}\norm{\tilde s}\tilde s=0.
  \label{eq:initial_cubic_optimality}
\end{equation}
Since $f$ is convex, $\nabla^2f(x_0)\succeq0$. Taking the inner product
of~\eqref{eq:initial_cubic_optimality} with $\tilde s$, we obtain
\begin{gather*}
  \frac{15\tilde{\M}}{2}\norm{\tilde s}^3
  \leq-\la\nabla f(x_0),\tilde s\ra
  \leq\norm{\nabla f(x_0)}\norm{\tilde s},
\end{gather*}
and hence
\begin{equation}
  \norm{\tilde s}^2
  \leq\frac{2\norm{\nabla f(x_0)}}{15\tilde{\M}}.
  \label{eq:initial_step_bound}
\end{equation}
Let $\tilde{\M}\geq T_0$, defined in~\eqref{eq:initial_safe_threshold}. Under~\eqref{as:lh},
equations~\eqref{eq:initial_safe_threshold}
and~\eqref{eq:initial_step_bound} give
\begin{equation*}
  \norm{\tilde s}^2
  \leq\frac{2\norm{\nabla f(x_0)}}{15\tilde{\M}}
  \leq\rho^2.
\end{equation*}
Therefore, the Taylor remainder at $x_0$ can be applied under both
assumptions. Using~\eqref{eq:initial_cubic_optimality}, we get
\begin{gather*}
  \norm{\nabla f(\tilde x)}
  \leq
  \norm{\nabla f(\tilde x)-\nabla f(x_0)
  -\nabla^2f(x_0)\tilde s}
  +\frac{15\tilde{\M}}{2}\norm{\tilde s}^2
  \\
  \leq
  \frac{M_0+M_1\norm{\nabla f(x_0)}+15\tilde{\M}}{2}
  \norm{\tilde s}^2
  \\
  \leq
  \ls
    1+\frac{M_0+M_1\norm{\nabla f(x_0)}}{15\tilde{\M}}
  \rs\norm{\nabla f(x_0)}
  \leq\frac{16}{15}\norm{\nabla f(x_0)}.
\end{gather*}
Consequently,
\begin{gather*}
  M_0+M_1\norm{\nabla f(\tilde x)}
  \leq M_0+2M_1\norm{\nabla f(x_0)}
  \leq \tilde{\M}.
\end{gather*}
Since $\norm{\tilde s}\leq\rho$ under~\eqref{as:lh}, the hypotheses of Lemma~\ref{lem:backward_progress} hold under both assumptions. Therefore,~\eqref{eq:progress_final} gives
\begin{equation*}
  \la\nabla f(\tilde x),-\tilde s\ra
  \geq\frac{2}{9\sqrt{\tilde{\M}}}\norm{\nabla f(\tilde x)}^{3/2},
\end{equation*}
which is condition~\eqref{eq:cubic_initial_test}. Hence, every trial with
$\tilde{\M}\geq T_0$ is accepted.

Starting from the input value $\tilde{\M}_0$, the $j$-th trial uses
$2^j\tilde{\M}_0$. For
\begin{equation*}
  j=\left\lceil
    \log_2\max\lb1,\frac{T_0}{\tilde{\M}_0}\rb
  \right\rceil,
\end{equation*}
we have $2^j\tilde{\M}_0\geq T_0$. Therefore, the backtracking rejects at most $j$
trials.

It remains to verify~\eqref{eq:inital_cnd2}. If
$\nabla f(\tilde x)\neq0$, condition~\eqref{eq:cubic_initial_test} and
convexity give
\begin{equation*}
  f(x_0)-f(\tilde x)
  \geq\la\nabla f(\tilde x),x_0-\tilde x\ra
  >0.
\end{equation*}
If $\nabla f(\tilde x)=0$, then $\tilde x$ is optimal. Since
$\nabla f(x_0)\neq0$, in this case we also have
$f(x_0)-f(\tilde x)>0$. Moreover, by convexity,
\begin{gather*}
  0< f(x_0)-f(\tilde x)
  \leq f(x_0)-f^\ast
  \leq\la\nabla f(x_0),x_0-x^\ast\ra
  \leq\norm{\nabla f(x_0)}R.
\end{gather*}
Thus,
\begin{equation*}
  0<\frac{f(x_0)-f(\tilde x)}{\norm{\nabla f(x_0)}}\leq R.
\end{equation*}
Using this estimate one more time, we obtain
\begin{gather*}
  f(x_0)-f^\ast
  \leq\norm{\nabla f(x_0)}R
  =\ls f(x_0)-f(\tilde x)\rs
  \frac{\norm{\nabla f(x_0)}R}{f(x_0)-f(\tilde x)}
  \\
  \leq\ls f(x_0)-f(\tilde x)\rs
  \ls
    \frac{\norm{\nabla f(x_0)}R}{f(x_0)-f(\tilde x)}
  \rs^3
  =\frac{\Gamma_0R^3}{3}.
\end{gather*}
This proves~\eqref{eq:inital_cnd2}.
\end{proof}

\subsubsection{Proof of Theorem~\ref{thm:fully_adaptive_convergence}}

Fix an iteration $k$. For a trial coefficient $\M>0$, Algorithm~\ref{alg:fully_adaptive_cubic} computes $\alpha_k$, $A_{k+1}$, $a_{k+1}$, $v_k$, and $w_{k+1}$ in the same way as Algorithm~\ref{alg:global_accelerated_cubic} with $\M_k=\M$. Hence, under the induction assumption~\eqref{eq:current_invariant}, Lemma~\ref{lem:one_safe_step} gives~\eqref{eq:raw_one_step} for every trial. If the trial passes~\eqref{eq:adaptive_certificate}, which is~\eqref{eq:progress_final} for $v=v_k$ and $\M=\M_k$, then the proof of Lemma~\ref{lem:one_safe_step} gives~\eqref{eq:strengthened_lower} and~\eqref{eq:invariant_update}.

We now prove that~\eqref{eq:adaptive_certificate} holds for every sufficiently large trial coefficient. Let $B$ be defined
by~\eqref{eq:adaptive_B}.

\begin{lemma}
\label{lem:adaptive_safe_trial}
Suppose that either~\eqref{as:gh} or~\eqref{as:lh} holds and that
\eqref{eq:current_invariant} is satisfied. Every trial with
$\M\geq\max\lb2M_0,B\Gamma_k\rb$
satisfies~\eqref{eq:adaptive_certificate}.
\end{lemma}

\begin{proof}
By~\eqref{eq:inital_cnd2}, condition~\eqref{eq:inital_cnd} holds with
$\bar R=R$. Under~\eqref{as:gh}, we set $\rho=+\infty$. The trial with coefficient $\M$ coincides with the virtual step $w(\M)$ defined after~\eqref{eq:interpolation_point_M}. Set $r(\M)\eqdef\norm{w(\M)-v(\M)}$ and
\begin{equation*}
  \mathcal U_k\eqdef\lb\M>0:r(\M)>\rho\text{ or }\M<L(w(\M))\rb.
\end{equation*}
For $\M\notin\mathcal U_k$, the hypotheses of Lemma~\ref{lem:backward_progress} hold for $v=v(\M)$, and~\eqref{eq:progress_final} gives~\eqref{eq:adaptive_certificate}. Thus, it suffices to show that every $\M\geq\max\lb2M_0,B\Gamma_k\rb$ satisfies $\M\notin\mathcal U_k$.

If $\mathcal U_k=\emptyset$, there is nothing to prove. Otherwise, Lemma~\ref{lem:virtual_step_continuity} shows that $\mathcal U_k$ is open and bounded above. Let $\M_k^\dagger\eqdef\sup\mathcal U_k<+\infty$. Then $\M\notin\mathcal U_k$ for every $\M\geq\M_k^\dagger$. By continuity,
\begin{equation}
  r^\dagger\leq\rho,
  \qquad
  \M_k^\dagger\geq L(w^\dagger),
  \qquad
  r^\dagger=\rho
  \text{ or }
  \M_k^\dagger=L(w^\dagger),
  \label{eq:adaptive_boundary}
\end{equation}
where a dagger denotes evaluation at $\M_k^\dagger$ and $g^\dagger\eqdef\grad f(w^\dagger)$.

By~\eqref{eq:adaptive_boundary}, the hypotheses of Lemma~\ref{lem:backward_progress} hold at $\M_k^\dagger$. Hence,~\eqref{eq:q_step},~\eqref{eq:inner_step}, and~\eqref{eq:progress_final} give
\begin{gather*}
  \la g^\dagger,v^\dagger-w^\dagger\ra
  -\frac{\norm{g^\dagger}^{3/2}}{9\sqrt{\M_k^\dagger}}
  \geq
  6\M_k^\dagger(r^\dagger)^3
  -\frac{\ls9\M_k^\dagger(r^\dagger)^2\rs^{3/2}}{9\sqrt{\M_k^\dagger}}
  =3\M_k^\dagger(r^\dagger)^3,\\
  \la g^\dagger,v^\dagger-w^\dagger\ra
  -\frac{\norm{g^\dagger}^{3/2}}{9\sqrt{\M_k^\dagger}}
  \geq
  \frac{\norm{g^\dagger}^{3/2}}{9\sqrt{\M_k^\dagger}}.
\end{gather*}
Lemma~\ref{lem:one_safe_step} applies to the virtual step, as explained before~\eqref{eq:virtual_safe_step}. Therefore,~\eqref{eq:raw_one_step} and evaluation of $\psi_+(\cdot;\M_k^\dagger)$ at $x^\ast$ give
\begin{gather*}
  A_+^\dagger f(w^\dagger)
  +A_+^\dagger\max\lb
    3\M_k^\dagger(r^\dagger)^3,
    \frac{\norm{g^\dagger}^{3/2}}{9\sqrt{\M_k^\dagger}}
  \rb
  \leq\psi_+^\ast(\M_k^\dagger)
  \leq\psi_+(x^\ast;\M_k^\dagger)
  \leq A_+^\dagger f^\ast
  +\Delta_0+\frac{\Gamma_0R^3}{3}.
\end{gather*}
Since $f(w^\dagger)\geq f^\ast$, $A_+^\dagger\geq A_k$, and
$\Gamma_k=\Gamma_0/A_k$, equation~\eqref{eq:upper_estimate_exact} with $\bar R=R$ gives
\begin{equation}
  \max\lb
    3\M_k^\dagger(r^\dagger)^3,
    \frac{\norm{g^\dagger}^{3/2}}{9\sqrt{\M_k^\dagger}}
  \rb
  \leq\frac{2}{3}\Gamma_kR^3.
  \label{eq:adaptive_boundary_budget}
\end{equation}
If $r^\dagger=\rho$, then~\eqref{eq:adaptive_boundary_budget} gives
\begin{equation*}
  \M_k^\dagger\leq\frac{2R^3}{9\rho^3}\Gamma_k\leq B\Gamma_k.
\end{equation*}
If $\M_k^\dagger=L(w^\dagger)$, then~\eqref{eq:L_max_bound} and~\eqref{eq:adaptive_boundary_budget} give
\begin{equation*}
  \M_k^\dagger
  \leq\max\lb
    2M_0,
    2\sqrt{2}M_1^{3/2}
    \frac{\norm{g^\dagger}^{3/2}}{\sqrt{\M_k^\dagger}}
  \rb
  \leq\max\lb
    2M_0,
    12\sqrt{2}M_1^{3/2}R^3\Gamma_k
  \rb
  \leq\max\lb2M_0,B\Gamma_k\rb.
\end{equation*}
In both cases, every $\M\geq\max\lb2M_0,B\Gamma_k\rb$ satisfies $\M\geq\M_k^\dagger$, and hence $\M\notin\mathcal U_k$.
\end{proof}

\begin{lemma}
\label{lem:adaptive_backtracking}
Suppose that either~\eqref{as:gh} or~\eqref{as:lh} holds. Every main
backtracking loop in Algorithm~\ref{alg:fully_adaptive_cubic} terminates.
For every accepted iteration $k\geq0$,
\begin{equation}
  \M_k
  \leq\max\lb
    4M_0,
    \max\lb1,2B\rb\Gamma_k
  \rb.
  \label{eq:adaptive_regularization_upper}
\end{equation}
For every $K\geq1$, the total number of rejected main trials during
iterations $k=0,\ldots,K-1$ is at most
\begin{equation}
  \left\lceil
    \log_2\max\lb
      1,
      B,
      \frac{2M_0}{\Gamma_{K-1}}
    \rb
  \right\rceil.
  \label{eq:adaptive_rejected_prefix}
\end{equation}
\end{lemma}

\begin{proof}
At $k=0$,~\eqref{eq:current_invariant} holds with equality. Suppose that it
holds at the beginning of iteration $k$. Repeated doubling sends $\M_k$ to
$+\infty$. Hence, Lemma~\ref{lem:adaptive_safe_trial} shows that
condition~\eqref{eq:adaptive_certificate} is satisfied after finitely many
trials. For the accepted trial,~\eqref{eq:invariant_update} gives
\eqref{eq:current_invariant} at iteration $k+1$. Thus, every backtracking loop
terminates by induction.

It remains to prove~\eqref{eq:adaptive_regularization_upper}. At $k=0$, if
the first trial is accepted, then
\begin{equation*}
  \M_0=\Gamma_0
  \leq\max\lb1,2B\rb\Gamma_0.
\end{equation*}
If at least one trial is rejected, the coefficient immediately preceding the
accepted one is $\M_0/2$. Lemma~\ref{lem:adaptive_safe_trial} gives
$
  \frac{\M_0}{2}
  <\max\lb2M_0,B\Gamma_0\rb
$.
Therefore,
$
  \M_0
  <\max\lb4M_0,2B\Gamma_0\rb
$,
and~\eqref{eq:adaptive_regularization_upper} follows.

Suppose that~\eqref{eq:adaptive_regularization_upper} holds at iteration
$k-1$. If the first trial is accepted at iteration $k$, then
\begin{gather*}
  \M_k
  =\frac{\Gamma_k}{\Gamma_{k-1}}\M_{k-1}
  \leq
  \frac{\Gamma_k}{\Gamma_{k-1}}
  \max\lb
    4M_0,
    \max\lb1,2B\rb\Gamma_{k-1}
  \rb
  \leq\max\lb
    4M_0,
    \max\lb1,2B\rb\Gamma_k
  \rb,
\end{gather*}
where we used $\Gamma_k\leq\Gamma_{k-1}$. If at least one trial is rejected,
the same argument as for $k=0$ proves the bound. This closes the induction.

Finally, the initialization and the update of $\M_k$ in
Algorithm~\ref{alg:fully_adaptive_cubic} imply that the ratio
$\M_k/\Gamma_k$ is initially equal to one, remains unchanged between
accepted iterations, and doubles after every rejected trial. Since
$\Gamma_k$ is decreasing, for every $k=0,\ldots,K-1$,
\begin{equation*}
  \max\lb
    B,
    \frac{2M_0}{\Gamma_k}
  \rb
  \leq\max\lb
    B,
    \frac{2M_0}{\Gamma_{K-1}}
  \rb.
\end{equation*}
By Lemma~\ref{lem:adaptive_safe_trial}, every trial satisfying
\begin{equation*}
  \frac{\M_k}{\Gamma_k}
  \geq\max\lb
    B,
    \frac{2M_0}{\Gamma_{K-1}}
  \rb
\end{equation*}
is accepted. After $J$ rejected main trials, the current ratio is $2^J$.
This proves~\eqref{eq:adaptive_rejected_prefix}.
\end{proof}

\begin{lemma}
\label{lem:adaptive_error_bound}
Suppose that $f$ is convex and either~\eqref{as:gh} or~\eqref{as:lh} holds.
For every accepted iterate of Algorithm~\ref{alg:fully_adaptive_cubic},
\begin{equation}
  f(w_k)-f^\ast
  \leq\frac{2}{3}\Gamma_kR^3.
  \label{eq:adaptive_error_bound}
\end{equation}
\end{lemma}

\begin{proof}
The induction in the proof of Lemma~\ref{lem:adaptive_backtracking} shows
that~\eqref{eq:current_invariant} holds for every accepted iteration. The
argument in the proof of Lemma~\ref{lem:global_error_bounds}, applied with
$\bar R=R$, proves~\eqref{eq:adaptive_error_bound}.
\end{proof}

For the scalar analysis, set
\begin{equation}
  \widehat B\eqdef\max\lb1,2B\rb,
  \qquad
  \Gamma_{\tr}^{\rm ad}\eqdef\frac{4M_0}{\widehat B}.
  \label{eq:adaptive_gamma_transition}
\end{equation}
The following two lemmas separate the convergence of the fully adaptive
method into the far and near regimes. This separation is used only in the
analysis and does not require the algorithm to know the transition level.

\begin{lemma}[Far regime]
\label{lem:adaptive_far_regime}
Suppose that $f$ is convex, either~\eqref{as:gh} or~\eqref{as:lh} holds, and
$\Gamma_k>\Gamma_{\tr}^{\rm ad}$. Then
\begin{equation}
  \alpha_k
  \geq
  \frac{1}{2+(144\widehat B)^{1/3}},
  \label{eq:adaptive_far_weight_bound}
\end{equation}
and
\begin{equation}
  \Gamma_{k+1}
  \leq
  \exp\ls
    -\frac{1}{2+(144\widehat B)^{1/3}}
  \rs\Gamma_k.
  \label{eq:adaptive_far_contraction}
\end{equation}
\end{lemma}

\begin{proof}
By~\eqref{eq:adaptive_gamma_transition},
$\widehat B\Gamma_k>4M_0$. Therefore,
equation~\eqref{eq:adaptive_regularization_upper} gives
\begin{equation*}
  \M_k\leq\widehat B\Gamma_k,
  \qquad
  \frac{\alpha_k^3}{1-\alpha_k}
  =\frac{\Gamma_k}{72\M_k}
  \geq\frac{1}{72\widehat B}.
\end{equation*}
If $\alpha_k\leq1/2$, then
\begin{equation*}
  \alpha_k^3\geq\frac{1}{144\widehat B}.
\end{equation*}
If $\alpha_k>1/2$, the bound $\alpha_k>1/2$ holds directly. Hence,
\begin{equation*}
  \alpha_k
  \geq
  \min\lb
    \frac12,
    \frac{1}{(144\widehat B)^{1/3}}
  \rb
  \geq
  \frac{1}{2+(144\widehat B)^{1/3}}.
\end{equation*}
This proves~\eqref{eq:adaptive_far_weight_bound}.
Equation~\eqref{eq:adaptive_far_contraction} follows from
$\Gamma_{k+1}=(1-\alpha_k)\Gamma_k$ and
$1-t\leq\exp(-t)$.
\end{proof}

\begin{lemma}[Near regime]
\label{lem:adaptive_near_regime}
Suppose that $f$ is convex, either~\eqref{as:gh} or~\eqref{as:lh} holds, and
$\Gamma_K\leq\Gamma_{\tr}^{\rm ad}$. Then, for every $j\geq K$,
\begin{equation}
  \M_j\leq4M_0.
  \label{eq:adaptive_near_regularization}
\end{equation}
Moreover,
\begin{equation}
  A_{K+n}^{1/3}
  \geq
  A_K^{1/3}
  +\frac{n}{3}\ls\frac{\Gamma_0}{288M_0}\rs^{1/3},
  \qquad n\geq1,
  \label{eq:adaptive_near_growth}
\end{equation}
and
\begin{gather}
  \Gamma_{K+n}\leq\frac{7776M_0}{n^3},
  \label{eq:adaptive_near_gamma_bound}\\
  f(w_{K+n})-f^\ast
  \leq\frac{5184M_0R^3}{n^3},
  \qquad n\geq1.
  \label{eq:adaptive_near_function_bound}
\end{gather}
\end{lemma}

\begin{proof}
The sequence $\{\Gamma_k\}_{k\geq0}$ is decreasing. Thus,
$\widehat B\Gamma_j\leq4M_0$ for every $j\geq K$, and
equation~\eqref{eq:adaptive_regularization_upper} proves
\eqref{eq:adaptive_near_regularization}. By~\eqref{eq:equivalent_weights},
\begin{equation*}
  a_{j+1}^3
  =\frac{\Gamma_0A_{j+1}^2}{72\M_j}
  \geq\frac{\Gamma_0A_{j+1}^2}{288M_0}.
\end{equation*}
Therefore,
\begin{gather*}
  A_{j+1}^{1/3}-A_j^{1/3}
  =\frac{a_{j+1}}
  {A_{j+1}^{2/3}+A_{j+1}^{1/3}A_j^{1/3}+A_j^{2/3}}\\
  \geq
  \frac{a_{j+1}}{3A_{j+1}^{2/3}}
  \geq
  \frac13\ls\frac{\Gamma_0}{288M_0}\rs^{1/3}.
\end{gather*}
Summing this inequality proves~\eqref{eq:adaptive_near_growth}. Hence,
\begin{equation*}
  A_{K+n}\geq\frac{n^3\Gamma_0}{7776M_0},
  \qquad
  \Gamma_{K+n}=\frac{\Gamma_0}{A_{K+n}}
  \leq\frac{7776M_0}{n^3}.
\end{equation*}
This proves~\eqref{eq:adaptive_near_gamma_bound}. Combining it with
\eqref{eq:adaptive_error_bound} proves
\eqref{eq:adaptive_near_function_bound}.
\end{proof}

\begin{proof}[Proof of Theorem~\ref{thm:fully_adaptive_convergence}]
If $\nabla f(x_0)=0$, then $x_0$ is optimal by convexity. Otherwise,
Lemma~\ref{lem:cubic_initialization} shows that the initialization
backtracking terminates, gives~\eqref{eq:inital_cnd2}, and proves the stated
bound on the number of rejected initialization trials.

  Lemma~\ref{lem:adaptive_error_bound} gives~\eqref{eq:adaptive_error_bound}.
It remains to consider $0<\e<f(w_0)-f^\ast$; otherwise the convergence
claim and the rejected-trial bound are immediate. In this case,
condition~\eqref{eq:inital_cnd2} implies
$\Gamma_0>3\e/R^3$, so the first index below is positive.
Let $N$ be the first index satisfying
\begin{equation*}
  \Gamma_N\leq\frac{3\e}{2R^3}.
\end{equation*}
Then~\eqref{eq:adaptive_error_bound} gives $f(w_N)-f^\ast\leq\e$.

As long as
  $\Gamma_k>\max\lb\Gamma_{\tr}^{\rm ad},3\e/(2R^3)\rb$,
Lemma~\ref{lem:adaptive_far_regime} gives geometric contraction. Therefore,
after at most
\begin{equation*}
  N_{\rm far}
  \eqdef
  \left\lceil
    \ls2+(144\widehat B)^{1/3}\rs
    \log\max\lb
      1,
      \frac{\Gamma_0}
      {\max\lb\Gamma_{\tr}^{\rm ad},3\e/(2R^3)\rb}
    \rb
  \right\rceil
\end{equation*}
iterations, either $\Gamma_k\leq3\e/(2R^3)$, or the near regime begins.
In the second case, Lemma~\ref{lem:adaptive_near_regime} shows that
\begin{equation*}
  N_{\rm near}
  \eqdef
  \left\lceil
    \ls\frac{5184M_0R^3}{\e}\rs^{1/3}
  \right\rceil
\end{equation*}
additional iterations suffice. Consequently, $N\leq N_{\rm far}+N_{\rm near}$.
By~\eqref{eq:adaptive_gamma_transition},
\begin{equation*}
  \frac{\Gamma_0}
  {\max\lb\Gamma_{\tr}^{\rm ad},3\e/(2R^3)\rb}
  =
  \min\lb
    \frac{\widehat B\Gamma_0}{4M_0},
    \frac{2\Gamma_0R^3}{3\e}
  \rb.
\end{equation*}
Thus,
\begin{equation}
\begin{aligned}
  N\leq
  \left\lceil
    \ls2+(144\widehat B)^{1/3}\rs
    \log\max\lb
      1,
      \min\lb
        \frac{\widehat B\Gamma_0}{4M_0},
        \frac{2\Gamma_0R^3}{3\e}
      \rb
    \rb
  \right\rceil
  +\left\lceil
    \ls\frac{5184M_0R^3}{\e}\rs^{1/3}
  \right\rceil.
\end{aligned}
  \label{eq:adaptive_explicit_complexity}
\end{equation}
Under~\eqref{as:gh}, definition~\eqref{eq:adaptive_B} gives
\begin{equation*}
  \widehat B^{1/3}
  =\cO\ls1+\sqrt{M_1}R\rs.
\end{equation*}
Under~\eqref{as:lh}, it gives
\begin{equation*}
  \widehat B^{1/3}
  =\cO\ls1+\sqrt{M_1}R+\frac{R}{\rho}\rs.
\end{equation*}
Substituting these estimates into~\eqref{eq:adaptive_explicit_complexity}
proves~\eqref{eq:adaptive_global_rate} and
\eqref{eq:adaptive_local_rate}.

The minimality of $N$ gives
\begin{equation*}
  \Gamma_{N-1}>\frac{3\e}{2R^3},
  \qquad
  \frac{2M_0}{\Gamma_{N-1}}
  <\frac{4M_0R^3}{3\e}.
\end{equation*}
Therefore,~\eqref{eq:adaptive_rejected_prefix} proves that
the main
backtracking loops reject at most
$
  \left\lceil
    \log_2\max\lb
      1,
      B,
      \tfrac{4M_0R^3}{3\e}
    \rb
  \right\rceil
$
trials.
\end{proof}

\section{Suboptimal A-NPE Method}
\label{app:ahpe_acceleration}

In this section, we propose suboptimal adaptive A-NPE method, which matches the convergence rate of Nesterov-type acceleration in exact Hessian case. For~\eqref{as:h} and~\eqref{as:inexact_hessian} it achieves convergence rate $\cO\ls \tfrac{LR^3}{N^3} + \tfrac{\delta R^2}{N^2} \rs$. 

We repeat the two backtracking tests for convenience:
\begin{equation}
  \PD(x,\tau): H(x)+\tau I\succ0, \quad
  \MS(x,\tau,s): \norm{\grad f(x+s)+\tau s}\leq\frac{\tau}{2}\norm{s}.
  \label{eq:ahpe_appendix_tests}
\end{equation}
For $w,v\in\R^d$ and $\tau>0$, we first check $\PD(w,\tau)$.
If it holds, we compute the trial step
\begin{equation}
  (H(w)+\tau I)s=-\grad f(v),
  \qquad
  y=v+s,
  \label{eq:ahpe_appendix_trial_step}
\end{equation}
and accept it when $\MS(v,\tau,s)$ holds.

The initialization backtracking keeps $x_0$, $\grad f(x_0)$,
and $H(x_0)$ fixed.
In the main loop, $H(w_k)$ remains fixed during backtracking, whereas
$v_k$ and $\grad f(v_k)$ are recomputed after each rejection.

\begin{algorithm}[!htb]
  \caption{Adaptive A-NPE with an inexact Hessian}
  \label{alg:adaptive_ahpe}
  \begin{algorithmic}[1]
    \STATE \textbf{Input:} $x_0\in\R^d$ and $\tau_0>0$.
    \STATE $(w_1, \tau_1)\eqdef\MSBacktrack(x_0,\tau_0)$.
    \STATE Set
      \begin{equation}
        A_1\eqdef\tfrac{1}{\tau_1},
        \qquad
        z_1\eqdef x_0-A_1\grad f(w_1).
        \label{eq:ahpe_initialization}
      \end{equation}
    \FOR{$k\geq1$}
      \STATE Query $H(w_k)$ once.
      \STATE Set
        \begin{equation}
          \alpha_k
          \eqdef
          \tfrac{2}{1+\sqrt{1+4A_k\tau_k}},
          \qquad
          A_{k+1}\eqdef\tfrac{A_k}{1-\alpha_k},
          \qquad
          a_{k+1}\eqdef\alpha_kA_{k+1}.
          \label{eq:ahpe_weights}
        \end{equation}
      \STATE Set $v_k\eqdef(1-\alpha_k)w_k+\alpha_kz_k$.
      \IF{$\PD(w_k,\tau_k)$}
        \STATE Compute $s_k$ and $w_{k+1}$ according to~\eqref{eq:ahpe_appendix_trial_step}
          with $w=w_k$, $v=v_k$, and $\tau=\tau_k$.
      \ENDIF
      \WHILE{\NOT $\PD(w_k,\tau_k)$
        \OR \NOT $\MS(v_k,\tau_k,s_k)$}
        \STATE Set $\tau_k:=2\tau_k$.
        \STATE Recompute $\alpha_k$, $A_{k+1}$, and $a_{k+1}$ according
          to~\eqref{eq:ahpe_weights}.
        \STATE Recompute $v_k\eqdef(1-\alpha_k)w_k+\alpha_kz_k$.
        \IF{$\PD(w_k,\tau_k)$}
          \STATE Recompute $s_k$ and $w_{k+1}$ according to~\eqref{eq:ahpe_appendix_trial_step}.
        \ENDIF
      \ENDWHILE
      \STATE Set
        \begin{equation*}
          z_{k+1}\eqdef z_k-a_{k+1}\grad f(w_{k+1}),
          \qquad
          \tau_{k+1}\eqdef(1-\alpha_k)\tau_k.
        \end{equation*}
    \ENDFOR
  \end{algorithmic}
\end{algorithm}

For the convergence bounds of Algorithm~\ref{alg:adaptive_ahpe}, we use
the convention $\rho=+\infty$ under~\eqref{as:gh} and define
\begin{equation}
  B\eqdef\max\lb1,8\sqrt{M_1}R,\tfrac{8R}{\rho}\rb.
  \label{eq:ahpe_B}
\end{equation}
If $R=0$, then $x_0$ is already a minimizer. For $R>0$, set
\begin{equation*}
  \e_{\tr}\eqdef
  \max\lb\tfrac{20M_0R^3}{B^3},\tfrac{4\delta R^2}{B^2}\rb.
\end{equation*}

\begin{theorem}
\label{thm:ahpe_convergence}
Suppose that $f$ is convex,~\eqref{as:inexact_hessian} holds, and
either~\eqref{as:gh} or~\eqref{as:lh} holds. Assume that $R>0$.
For every $\e>0$, the method generates an iterate $w_N$ satisfying
$f(w_N)-f^\ast\leq\e$ within
\begin{equation*}
  N=\cO\ls
    B\ls1+\log_+\tfrac{R^2}
                            {2A_1\max\lb\e,\e_{\tr}\rb}\rs
    +\ls\tfrac{M_0R^3}{\e}\rs^{1/3}
    +\sqrt{\tfrac{\delta R^2}{\e}}
  \rs
\end{equation*}
iterations. The number of rejections is at most
\begin{equation*}
  \left\lceil\log_2\max\lb
  1,\tfrac{4\delta}{\tau_0},
  \tfrac{2\sqrt{L(x_0)\norm{\grad f(x_0)}}}{\tau_0},
  \tfrac{2\norm{\grad f(x_0)}}{\rho\tau_0}
  \rb\right\rceil\allowbreak
  +\left\lceil\log_2\max\lb
  B^2,\tfrac{4\delta R^2}{\e},
  \ls\tfrac{20M_0R^3}{\e}\rs^{2/3}
  \rb\right\rceil.
\end{equation*}
Moreover, there exists an integer $K\geq1$ with
$K=\cO\ls B\max\lb1,\log\tfrac{R^2}{2A_1\e_{\tr}}\rb\rs$ such that
\begin{gather*}
  f(w_k)-f^\ast
  =\cO\ls\tfrac{R^2}{A_1}\exp\ls-\tfrac{c(k-1)}B\rs\rs,
  \qquad 1\leq k\leq K,\\
  f(w_{K+n})-f^\ast
  =\cO\ls\tfrac{M_0R^3}{(B+n)^3}+\tfrac{\delta R^2}{(B+n)^2}\rs,
  \qquad n=0,1,\ldots.
\end{gather*}
Here $c>0$ and the constants in $\cO(\cdot)$ are universal.
\end{theorem}

\subsection{Auxiliary Lemmas}

\begin{lemma}
\label{lem:ahpe_weight_relations}
For every main-loop trial of Algorithm~\ref{alg:adaptive_ahpe},
the parameters defined in~\eqref{eq:ahpe_weights} satisfy
\begin{equation*}
  A_{k+1}=A_k+a_{k+1},
  \quad \frac{a_{k+1}^2}{A_{k+1}}=\frac1{\tau_k},
  \quad
  v_k=\frac{A_kw_k+a_{k+1}z_k}{A_{k+1}},
  \quad A_{k+1}\tau_k=\frac1{\alpha_k^2}.
\end{equation*}
\end{lemma}

\begin{proof}
By~\eqref{eq:ahpe_weights}, we have
\begin{equation*}
  \alpha_k
  =\frac{2}{1+\sqrt{1+4A_k\tau_k}}
  =\frac{\sqrt{1+4A_k\tau_k}-1}{2A_k\tau_k}.
\end{equation*}
Thus $\alpha_k$ is the positive root of the quadratic equation
\begin{equation}
  A_k\tau_k\alpha_k^2+\alpha_k=1,
  \qquad \frac{\alpha_k^2}{1-\alpha_k}=\frac1{A_k\tau_k}.
  \label{eq:ahpe_alpha_identity}
\end{equation}
Using~\eqref{eq:ahpe_weights} and~\eqref{eq:ahpe_alpha_identity} gives
\begin{equation*}
  A_{k+1}-a_{k+1}
  =(1-\alpha_k)A_{k+1}=A_k,\qquad
  \frac{a_{k+1}^2}{A_{k+1}}
  =\alpha_k^2A_{k+1}
  =\frac{\alpha_k^2A_k}{1-\alpha_k}
  =\frac1{\tau_k}.
\end{equation*}
The definition of $v_k$ and the identities above imply
\begin{gather*}
  v_k=(1-\alpha_k)w_k+\alpha_kz_k
  =\frac{A_k}{A_{k+1}}w_k+\frac{a_{k+1}}{A_{k+1}}z_k
  =\frac{A_kw_k+a_{k+1}z_k}{A_{k+1}},\\
  A_{k+1}\tau_k
  =\frac{A_{k+1}^2}{a_{k+1}^2}
  =\frac1{\alpha_k^2}.
\end{gather*}
\end{proof}

\begin{lemma}
\label{lem:ahpe_frozen_remainder}
Suppose that either~\eqref{as:gh} or~\eqref{as:lh} holds.
Fix $w,v,s\in\R^d$, set $y\eqdef v+s$, and let $H(w)$
satisfy~\eqref{as:inexact_hessian} at $w$. Under~\eqref{as:lh}, assume additionally that
\begin{equation*}
  \norm{v-w}+\norm{s}\leq\rho.
\end{equation*}
Then
\begin{equation}
  \norm{\grad f(y)-\grad f(v)-H(w)s}
  \leq\ls\delta+L(y)\ls\norm{v-w}+\frac32\norm{s}\rs\rs\norm{s}.
  \label{eq:ahpe_frozen_remainder}
\end{equation}
\end{lemma}

\begin{proof}
By Lemmas~\ref{lem:global_taylor_remainders}~\ref{lem:local_taylor_remainders} and assumption of the lemma
\begin{equation*}
  \norm{\grad f(y)-\grad f(v)-\nabla^2f(y)s}
  =\norm{\grad f(y-s)-\grad f(y)+\nabla^2f(y)s}
  \leq\frac{L(y)}2\norm{s}^2.
\end{equation*}
Also,
\begin{equation*}
  \norm{y-w}=\norm{v+s-w}\leq\norm{v-w}+\norm{s},
\end{equation*}
which is at most $\rho$ under~\eqref{as:lh}.
Applying~\eqref{as:gh} or~\eqref{as:lh} with base point $y$,
together with~\eqref{as:inexact_hessian} at $w$, gives
\begin{gather*}
  \norm{\nabla^2f(y)-H(w)}
  \leq\norm{\nabla^2f(y)-\nabla^2f(w)}
       +\norm{\nabla^2f(w)-H(w)}\\
  \leq L(y)\norm{y-w}+\delta
  \leq\delta+L(y)\ls\norm{v-w}+\norm{s}\rs.
\end{gather*}
Adding and subtracting $\nabla^2f(y)s$, we obtain
\begin{gather*}
  \norm{\grad f(y)-\grad f(v)-H(w)s}
  \leq\norm{\grad f(y)-\grad f(v)-\nabla^2f(y)s}
       +\norm{\nabla^2f(y)-H(w)}\norm{s}\\
  \leq\frac{L(y)}2\norm{s}^2
       +\ls\delta+L(y)\ls\norm{v-w}+\norm{s}\rs\rs\norm{s}
  =\ls\delta+L(y)\ls\norm{v-w}+\frac32\norm{s}\rs\rs\norm{s}.
\end{gather*}
This proves~\eqref{eq:ahpe_frozen_remainder}.
\end{proof}

For a trial step~\eqref{eq:ahpe_appendix_trial_step}, the linear system gives
\begin{equation*}
  \grad f(y)+\tau s
  =\grad f(y)-\grad f(v)-H(w)s.
\end{equation*}

\subsection{Step Estimates}
\begin{lemma}
\label{lem:ahpe_energy}
Let $f$ be convex and differentiable.
Fix $w,z\in\R^d$ and $A\geq0$. Choose $a>0$ and $\tau>0$ such that
\begin{equation*}
  A_+\eqdef A+a,
  \qquad \frac{a^2}{A_+}=\frac1\tau,
  \qquad v\eqdef\frac{A}{A_+}w+\frac{a}{A_+}z.
\end{equation*}
Let $\MS(v,\tau,s)$ in~\eqref{eq:ahpe_appendix_tests} hold, i.e.
\begin{equation}
  \label{eq:app_test}
  \norm{\grad f(v+s)+\tau s}\leq\frac{\tau}{2}\norm{s},
\end{equation}
and define
\begin{equation*}
  y\eqdef v+s,
  \qquad z_+\eqdef z-a\grad f(y).
\end{equation*}
Then
\begin{equation}
\begin{aligned}
  A_+(f(y)-f^\ast)+\frac12\norm{z_+-x^\ast}^2
  \leq A(f(w)-f^\ast)+\frac12\norm{z-x^\ast}^2 -\frac38 A_+\tau\norm{s}^2.
\end{aligned}
\label{eq:ahpe_energy_step}
\end{equation}
\end{lemma}

\begin{proof}
Set $g\eqdef\grad f(y)$. Squaring~\eqref{eq:app_test} gives
\begin{equation*}
  \norm{g}^2+2\tau\la g,s\ra+\tau^2\norm{s}^2
  =\norm{g+\tau s}^2
  \leq\frac{\tau^2}{4}\norm{s}^2.
\end{equation*}
Rearranging and dividing by $2\tau>0$, we obtain
\begin{equation}
  \la g,s\ra
  \leq-\frac{1}{2\tau}\norm{g}^2-\frac{3\tau}{8}\norm{s}^2.
  \label{eq:ahpe_residual_progress}
\end{equation}
Convexity at $y$ gives
\begin{equation*}
  f(w)\geq f(y)+\la g,w-y\ra,
  \qquad f^\ast\geq f(y)+\la g,x^\ast-y\ra.
\end{equation*}
Multiply these inequalities by $A$ and $a$, respectively, and add them:
\begin{equation*}
  Af(w)+af^\ast
  \geq A_+f(y)+\la g,Aw+ax^\ast-A_+y\ra.
\end{equation*}
The identities $A_+v=Aw+az$ and $y=v+s$ imply
\begin{equation*}
  Aw+ax^\ast-A_+y
  =Aw+ax^\ast-(Aw+az)-A_+s
  =a(x^\ast-z)-A_+s.
\end{equation*}
Substituting this identity and using $A_+=A+a$ yields
\begin{equation}
  A_+(f(y)-f^\ast)
  \leq A(f(w)-f^\ast)+a\la g,z-x^\ast\ra+A_+\la g,s\ra.
  \label{eq:ahpe_convexity_step}
\end{equation}
The update $z_+=z-ag$ gives
\begin{equation}
  \frac12\norm{z_+-x^\ast}^2
  =\frac12\norm{z-x^\ast-ag}^2
  =\frac12\norm{z-x^\ast}^2-a\la g,z-x^\ast\ra
    +\frac{a^2}{2}\norm{g}^2.
\label{eq:ahpe_distance_step}
\end{equation}
Adding~\eqref{eq:ahpe_convexity_step} and~\eqref{eq:ahpe_distance_step} gives 
\begin{equation*}
  A_+(f(y)-f^\ast)+\frac12\norm{z_+-x^\ast}^2
  \leq A(f(w)-f^\ast)+\frac12\norm{z-x^\ast}^2
  +A_+\la g,s\ra+\frac{a^2}{2}\norm{g}^2.
\end{equation*}
Finally,~\eqref{eq:ahpe_residual_progress} and $a^2=A_+/\tau$ give
\begin{equation*}
  A_+\la g,s\ra+\frac{a^2}{2}\norm{g}^2
  \leq\ls\frac{a^2}{2}-\frac{A_+}{2\tau}\rs\norm{g}^2
       -\frac38 A_+\tau\norm{s}^2
  =-\frac38 A_+\tau\norm{s}^2.
\end{equation*}
Substitution proves~\eqref{eq:ahpe_energy_step}.
\end{proof}

Lemma~\ref{lem:ahpe_weight_relations} shows that the weights used in
the main loop of Algorithm~\ref{alg:adaptive_ahpe} satisfy the relations
assumed in Lemma~\ref{lem:ahpe_energy}.

To bound the potential after initialization step in Algorithm~\ref{alg:adaptive_ahpe}, we apply
Lemma~\ref{lem:ahpe_energy} with $A=0$ and $w=z=x_0$.
For the accepted step, take $a=A_+=A_1=1/\tau_1$
and $\tau=\tau_1$. These choices satisfy the weight condition,
and the updates in~\eqref{eq:ahpe_initialization} give $y=w_1$ and
$z_+=z_1$. Since $A_1\tau_1=1$, the lemma yields
\begin{equation}
  A_1(f(w_1)-f^\ast)+\frac12\norm{z_1-x^\ast}^2
  +\frac38\norm{w_1-x_0}^2
  \leq\frac{R^2}{2}.
  \label{eq:ahpe_initialization_bound}
\end{equation}
For the accepted iterates, write
\begin{equation*}
  \ECal_k\eqdef A_k(f(w_k)-f^\ast)+\frac12\norm{z_k-x^\ast}^2.
\end{equation*}
For any accepted iterate $w_k$, summing~\eqref{eq:ahpe_energy_step}
over $i=1,\ldots,k-1$ and using the initialization bound gives
\begin{equation}
  \ECal_k+\frac38\sum_{i=1}^{k-1}A_{i+1}\tau_i\norm{s_i}^2
  \leq\ECal_1\leq\frac{R^2}{2},
  \label{eq:ahpe_energy_budget}
\end{equation}
where $s_i=w_{i+1}-v_i$ and $\tau_i$ is accepted on iteration $i$.
Since all terms on the left-hand side of~\eqref{eq:ahpe_energy_budget}
are nonnegative,
\begin{equation*}
  A_k(f(w_k)-f^\ast) \leq \ECal_k\leq\frac{R^2}{2},
  \qquad
  \sum_{i=1}^{k-1}A_{i+1}\tau_i\norm{s_i}^2\leq\frac43R^2.
\end{equation*}
Consider any trial satisfying the
conditions of Lemma~\ref{lem:ahpe_energy} with $A=A_k$, $w=w_k$, and
$z=z_k$, and write $\alpha=a/A_+$. Let $y\eqdef v+s$ be the trial point.
For an accepted main-loop step of Algorithm~\ref{alg:adaptive_ahpe},
$y=v_k+s_k=w_{k+1}$. The left-hand side
of~\eqref{eq:ahpe_energy_step} is nonnegative. Hence,
by~\eqref{eq:ahpe_energy_budget},
\begin{equation*}
  \frac38 A_+\tau\norm{s}^2\leq\ECal_k\leq\frac{R^2}{2}.
\end{equation*}
The MS test~\eqref{eq:app_test} also gives 
\begin{equation*}
  \norm{\grad f(y)}
  \leq\norm{\grad f(y)+\tau s}+\tau\norm{s}
  \leq\frac{\tau}{2}\norm{s}+\tau\norm{s}
  =\frac{3\tau}{2}\norm{s}.
\end{equation*}
Using $A_+\tau=A_+^2/a^2=1/\alpha^2$, we obtain
\begin{equation}
  \norm{s}\leq\frac{2}{\sqrt3}\alpha R,
  \qquad
  \norm{\grad f(y)}\leq\frac32\tau\norm{s}.
  \label{eq:ahpe_certified_bounds}
\end{equation}

\begin{lemma}
\label{lem:ahpe_localization}
For every $k\geq1$, the accepted iterates of
Algorithm~\ref{alg:adaptive_ahpe} satisfy
\begin{equation}
  \norm{z_k-x^\ast}\leq R,
  \label{eq:ahpe_z_bound}
\end{equation}
\begin{equation}
  \norm{w_k-x^\ast}\leq\ls1+\frac{2}{\sqrt3}\rs R,
  \qquad
  \norm{z_k-w_k}\leq\ls2+\frac{2}{\sqrt3}\rs R.
  \label{eq:ahpe_localization}
\end{equation}
\end{lemma}

\begin{proof}
By~\eqref{eq:ahpe_energy_budget},
\begin{equation*}
  \frac12\norm{z_k-x^\ast}^2\leq\ECal_k\leq\frac{R^2}{2},
\end{equation*}
which gives~\eqref{eq:ahpe_z_bound}.
We prove the bound on $\norm{w_k-x^\ast}$ by induction.
The initialization bound~\eqref{eq:ahpe_initialization_bound} gives
\begin{gather*}
  \norm{w_1-x_0}\leq\frac{2}{\sqrt3}R,\\
  \norm{w_1-x^\ast}
  \leq\norm{w_1-x_0}+\norm{x_0-x^\ast}
  \leq\ls1+\frac{2}{\sqrt3}\rs R.
\end{gather*}
Suppose the bound holds for $w_k$. For the next accepted step,
$w_{k+1}=(1-\alpha_k)w_k+\alpha_kz_k+s_k$.
Using~\eqref{eq:ahpe_certified_bounds} and~\eqref{eq:ahpe_z_bound}, we obtain
\begin{gather*}
  \norm{w_{k+1}-x^\ast}
  \leq(1-\alpha_k)\norm{w_k-x^\ast}
       +\alpha_k\norm{z_k-x^\ast}+\norm{s_k}\\
  \leq(1-\alpha_k)\ls1+\frac{2}{\sqrt3}\rs R
       +\alpha_k R+\frac{2}{\sqrt3}\alpha_k R
  =\ls1+\frac{2}{\sqrt3}\rs R.
\end{gather*}
This completes the induction. The bound on $\norm{z_k-w_k}$ follows from
\begin{equation*}
  \norm{z_k-w_k}
  \leq\norm{z_k-x^\ast}+\norm{w_k-x^\ast}
  \leq\ls2+\frac{2}{\sqrt3}\rs R.
\end{equation*}
\end{proof}

Fix already accepted $w_k,z_k$. These points remain fixed during
backtracking. For $v=(1-\alpha)w_k+\alpha z_k$ with any $\alpha\in(0,1)$,
inequality~\eqref{eq:ahpe_localization} gives
\begin{equation*}
  \norm{v-w_k}=\alpha\norm{z_k-w_k}
  \leq\ls2+\frac{2}{\sqrt3}\rs\alpha R.
\end{equation*}
This bound holds for every trial, whether or not the MS test passes.

\subsection{Explicit Choice of Regularization Parameter}

With $B$ defined in~\eqref{eq:ahpe_B}, set
\begin{equation}
  Q(A)\eqdef\max\lb B^2,8\delta A,(40M_0RA)^{2/3}\rb,
  \qquad A>0.
  \label{eq:ahpe_safe_threshold}
\end{equation}

\begin{lemma}
\label{lem:ahpe_safe_parameter}
Suppose that $f$ is convex and either~\eqref{as:gh} or~\eqref{as:lh} holds.
Fix an accepted state $A_k,w_k,z_k$ of Algorithm~\ref{alg:adaptive_ahpe}
and a matrix $H(w_k)$ satisfying~\eqref{as:inexact_hessian} at $w_k$.
Every main-loop trial with
\begin{equation}
  \tau\geq\frac{Q(A_k)}{A_k}
  \label{eq:ahpe_safe_parameter}
\end{equation}
satisfies $\PD(w_k,\tau)$ and $\MS(v,\tau,s)$
in~\eqref{eq:ahpe_appendix_tests}, where the trial weights are defined
by~\eqref{eq:ahpe_weights} with $\tau_k=\tau$,
$v=(1-\alpha_k)w_k+\alpha_kz_k$, and $s$ is given
by~\eqref{eq:ahpe_appendix_trial_step} with $w=w_k$.
\end{lemma}

\begin{proof}
Keep $A_k,w_k,z_k,H(w_k)$ fixed and vary $q\geq Q(A_k)$ continuously.
For each $q$, use the weights in~\eqref{eq:ahpe_weights} with
$\tau_k=\tau(q)$, where
\begin{gather*}
  \tau(q)\eqdef\frac{q}{A_k},
  \quad
  \alpha(q)\eqdef\frac{2}{1+\sqrt{1+4q}}, \quad
  v(q)\eqdef(1-\alpha(q))w_k+\alpha(q)z_k.
\end{gather*}
Convexity and~\eqref{as:inexact_hessian} give $H(w_k)\succeq-\delta I$.
Since $q\geq Q(A_k)\geq\max\lb1,8\delta A_k\rb$, we have
\begin{equation*}
  H(w_k)+\tau(q)I
  \succeq(\tau(q)-\delta)I
  \succeq\frac78\tau(q)I\succ0.
\end{equation*}
Thus, $\PD(w_k,\tau(q))$ holds throughout this range, and
we can define
\begin{gather*}
  s(q)\eqdef-\ls H(w_k)+\tau(q)I\rs^{-1}\grad f(v(q)),
  \qquad y(q)\eqdef v(q)+s(q),\\
  d(q)\eqdef\norm{v(q)-w_k}=\alpha(q)\norm{z_k-w_k},
  \qquad r(q)\eqdef\norm{s(q)}.
\end{gather*}
All these quantities are continuous because $H(w_k)$ is fixed and the
inverse exists. As $q\to+\infty$,
\begin{gather*}
  \tau(q)\to+\infty,~~\alpha(q)\to0,~~v(q)\to w_k,~~ d(q)\to0,~~
  r(q)\leq\frac{\norm{\grad f(v(q))}}{\tau(q)-\delta}\to0,~~ y(q)\to w_k.
\end{gather*}
By Lemma~\ref{lem:ahpe_frozen_remainder} and~\eqref{eq:ahpe_appendix_trial_step},
the following conditions are sufficient for the test~\eqref{eq:ahpe_appendix_tests}:
\begin{equation}
\begin{aligned}
  d(q)+r(q) \leq\rho~~\text{only under~\eqref{as:lh}},\quad
  \delta+L(y(q))\ls d(q)+\frac32r(q)\rs \leq\frac{\tau(q)}2.
\end{aligned}
  \label{eq:ahpe_sufficient_conditions}
\end{equation}
Indeed, if~\eqref{eq:ahpe_sufficient_conditions} holds, then
Lemma~\ref{lem:ahpe_frozen_remainder} and~\eqref{eq:ahpe_appendix_trial_step} give
\begin{gather*}
  \norm{\grad f(y(q))+\tau(q)s(q)}
  =\norm{\grad f(y(q))-\grad f(v(q))-H(w_k)s(q)}\\
  \leq\ls\delta+L(y(q))\ls d(q)+\frac32r(q)\rs\rs r(q)
  \leq\frac{\tau(q)}2\norm{s(q)}.
\end{gather*}
Since $L(y(q))\to L(w_k)$,
\begin{equation*}
  \frac{\delta+L(y(q))\ls d(q)+\frac32r(q)\rs}{\tau(q)}\to0.
\end{equation*}
Hence, all sufficiently large $q$ satisfy~\eqref{eq:ahpe_sufficient_conditions}
strictly.

Let $\mathcal U_k$ be the set of $q\geq Q(A_k)$ for which at least one
condition in~\eqref{eq:ahpe_sufficient_conditions} fails.
If $\mathcal U_k=\emptyset$, the statement follows.
Otherwise, set $q^\dagger\eqdef\sup\mathcal U_k$.
Since the conditions in~\eqref{eq:ahpe_sufficient_conditions} hold for all
sufficiently large $q$, the set $\mathcal U_k$ is bounded above.
Hence, $q^\dagger<+\infty$.
By continuity, $q^\dagger$ lies on the boundary of the admissible set.
Hence, the conditions in~\eqref{eq:ahpe_sufficient_conditions} hold at
$q^\dagger$, and at least one is an equality.
A dagger below denotes evaluation at $q^\dagger$.

The boundary trial therefore satisfies the test~\eqref{eq:ahpe_sufficient_conditions}.
Moreover,~\eqref{eq:ahpe_energy_budget} gives $\ECal_k\leq\frac{R^2}{2}$.
Consequently, the bounds~\eqref{eq:ahpe_certified_bounds} apply to this trial.
Together with~\eqref{eq:ahpe_localization}, they give
\begin{equation}
  r^\dagger\leq\frac{2}{\sqrt3}\alpha^\dagger R, 
  \quad
  \norm{\grad f(y^\dagger)}
  \leq\frac32\tau^\dagger r^\dagger
  \leq\sqrt3\tau^\dagger\alpha^\dagger R, \quad
  d^\dagger=\alpha^\dagger\norm{z_k-w_k}
  \leq\ls2+\frac{2}{\sqrt3}\rs\alpha^\dagger R.
  \label{eq:ahpe_boundary_bounds}
\end{equation}
In particular,
\begin{equation}
  d^\dagger+r^\dagger
  \leq d^\dagger+\frac32r^\dagger
  \leq\ls2+\frac{5}{\sqrt3}\rs\alpha^\dagger R
  \leq5\alpha^\dagger R.
  \label{eq:ahpe_boundary_distances}
\end{equation}
The definition of $\alpha(q)$ gives
\begin{equation*}
  \alpha^\dagger=\frac{2}{1+\sqrt{1+4q^\dagger}}
  \leq\frac1{\sqrt{q^\dagger}}.
\end{equation*}
By~\eqref{eq:ahpe_safe_threshold}, we have $q^\dagger\geq Q(A_k)\geq B^2$.
Under~\eqref{as:lh}, definition~\eqref{eq:ahpe_B} also gives
$B\geq\frac{8R}{\rho}$.
Together with~\eqref{eq:ahpe_boundary_distances}, this gives
\begin{equation*}
  d^\dagger+r^\dagger
  \leq\frac{5R}{\sqrt{q^\dagger}}
  \leq\frac{5R}{B}
  \leq\frac58\rho<\rho.
\end{equation*}
Using $L(y^\dagger)=M_0+M_1\norm{\grad f(y^\dagger)}$,
equations~\eqref{eq:ahpe_boundary_bounds} and~\eqref{eq:ahpe_boundary_distances}
give
\begin{gather*}
  \frac{\delta+L(y^\dagger)\ls d^\dagger+\frac32r^\dagger\rs}
       {\tau^\dagger}
  \leq\frac{\delta}{\tau^\dagger}
       +\frac{5M_0R\alpha^\dagger}{\tau^\dagger}
       +5\sqrt3M_1R^2(\alpha^\dagger)^2\\
  \leq\frac{\delta A_k}{q^\dagger}
       +\frac{5M_0RA_k}{(q^\dagger)^{3/2}}
       +\frac{5\sqrt3M_1R^2}{q^\dagger}
  \leq\frac18+\frac18+\frac{5\sqrt3}{64}<\frac12,
\end{gather*}
where we used $q^\dagger\geq Q(A_k)$,~\eqref{eq:ahpe_safe_threshold},
and~\eqref{eq:ahpe_B}.
Thus, all conditions in~\eqref{eq:ahpe_sufficient_conditions} are strict
at $q^\dagger$, contradicting its boundary property.
Therefore, $\mathcal U_k=\emptyset$, and every trial
satisfying~\eqref{eq:ahpe_safe_parameter} passes both tests.
\end{proof}

\subsection{Main-Loop Backtracking}

\begin{lemma}
\label{lem:ahpe_backtracking}
Suppose that $f$ is convex,~\eqref{as:inexact_hessian} holds, and
either~\eqref{as:gh} or~\eqref{as:lh} holds.
Every main backtracking loop of Algorithm~\ref{alg:adaptive_ahpe} terminates.
For every $k\geq1$, let $\tau_k$ denote the parameter accepted at iteration $k$.
Then
\begin{equation}
  A_k\tau_k\leq2Q(A_k).
  \label{eq:ahpe_parameter_upper}
\end{equation}
For every $K\geq1$, let $J_K$ be the total number of rejected main trials
during iterations $1,\ldots,K$, excluding initialization. Then
\begin{equation}
  J_K\leq\left\lceil\log_2 Q(A_K)\right\rceil.
  \label{eq:ahpe_rejected_prefix}
\end{equation}
\end{lemma}

\begin{proof}
We prove by induction on $k$ that iteration $k$ is reached after
finitely many trials. The case $k=1$ follows from
Lemma~\ref{lem:pd_ms_backtracking} with $(w_1, \tau_1) = \MSBacktrack(x_0, \tau_0)$ and $\overline{\tau}_0 \eqdef \overline{\tau}(x_0)$.
Suppose this holds for some $k\geq1$.
Then~\eqref{eq:ahpe_energy_budget} and~\eqref{eq:ahpe_localization} hold.
Since $A_k$ is fixed, repeated doubling reaches
$\tau\geq\frac{Q(A_k)}{A_k}$, and Lemma~\ref{lem:ahpe_safe_parameter}
guarantees acceptance. Thus, iteration $k+1$ is reached after finitely
many trials, completing the induction.

We next track the product of the current weight and trial parameter.
By~\eqref{eq:ahpe_initialization}, before the first main trial this product is
\begin{equation*}
  A_1\tau_1=\frac1{\tau_1}{\tau_1}=1.
\end{equation*}
After accepting iteration $k$, the first trial at iteration $k+1$ uses
$\tau=(1-\alpha_k)\tau_k$. By~\eqref{eq:ahpe_weights}, its product is
\begin{equation*}
  A_{k+1}(1-\alpha_k)\tau_k
  =\frac{A_k}{1-\alpha_k}(1-\alpha_k)\tau_k
  =A_k\tau_k.
\end{equation*}
Every rejected trial doubles $\tau_k$, while $A_k$ remains unchanged.
Therefore,
\begin{equation}
  A_k\tau_k=2^{J_k},\qquad k\geq1.
  \label{eq:ahpe_rejection_product}
\end{equation}
Fix $K\geq1$. The weights increase by~\eqref{eq:ahpe_weights}, and
$Q$ is nondecreasing by~\eqref{eq:ahpe_safe_threshold}. Hence,
$Q(A_i)\leq Q(A_K)$ for every $i\leq K$.
If $J_K=0$, then~\eqref{eq:ahpe_rejection_product} and $Q(A_K)\geq B^2\geq1$
give both~\eqref{eq:ahpe_parameter_upper} and~\eqref{eq:ahpe_rejected_prefix}.

Suppose that $J_K\geq1$, and let $i\leq K$ be the iteration at which the
last rejection occurs. For this trial, the product of $A_i$ and the trial
parameter is $2^{J_K-1}$. By Lemma~\ref{lem:ahpe_safe_parameter}, the trial
would be accepted if this product were at least $Q(A_i)$. Therefore,
\begin{equation*}
  2^{J_K-1}<Q(A_i)\leq Q(A_K).
\end{equation*}
Multiplying by $2$ and using~\eqref{eq:ahpe_rejection_product}, we obtain
\begin{equation*}
  A_K\tau_K=2^{J_K}<2Q(A_K),
\end{equation*}
which proves~\eqref{eq:ahpe_parameter_upper}.
Taking logarithms in $2^{J_K-1}<Q(A_K)$ gives
$J_K<1+\log_2 Q(A_K)$. Since $J_K$ is an integer,
\begin{equation*}
  J_K\leq\left\lceil\log_2 Q(A_K)\right\rceil,
\end{equation*}
which proves~\eqref{eq:ahpe_rejected_prefix}.
\end{proof}

\subsection{Proof of Theorem~\ref{thm:ahpe_convergence}}

For the next lemma, recall the definitions
in~\eqref{eq:ahpe_B} and~\eqref{eq:ahpe_safe_threshold}:
\begin{gather*}
  B=
  \begin{cases}
    \max\lb1,8\sqrt{M_1}R\rb,
    &\text{under~\eqref{as:gh}},\\
    \max\lb1,8\sqrt{M_1}R,\frac{8R}{\rho}\rb,
    &\text{under~\eqref{as:lh}},
  \end{cases}\\
  Q(A)=\max\lb B^2,8\delta A,(40M_0RA)^{2/3}\rb,
  \qquad A>0.
\end{gather*}
We also recall the error threshold
\begin{equation}
  \e_{\tr}=
  \max\lb\frac{20M_0R^3}{B^3},\frac{4\delta R^2}{B^2}\rb.
  \label{eq:app_ahpe_transition_error}
\end{equation}
For $R>0$, set
\begin{equation}
  A_{\tr}\eqdef
  \min\lb\frac{B^3}{40M_0R},\frac{B^2}{8\delta}\rb
  =\frac{R^2}{2\e_{\tr}}.
  \label{eq:ahpe_transition_weight}
\end{equation}
Here $\frac{B^2}{8\delta}=+\infty$ when $\delta=0$.
In particular,
\begin{equation*}
  \sqrt{Q(A)}
  =\max\lb B,(40M_0RA)^{1/3},\sqrt{8\delta A}\rb,
  \qquad A>0.
\end{equation*}
The threshold $A_{\tr}$ in~\eqref{eq:ahpe_transition_weight}
is the largest weight for which $Q(A)=B^2$.

\begin{lemma}
\label{lem:ahpe_weight_growth}
Under the assumptions of Lemma~\ref{lem:ahpe_backtracking}, suppose that $R>0$.
For the accepted steps of Algorithm~\ref{alg:adaptive_ahpe}, the following
bounds hold for every $k\geq1$.
If $A_k\leq A_{\tr}$, then
\begin{equation}
  A_{k+1}\geq\exp\ls\frac1{3B}\rs A_k.
  \label{eq:ahpe_geometric_growth}
\end{equation}
If $A_k\geq A_{\tr}$, then
\begin{equation}
  \sqrt{Q(A_{k+1})}-\sqrt{Q(A_k)}\geq\frac19.
  \label{eq:ahpe_scale_growth}
\end{equation}
In particular, for any $K\geq1$ with $A_K\geq A_{\tr}$ and
every integer $n\geq0$,
\begin{equation}
  \sqrt{Q(A_{K+n})}\geq\sqrt{Q(A_K)}+\frac n9\geq B+\frac n9.
  \label{eq:ahpe_cumulative_growth}
\end{equation}
\end{lemma}

\begin{proof}
By~\eqref{eq:ahpe_parameter_upper}, the accepted parameter satisfies
$A_k\tau_k\leq2Q(A_k)$. Since $\sqrt{Q(A_k)}\geq B\geq1$,
the definition of $\alpha_k$ in~\eqref{eq:ahpe_weights} gives
\begin{equation*}
  \frac1{\alpha_k}
  =\frac{1+\sqrt{1+4A_k\tau_k}}2
  \leq\frac{1+\sqrt{1+8Q(A_k)}}2
  \leq3\sqrt{Q(A_k)}.
\end{equation*}
Hence,
\begin{equation}
  \alpha_k\geq\frac1{3\sqrt{Q(A_k)}},\qquad k\geq1.
  \label{eq:ahpe_alpha_lower}
\end{equation}
If $A_k\leq A_{\tr}$,
$\max\lb(40M_0RA_k)^{1/3},\sqrt{8\delta A_k}\rb \leq B$. Thus, by definition of $Q(A_k)$, we have $B=\sqrt{Q(A_k)}$. By~\eqref{eq:ahpe_alpha_lower}, $\alpha_k\geq\frac1{3B}$.
Using~\eqref{eq:ahpe_weights} and
$1-t\leq\exp(-t)$ for $0\leq t<1$, we obtain
\begin{equation*}
  A_{k+1}=\frac{A_k}{1-\alpha_k}
  \geq\exp(\alpha_k)A_k
  \geq\exp\ls\frac1{3B}\rs A_k,
\end{equation*}
which proves~\eqref{eq:ahpe_geometric_growth}.

Since $A\geq A_{\tr}$, $B\leq\max\lb(40M_0RA)^{1/3},\sqrt{8\delta A}\rb=\sqrt{Q(A)}$.
Suppose that $A_k\geq A_{\tr}$.
Since $A_{k+1}>A_k\geq A_{\tr}$, the identity holds
for $A=A_k$ and $A=A_{k+1}$.
Since $0<1-\alpha_k<1$, substitution
of~\eqref{eq:ahpe_weights} gives
\begin{gather*}
  \sqrt{Q(A_{k+1})}
  =\max\lb
    \frac{(40M_0RA_k)^{1/3}}{(1-\alpha_k)^{1/3}},
    \frac{\sqrt{8\delta A_k}}{(1-\alpha_k)^{1/2}}
  \rb\\
  \geq\frac{\max\lb(40M_0RA_k)^{1/3},\sqrt{8\delta A_k}\rb}
               {(1-\alpha_k)^{1/3}}
  =\frac{\sqrt{Q(A_k)}}{(1-\alpha_k)^{1/3}}.
\end{gather*}
The inequality $(1-t)^{-1/3}\geq1+\frac t3$ for $0\leq t<1$
and~\eqref{eq:ahpe_alpha_lower} now imply
\begin{equation*}
  \sqrt{Q(A_{k+1})}-\sqrt{Q(A_k)}
  \geq\sqrt{Q(A_k)}\ls(1-\alpha_k)^{-1/3}-1\rs
  \geq\frac{\alpha_k\sqrt{Q(A_k)}}3
  \geq\frac19.
\end{equation*}
This proves~\eqref{eq:ahpe_scale_growth}.
Finally, if $A_K\geq A_{\tr}$, then $A_j\geq A_{\tr}$
for every $j\geq K$. Summing~\eqref{eq:ahpe_scale_growth} gives
\begin{equation*}
  \sqrt{Q(A_{K+n})}
  =\sqrt{Q(A_K)}
  +\sum_{j=K}^{K+n-1}\ls\sqrt{Q(A_{j+1})}-\sqrt{Q(A_j)}\rs
  \geq\sqrt{Q(A_K)}+\frac n9
  \geq B+\frac n9,
\end{equation*}
which proves~\eqref{eq:ahpe_cumulative_growth}.
\end{proof}

We restate Theorem~\ref{thm:ahpe_convergence} for convenience.

\begin{theorem}
\label{thm:ahpe_convergence_appendix}
Suppose that $f$ is convex,~\eqref{as:inexact_hessian} holds, and
either~\eqref{as:gh} or~\eqref{as:lh} holds. Assume that $R>0$.
All backtracking loops in Algorithm~\ref{alg:adaptive_ahpe} terminate.
For every $\e>0$, the method generates an iterate $w_N$ satisfying
$f(w_N)-f^\ast\leq\e$ within
\begin{equation}
  N=\cO\ls
    B\ls1+\log_+\frac{R^2}
                            {2A_1\max\lb\e,\e_{\tr}\rb}\rs
    +\ls\frac{M_0R^3}{\e}\rs^{1/3}
    +\sqrt{\frac{\delta R^2}{\e}}
  \rs
  \label{eq:app_ahpe_iteration_complexity}
\end{equation}
iterations. Under~\eqref{as:gh}, initialization rejects at most
\begin{equation}
  \left\lceil\log_2\max\lb
    1,\frac{4\delta}{\tau_0},
    \frac{2\sqrt{L(x_0)\norm{\grad f(x_0)}}}{\tau_0}
  \rb\right\rceil
  \label{eq:app_ahpe_initialization_cost}
\end{equation}
trials. Under~\eqref{as:lh}, the maximum additionally includes
$\frac{2\norm{\grad f(x_0)}}{\rho\tau_0}$.
The total number of rejected main-loop trials before producing $w_N$
is at most
\begin{equation}
  \left\lceil\log_2\max\lb
    B^2,\frac{4\delta R^2}{\e},
    \ls\frac{20M_0R^3}{\e}\rs^{2/3}
  \rb\right\rceil.
  \label{eq:app_ahpe_accuracy_rejections}
\end{equation}
Moreover, there exists an integer $K\geq1$ with
\begin{equation}
  K\leq1+\left\lceil
    3B\log_+\frac{R^2}{2A_1\e_{\tr}}
  \right\rceil.
  \label{eq:app_ahpe_transition_index}
\end{equation}
For this $K$,
\begin{gather}
  f(w_k)-f^\ast
  \leq\frac{R^2}{2A_1}\exp\ls-\frac{k-1}{3B}\rs,
  \qquad 1\leq k\leq K,
  \label{eq:app_ahpe_far_rate}\\
  f(w_{K+n})-f^\ast
  \leq\max\lb
    \frac{20M_0R^3}{(B+n/9)^3},
    \frac{4\delta R^2}{(B+n/9)^2}
  \rb,
  \qquad n=0,1,\ldots.
  \label{eq:app_ahpe_near_rate}
\end{gather}
\end{theorem}

\begin{proof}
Lemmas~\ref{lem:pd_ms_backtracking}
and~\ref{lem:ahpe_backtracking} show that all backtracking loops terminate.
By~\eqref{eq:ahpe_energy_budget},
\begin{equation}
  f(w_k)-f^\ast\leq\frac{\ECal_k}{A_k}\leq\frac{R^2}{2A_k},
  \qquad k\geq1.
  \label{eq:ahpe_error_certificate}
\end{equation}
Set $K\eqdef\min\lb k\geq1:A_k\geq A_{\tr}\rb$, with $K=+\infty$
if this set is empty.
If $A_1\geq A_{\tr}$, then $K=1$,
and~\eqref{eq:app_ahpe_transition_index} holds.
Otherwise, set
\begin{equation*}
  m\eqdef\left\lceil3B\log\frac{A_{\tr}}{A_1}\right\rceil.
\end{equation*}
If $A_j<A_{\tr}$ for $j=1,\ldots,m$,
then~\eqref{eq:ahpe_geometric_growth} yields
\begin{equation*}
  A_{m+1}\geq A_1\exp\ls\frac{m}{3B}\rs
  \geq A_1\exp\ls\log\frac{A_{\tr}}{A_1}\rs
  =A_{\tr}.
\end{equation*}
The threshold is therefore reached by index $m+1$.
Since $A_{\tr}=\frac{R^2}{2\e_{\tr}}$,
this proves~\eqref{eq:app_ahpe_transition_index}.

Fix $1\leq k\leq K$.
By the definition of $K$, we have $A_j<A_{\tr}$ for $1\leq j<k$.
Applying~\eqref{eq:ahpe_geometric_growth} successively and
using~\eqref{eq:ahpe_error_certificate} gives
\begin{gather*}
  A_k\geq A_1\exp\ls\frac{k-1}{3B}\rs,\\
  f(w_k)-f^\ast\leq\frac{R^2}{2A_k}
  \leq\frac{R^2}{2A_1}\exp\ls-\frac{k-1}{3B}\rs.
\end{gather*}
This proves~\eqref{eq:app_ahpe_far_rate}.

For every $A\geq A_{\tr}$,
\begin{equation*}
  \sqrt{Q(A)}=\max\lb(40M_0RA)^{1/3},\sqrt{8\delta A}\rb.
\end{equation*}
Hence, $Q(A)^{3/2}\geq40M_0RA$ and $Q(A)\geq8\delta A$, so
\begin{equation*}
  \frac{20M_0R^3}{Q(A)^{3/2}}\leq\frac{R^2}{2A},
  \qquad
  \frac{4\delta R^2}{Q(A)}\leq\frac{R^2}{2A}.
\end{equation*}
At least one of these inequalities is an equality, because $\sqrt{Q(A)}$ equals
at least one of the two terms defining it. Therefore,
\begin{equation}
  \frac{R^2}{2A}
  =\max\lb
    \frac{20M_0R^3}{Q(A)^{3/2}},\frac{4\delta R^2}{Q(A)}
  \rb,
  \qquad A\geq A_{\tr}.
  \label{eq:ahpe_error_scale}
\end{equation}
For $A=A_{K+n}$,
combine~\eqref{eq:ahpe_error_certificate},~\eqref{eq:ahpe_error_scale},
and~\eqref{eq:ahpe_cumulative_growth} to obtain
\begin{gather*}
  f(w_{K+n})-f^\ast\leq\frac{R^2}{2A_{K+n}}
  =\max\lb
    \frac{20M_0R^3}{Q(A_{K+n})^{3/2}},
    \frac{4\delta R^2}{Q(A_{K+n})}
  \rb\\
  \leq\max\lb
    \frac{20M_0R^3}{(B+n/9)^3},
    \frac{4\delta R^2}{(B+n/9)^2}
  \rb.
\end{gather*}
This proves~\eqref{eq:app_ahpe_near_rate}.

It remains to bound the iteration and backtracking counts.
Fix $\e>0$.
Set
\begin{equation*}
  A_{\e}\eqdef\frac{R^2}{2\e},\qquad
  N\eqdef\min\lb k\geq1:A_k\geq A_{\e}\rb.
\end{equation*}
By~\eqref{eq:ahpe_cumulative_growth} and~\eqref{eq:ahpe_error_scale},
$\frac{R^2}{2A_{K+n}}\to0$ as $n\to+\infty$.
Thus, $N$ is finite, and~\eqref{eq:ahpe_error_certificate} gives
$f(w_N)-f^\ast\leq\e$.

First, suppose that $\e\geq\e_{\tr}$.
Then $A_{\e}\leq A_{\tr}$, so every iteration with
$A_k<A_{\e}$ satisfies~\eqref{eq:ahpe_geometric_growth}.
The same argument used to prove~\eqref{eq:app_ahpe_transition_index},
with $A_{\e}$ in place of $A_{\tr}$, gives
\begin{equation*}
  N\leq1+\left\lceil3B\log_+\frac{A_{\e}}{A_1}\right\rceil
  =1+\left\lceil3B\log_+\frac{R^2}{2A_1\e}\right\rceil.
\end{equation*}
Now suppose that $0<\e<\e_{\tr}$. Set
\begin{equation*}
  n\eqdef\left\lceil9\max\lb
    \ls\frac{20M_0R^3}{\e}\rs^{1/3},
    2\sqrt{\frac{\delta R^2}{\e}}
  \rb\right\rceil.
\end{equation*}
Equations~\eqref{eq:ahpe_error_scale} and~\eqref{eq:ahpe_cumulative_growth}
give
\begin{equation*}
  \frac{R^2}{2A_{K+n}}
  \leq\max\lb
    \frac{20M_0R^3}{(B+n/9)^3},
    \frac{4\delta R^2}{(B+n/9)^2}
  \rb
  \leq\e.
\end{equation*}
Hence $A_{K+n}\geq A_{\e}$, and $N\leq K+n$.
Using~\eqref{eq:app_ahpe_transition_index}, we obtain
\begin{equation*}
  N\leq1+
  \left\lceil3B\log_+\frac{A_{\tr}}{A_1}\right\rceil
  +\left\lceil9\max\lb
    \ls\frac{20M_0R^3}{\e}\rs^{1/3},
    2\sqrt{\frac{\delta R^2}{\e}}
  \rb\right\rceil.
\end{equation*}
Combining the two cases and using~\eqref{eq:app_ahpe_transition_error} and $B\geq1$
proves~\eqref{eq:app_ahpe_iteration_complexity}.

We now bound the number of rejected trials before producing $w_N$.
If $N=1$, no main trial is needed.
If $N\geq2$, the main iterations needed to compute $w_N$ are
$1,\ldots,N-1$. The minimality of $N$ gives $A_{N-1}<A_{\e}$.
Since $Q$ is nondecreasing, Lemma~\ref{lem:ahpe_backtracking} yields
\begin{equation*}
  J_{N-1}\leq\left\lceil\log_2Q(A_{N-1})\right\rceil
  \leq\left\lceil\log_2Q(A_{\e})\right\rceil.
\end{equation*}
Substituting $A_{\e}=\frac{R^2}{2\e}$
into~\eqref{eq:ahpe_safe_threshold} gives
\begin{equation*}
  Q(A_{\e})
  =\max\lb
    B^2,\frac{4\delta R^2}{\e},
    \ls\frac{20M_0R^3}{\e}\rs^{2/3}
  \rb.
\end{equation*}
This proves the rejection bound~\eqref{eq:app_ahpe_accuracy_rejections}.
By Lemma~\ref{lem:pd_ms_backtracking}, initialization rejects
at most
\begin{equation*}
  \left\lceil
    \log_2\max\lb1,\frac{\overline\tau_0}{\tau_0}\rb
  \right\rceil
\end{equation*}
trials, where $\overline\tau_0 = \overline\tau (x_0)$ is defined
in~\eqref{eq:backtrack_upper_bound}.
Substituting~\eqref{eq:backtrack_upper_bound} gives
the initialization bound~\eqref{eq:app_ahpe_initialization_cost}
and its counterpart under~\eqref{as:lh}.

\end{proof}

\subsection{Oracle Complexity}

\begin{corollary}
\label{cor:ahpe_oracle_complexity}
Under the assumptions of Theorem~\ref{thm:ahpe_convergence_appendix},
fix $\e>0$ and let $N$ be the first index with $A_N\geq R^2/(2\e)$.
Let $\overline\tau_0\eqdef\overline\tau(x_0)$, where $\overline\tau(x)$ is defined in~\eqref{eq:backtrack_upper_bound}.
Then Algorithm~\ref{alg:adaptive_ahpe} computes $w_N$ using at most $N$ Hessian evaluations and at most
\begin{equation*}
  2N+3+\log_{2,+}\frac{\overline\tau_0}{\tau_0}
  +2\log_2\max\lb B^2,\frac{4\delta R^2}{\e},\ls\frac{20M_0R^3}{\e}\rs^{2/3}\rb
\end{equation*}
gradient evaluations.
\end{corollary}
\begin{proof}
Initialization uses one Hessian, the gradient at $x_0$, and at most one gradient per trial.
By Lemma~\ref{lem:pd_ms_backtracking}, it has at most
$2+\log_{2,+}\frac{\overline\tau_0}{\tau_0}$ trials.
Each of the $N-1$ main iterations uses one Hessian $H(w_k)$ and at most two gradients per trial,
$\grad f(v_k)$ and $\grad f(v_k+s_k)$, since $v_k$ changes after each rejection.
These iterations have $N-1$ accepted trials, and~\eqref{eq:app_ahpe_accuracy_rejections}
bounds the number of rejected ones.
The updates reuse the gradient at the accepted point.
Hence, the number of Hessian evaluations is at most $1+(N-1)=N$, and the number of
gradient evaluations is at most
\begin{gather*}
  1+\ls2+\log_{2,+}\frac{\overline\tau_0}{\tau_0}\rs
  +2\ls N-1+\left\lceil\log_2\max\lb B^2,\frac{4\delta R^2}{\e},\ls\frac{20M_0R^3}{\e}\rs^{2/3}\rb\right\rceil\rs\\
  \leq2N+3+\log_{2,+}\frac{\overline\tau_0}{\tau_0}
  +2\log_2\max\lb B^2,\frac{4\delta R^2}{\e},\ls\frac{20M_0R^3}{\e}\rs^{2/3}\rb.
\end{gather*}
\end{proof}

\section{Proof of Algorithm~\ref{alg:ahpe_damped}}
\label{app:ahpe_damped}

We restate Algorithm~\ref{alg:ahpe_damped}.

\begin{algorithm}[!htb]
  \caption{Adaptive damped A-NPE with an inexact Hessian}
  \label{alg:ahpe_damped_appendix}
  \begin{algorithmic}[1]
    \STATE \textbf{Input:} $x_0\in\R^d$ and $\eta_0>0$.
    \STATE $(w_1,\eta_1)\eqdef\MSBacktrack(x_0,\eta_0)$.\label{line:ahpe_damped_appendix_initialization}
    \STATE Set
      \begin{equation}
        A_1\eqdef\frac{1}{\eta_1},
        \qquad
        z_1\eqdef x_0-A_1\grad f(w_1).
        \label{eq:ahpe_damped_initialization}
      \end{equation}
    \FOR{$k\geq1$}
      \STATE Set
        \begin{equation}
          a'_k\eqdef\frac{1+\sqrt{1+4\eta_kA_k}}{2\eta_k},
          \qquad
          A'_k\eqdef A_k+a'_k,
          \qquad
          v_k\eqdef\frac{A_kw_k+a'_kz_k}{A'_k}.
          \label{eq:ahpe_damped_trial_weights}
        \end{equation}
      \STATE $(y_k,\tau_k)\eqdef\MSBacktrack(v_k,\eta_k)$.
      \STATE Set
        \begin{equation}
          \gamma_k\eqdef\frac{\eta_k}{\tau_k},
          \qquad
          a_{k+1}\eqdef\gamma_ka'_k,
          \qquad
          A_{k+1}\eqdef A_k+a_{k+1}.
          \label{eq:ahpe_damped_accepted_weights}
        \end{equation}
      \STATE Set
        \begin{equation}
          w_{k+1}\eqdef\frac{(1-\gamma_k)A_kw_k+\gamma_kA'_ky_k}{A_{k+1}},
          \qquad
          z_{k+1}\eqdef z_k-a_{k+1}\grad f(y_k).
          \label{eq:ahpe_damped_iterate_updates}
        \end{equation}
      \IF{$\tau_k=\eta_k$}
        \STATE Set $\eta_{k+1}\eqdef\frac{\eta_k}{2}$.
      \ELSE
        \STATE Set $\eta_{k+1}\eqdef2\eta_k$.
      \ENDIF
    \ENDFOR
  \end{algorithmic}
\end{algorithm}

\subsection{Initialization and Backtracking}

Lemma~\ref{lem:pd_ms_backtracking} guarantees termination of
all calls to $\MSBacktrack$ in Algorithm~\ref{alg:ahpe_damped_appendix}.

\begin{lemma}
\label{lem:ahpe_damped_initialization}
Suppose that $f$ is convex,~\eqref{as:inexact_hessian} holds, and
either~\eqref{as:gh} or~\eqref{as:lh} holds. Initialization in
Algorithm~\ref{alg:ahpe_damped_appendix} (line~\ref{line:ahpe_damped_appendix_initialization}) uses at most
\begin{equation*}
  1+\left\lceil
    \log_2\max\lb1,\frac{\overline\tau_0}{\eta_0}\rb
  \right\rceil
\end{equation*}
trials, where $\overline\tau_0 = \overline\tau(x_0)$ from~\eqref{eq:backtrack_upper_bound}, and accepted parameter satisfies
\begin{equation*}
  \eta_0\leq\eta_1\leq\max\lb\eta_0,2\overline\tau_0\rb.
\end{equation*} 
The algorithm either returns a minimizer or produces an initial state
$A_1,w_1,z_1,\eta_1$ satisfying $A_1\eta_1=1$ and
\begin{equation}
  A_1(f(w_1)-f^\ast)+\frac12\norm{z_1-x^\ast}^2
  +\frac38\norm{w_1-x_0}^2
  \leq\frac{R^2}{2}.
  \label{eq:ahpe_damped_initialization_bound}
\end{equation}
\end{lemma}

\begin{proof}
If $\grad f(x_0)=0$, the algorithm returns $x_0$ without any trials.
Otherwise,
Lemma~\ref{lem:pd_ms_backtracking} with $x=x_0$,
$\eta=\eta_0$, and $\overline\tau(x_0)=\overline\tau_0$ proves the
bound.
For the initialization~\eqref{eq:ahpe_damped_initialization},
Lemma~\ref{lem:ahpe_energy} applies with $A=0$, $w=z=x_0$,
$a=A_+=A_1=1/\eta_1$, and $\tau=\eta_1$.
The same updates give $y=w_1$ and $z_+=z_1$.
Since $A_1\tau=1$, the lemma yields~\eqref{eq:ahpe_damped_initialization_bound}.
\end{proof}

\subsection{Step Estimates}

\begin{lemma}
\label{lem:ahpe_damped_energy}
Let $f$ be convex.
Every completed main iteration of Algorithm~\ref{alg:ahpe_damped_appendix} satisfies
\begin{equation}
  \ECal_{k+1}\leq\ECal_k-\frac38A'_k\eta_k\norm{s_k}^2,
  \label{eq:ahpe_damped_energy_step}
\end{equation}
where $\ECal_k\eqdef A_k(f(w_k)-f^\ast)+\frac12\norm{z_k-x^\ast}^2$.
\end{lemma}

\begin{proof}
By~\eqref{eq:ahpe_damped_trial_weights} and~\eqref{eq:ahpe_damped_accepted_weights},
$(a'_k)^2=A'_k/\eta_k$ and $\gamma_k\tau_k=\eta_k$, so
\begin{equation*}
  \frac{a_{k+1}^2}{\gamma_kA'_k}=\frac1{\tau_k},
  \qquad
  v_k=\frac{\gamma_kA_kw_k+a_{k+1}z_k}{\gamma_kA'_k}.
\end{equation*}
Apply Lemma~\ref{lem:ahpe_energy} with $A=\gamma_kA_k$,
$a=a_{k+1}$, $A_+=\gamma_kA'_k$, $w=w_k$, $z=z_k$,
and $\tau=\tau_k$. By~\eqref{eq:ahpe_damped_iterate_updates},
$z_+=z_{k+1}$, so it gives
\begin{equation}
  \gamma_kA'_k(f(y_k)-f^\ast)
  +\frac12\norm{z_{k+1}-x^\ast}^2
  \leq\gamma_kA_k(f(w_k)-f^\ast)
  +\frac12\norm{z_k-x^\ast}^2
  -\frac38A'_k\eta_k\norm{s_k}^2.
  \label{eq:ahpe_damped_scaled_energy}
\end{equation}
Since $0<\gamma_k\leq1$ and
\begin{equation*}
  \frac{(1-\gamma_k)A_k}{A_{k+1}}
  +\frac{\gamma_kA'_k}{A_{k+1}}
  =\frac{A_k+\gamma_ka'_k}{A_{k+1}}=1,
\end{equation*}
convexity of $f$ and~\eqref{eq:ahpe_damped_iterate_updates} give
\begin{equation}
  A_{k+1}(f(w_{k+1})-f^\ast)
  \leq(1-\gamma_k)A_k(f(w_k)-f^\ast)
  +\gamma_kA'_k(f(y_k)-f^\ast).
  \label{eq:ahpe_damped_convexity}
\end{equation}
Using~\eqref{eq:ahpe_damped_convexity} and then~\eqref{eq:ahpe_damped_scaled_energy}, we obtain
\begin{gather*}
  \ECal_{k+1}
  =A_{k+1}(f(w_{k+1})-f^\ast)
  +\frac12\norm{z_{k+1}-x^\ast}^2\\
  \leq(1-\gamma_k)A_k(f(w_k)-f^\ast)
  +\gamma_kA'_k(f(y_k)-f^\ast)
  +\frac12\norm{z_{k+1}-x^\ast}^2\\
  \leq(1-\gamma_k)A_k(f(w_k)-f^\ast)
  +\gamma_kA_k(f(w_k)-f^\ast)
  +\frac12\norm{z_k-x^\ast}^2
  -\frac38A'_k\eta_k\norm{s_k}^2\\
  =\ECal_k-\frac38A'_k\eta_k\norm{s_k}^2.
\end{gather*}
\end{proof}

Summing~\eqref{eq:ahpe_damped_energy_step} over $i=1,\ldots,k-1$
and using~\eqref{eq:ahpe_damped_initialization_bound}, we obtain
\begin{equation}
  \ECal_k+\frac38\sum_{i=1}^{k-1}A'_i\eta_i\norm{s_i}^2
  \leq\ECal_1\leq\frac{R^2}{2}.
  \label{eq:ahpe_damped_energy_budget}
\end{equation}
Since all terms on the left-hand side of~\eqref{eq:ahpe_damped_energy_budget}
are nonnegative,
\begin{equation*}
  f(w_k)-f^\ast\leq\frac{R^2}{2A_k},
  \qquad
  \sum_{i=1}^{k-1}A'_i\eta_i\norm{s_i}^2\leq\frac43R^2.
\end{equation*}

For every completed main iteration,~\eqref{eq:ahpe_damped_energy_step}
and~\eqref{eq:ahpe_damped_energy_budget} give
\begin{equation*}
  \frac38A'_k\eta_k\norm{s_k}^2
  \leq\ECal_k-\ECal_{k+1}
  \leq\ECal_k\leq\frac{R^2}{2}.
\end{equation*}
The MS test~\eqref{eq:app_test} also gives
\begin{equation*}
  \norm{\grad f(y_k)}
  \leq\norm{\grad f(y_k)+\tau_ks_k}+\tau_k\norm{s_k}
  \leq\frac32\tau_k\norm{s_k}.
\end{equation*}
Using $(a'_k)^2=A'_k/\eta_k$, we obtain
\begin{equation}
  \norm{s_k}\leq\frac{2R}{\sqrt{3A'_k\eta_k}}
  =\frac{2a'_k}{\sqrt3A'_k}R,
  \qquad
  \norm{\grad f(y_k)}\leq\frac32\tau_k\norm{s_k}.
  \label{eq:ahpe_damped_certified_bounds}
\end{equation}

\begin{lemma}
\label{lem:ahpe_damped_localization}
For every $k\geq1$, the points $w_k,z_k,v_k$ generated by
Algorithm~\ref{alg:ahpe_damped_appendix} satisfy
\begin{equation}
  \norm{z_k-x^\ast}\leq R,
  \label{eq:ahpe_damped_z_bound}
\end{equation}
\begin{equation}
  \norm{w_k-x^\ast}\leq\ls1+\frac2{\sqrt3}\rs R,
  \qquad
  \norm{v_k-x^\ast}\leq\ls1+\frac2{\sqrt3}\rs R.
  \label{eq:ahpe_damped_localization}
\end{equation}
\end{lemma}

\begin{proof}
By~\eqref{eq:ahpe_damped_energy_budget},
\begin{equation*}
  \frac12\norm{z_k-x^\ast}^2\leq\ECal_k\leq\frac{R^2}{2},
\end{equation*}
which gives~\eqref{eq:ahpe_damped_z_bound}.
We prove the bound on $\norm{w_k-x^\ast}$ by induction.
The initialization bound~\eqref{eq:ahpe_damped_initialization_bound} gives
\begin{equation*}
  \norm{w_1-x^\ast}
  \leq\norm{w_1-x_0}+\norm{x_0-x^\ast}
  \leq\ls1+\frac2{\sqrt3}\rs R.
\end{equation*}
Suppose the bound holds for $w_k$.
Substituting $y_k=v_k+s_k$ and
$A'_kv_k=A_kw_k+a'_kz_k$ into~\eqref{eq:ahpe_damped_iterate_updates} yields
\begin{equation*}
  w_{k+1}
  =\frac{(1-\gamma_k)A_kw_k+\gamma_kA'_k(v_k+s_k)}{A_{k+1}}
  =\frac{A_k}{A_{k+1}}w_k
  +\frac{a_{k+1}}{A_{k+1}}z_k
  +\frac{\gamma_kA'_k}{A_{k+1}}s_k.
\end{equation*}
Using~\eqref{eq:ahpe_damped_certified_bounds},~\eqref{eq:ahpe_damped_z_bound},
and $\gamma_ka'_k=a_{k+1}$, we obtain
\begin{gather*}
  \norm{w_{k+1}-x^\ast}
  \leq\frac{A_k}{A_{k+1}}\norm{w_k-x^\ast}
  +\frac{a_{k+1}}{A_{k+1}}\norm{z_k-x^\ast}
  +\frac{\gamma_kA'_k}{A_{k+1}}\norm{s_k}\\
  \leq\frac{A_k}{A_{k+1}}\ls1+\frac2{\sqrt3}\rs R
  +\frac{a_{k+1}}{A_{k+1}}R
  +\frac{2a_{k+1}}{\sqrt3A_{k+1}}R
  =\ls1+\frac2{\sqrt3}\rs R.
\end{gather*}
This completes the induction. Finally,
\begin{equation*}
  \norm{v_k-x^\ast}
  \leq\frac{A_k}{A'_k}\norm{w_k-x^\ast}
  +\frac{a'_k}{A'_k}\norm{z_k-x^\ast}
  \leq\frac{A_k}{A'_k}\ls1+\frac2{\sqrt3}\rs R
  +\frac{a'_k}{A'_k}R
  \leq\ls1+\frac2{\sqrt3}\rs R,
\end{equation*}
which proves~\eqref{eq:ahpe_damped_localization}.
\end{proof}

For iterations of Algorithm~\ref{alg:ahpe_damped_appendix} with $\gamma_k<1$,
apply Lemma~\ref{lem:ahpe_damped_movement} with $x=v_k$, $y=y_k$,
$\eta=\eta_k$, and $\tau=\tau_k$.
If $\tau_k\geq16\delta$,~\eqref{eq:ahpe_damped_movement_lower_bound}
bounds the step length from below.
As in~\citet[Appendix~A]{carmon2022optimal}, we combine this bound
with~\eqref{eq:ahpe_damped_energy_budget} to control the number of these iterations.

\subsection{Proof of Theorem~\ref{thm:ahpe_damped_convergence}}

Define
\begin{equation}
\begin{aligned}
  B&\eqdef
  \begin{cases}
    \max\lb1,\ls6\sqrt{M_1}R\rs^{2/3}\rb,
    & \text{under~\eqref{as:gh}},\\
    \max\lb1,\ls6\sqrt{M_1}R\rs^{2/3},\ls\frac{2R}{\rho}\rs^{2/3}\rb,
    & \text{under~\eqref{as:lh}},
  \end{cases}\\
  \e_{\tr}&\eqdef
  \max\lb\frac{48M_0R^3}{B^{7/2}},\frac{16\delta R^2}{B^2}\rb.
\end{aligned}
  \label{eq:app_ahpe_damped_scales}
\end{equation}
Equivalently, with $r_0$ defined in~\eqref{eq:ahpe_damped_movement_radius},$B=\max\lb1,\ls\frac{R}{r_0}\rs^{2/3}\rb$, where $R/(+\infty)=0$.
For $A>0$, define
\begin{equation}
  \Lambda(A)\eqdef
  \max\lb32\delta,(96M_0R)^{4/7}A^{-3/7},\frac{B^2}{A}\rb.
  \label{eq:ahpe_damped_threshold}
\end{equation}
The function $\Lambda$ is nonincreasing, and $A\Lambda(A)$ is nondecreasing with
$A\Lambda(A)\geq B^2\geq1$.

\begin{lemma}
\label{lem:ahpe_damped_iteration_count}
Suppose that $f$ is convex,~\eqref{as:inexact_hessian} holds, and
either~\eqref{as:gh} or~\eqref{as:lh} holds. Assume that $R>0$ and set
\begin{equation}
  A_{\tr}\eqdef
  \min\lb\frac{B^{7/2}}{96M_0R},\frac{B^2}{32\delta}\rb,
  \label{eq:ahpe_damped_transition_weight}
\end{equation}
where $B^2/(32\delta)=+\infty$ when $\delta=0$.
The weights $A_k$, $k\geq1$, generated by Algorithm~\ref{alg:ahpe_damped_appendix} satisfy
\begin{equation}
  \begin{aligned}
    k-1\leq{} \frac{16}{3}B+2\log_2B
    +\ls6B+\frac32\rs
      \log_+\frac{\min\lb A_k,A_{\tr}\rb}{A_1}
    +35\ls\sqrt{A_k\Lambda(A_k)}-B\rs.
  \end{aligned}
  \label{eq:ahpe_damped_iteration_count}
\end{equation}
\end{lemma}
\begin{proof}
Fix $k\geq1$ and define
\begin{gather*}
  \cS_{k-1}^{=}\eqdef
  \lb i\in\{1,\ldots,k-1\}:\tau_i=\eta_i,
       \ \eta_i\leq2\Lambda(A_i)\rb,\\
  \cS_{k-1}^{>}\eqdef
  \lb i\in\{1,\ldots,k-1\}:\tau_i>\eta_i,
       \ \eta_i\geq\Lambda(A_i)/2\rb.
\end{gather*}
We first bound $|\cS_{k-1}^{=}|$ in terms of $A_k$.
On these iterations, the bound $\eta_i\leq2\Lambda(A_i)$
gives an increase of $A_i$.
For $i\in\cS_{k-1}^{=}$, we have $\gamma_i=1$, $A_{i+1}=A'_i$,
and $(a'_i)^2=A'_i/\eta_i$. Hence,
\begin{equation*}
  \frac{a'_i}{A_i}
  =\sqrt{\frac{A'_i}{\eta_iA_i^2}}
  \geq\frac1{\sqrt{\eta_iA_i}}
  \geq\frac1{\sqrt{2A_i\Lambda(A_i)}}.
\end{equation*}
Using $\log(1+t)\geq t/(1+t)$ for $t\geq0$ gives
\begin{gather*}
  \log\frac{A_{i+1}}{A_i}
  =\log\ls1+\frac{a'_i}{A_i}\rs
  \geq\log\ls1+\frac1{\sqrt{2A_i\Lambda(A_i)}}\rs
  \geq\frac1{\sqrt{2A_i\Lambda(A_i)}+1}
  \geq\frac1{3\sqrt{A_i\Lambda(A_i)}}.
\end{gather*}
Since $\sqrt{t\Lambda(t)}$ is nondecreasing, each $i\in\cS_{k-1}^{=}$
therefore satisfies
\begin{gather*}
  \int_{A_i}^{A_{i+1}}\sqrt{\frac{\Lambda(t)}{t}}dt
  =\int_{A_i}^{A_{i+1}}\frac{\sqrt{t\Lambda(t)}}{t}dt
  \geq\sqrt{A_i\Lambda(A_i)}\log\frac{A_{i+1}}{A_i}\geq\frac13.
\end{gather*}
Summing over $\cS_{k-1}^{=}$ and using the nonnegativity of the integrand,
we obtain
\begin{equation}
    |\cS_{k-1}^{=}|=\sum_{i\in\cS_{k-1}^{=}}1
    \leq3\sum_{i\in\cS_{k-1}^{=}}\int_{A_i}^{A_{i+1}}\sqrt{\frac{\Lambda(t)}{t}}dt 
    \leq3\sum_{i=1}^{k-1}\int_{A_i}^{A_{i+1}}\sqrt{\frac{\Lambda(t)}{t}}dt 
    =3\int_{A_1}^{A_k}\sqrt{\frac{\Lambda(t)}{t}}dt.
  \label{eq:ahpe_damped_down_count}
\end{equation}
We next bound $|\cS_{k-1}^{>}|$ using~\eqref{eq:ahpe_damped_energy_budget}.
We first bound the term $\frac38A'_i\eta_i\norm{s_i}^2$ from below.
For $i\in\cS_{k-1}^{>}$, we have $\tau_i>\eta_i\geq\frac{\Lambda(A_i)}{2}\geq16\delta$. Thus, by Lemma~\ref{lem:ahpe_damped_movement},~\eqref{eq:ahpe_damped_movement_lower_bound} gives
\begin{equation*}
  \norm{s_i}\geq\min\lb r_0,\frac{\tau_i}{48M_0}\rb
  \geq\min\lb r_0,\frac{\Lambda(A_i)}{96M_0}\rb.
\end{equation*}
Since $A'_i\geq A_i$ and $\eta_i\geq\Lambda(A_i)/2$, it follows that
\begin{equation}
\begin{aligned}
  \frac38A'_i\eta_i\norm{s_i}^2
  &\geq\frac3{16}A_i\Lambda(A_i)
  \min\lb r_0^2,\frac{\Lambda(A_i)^2}{(96M_0)^2}\rb\\
  &=\frac{3}{16\sqrt{A_i\Lambda(A_i)}}
  \min\lb \ls A_i\Lambda(A_i)\rs^{3/2}r_0^2,
       \frac{A_i^{3/2}\Lambda(A_i)^{7/2}}{(96M_0)^2}\rb.
\end{aligned}
  \label{eq:ahpe_damped_down_bound}
\end{equation}

We now use the definition of $\Lambda$~\eqref{eq:ahpe_damped_threshold} to bound both terms in the last minimum.
For $r_0<+\infty$, the definition of $B$ and~\eqref{eq:ahpe_damped_threshold} give
\begin{equation*}
  \ls A_i\Lambda(A_i)\rs^{3/2}r_0^2
  \geq B^3r_0^2
  \geq\ls\frac{R}{r_0}\rs^2r_0^2=R^2.
\end{equation*}
Also,~\eqref{eq:ahpe_damped_threshold} gives
\begin{equation*}
  \frac{A_i^{3/2}\Lambda(A_i)^{7/2}}{(96M_0)^2}
  \geq\frac{A_i^{3/2}\ls(96M_0R)^{4/7}A_i^{-3/7}\rs^{7/2}}{(96M_0)^2}
  =R^2.
\end{equation*}
Substituting these bounds in~\eqref{eq:ahpe_damped_down_bound} and using that $A_i\leq A_k$
and $A\Lambda(A)$ is nondecreasing, we obtain
\begin{equation*}
  \frac38A'_i\eta_i\norm{s_i}^2
  \geq\frac{3R^2}{16\sqrt{A_i\Lambda(A_i)}}
  \geq\frac{3R^2}{16\sqrt{A_k\Lambda(A_k)}}.
\end{equation*}
Summing over $i\in\cS_{k-1}^{>}$ and using~\eqref{eq:ahpe_damped_energy_budget}, we obtain
\begin{equation*}
  \frac{3R^2}{16\sqrt{A_k\Lambda(A_k)}}|\cS_{k-1}^{>}|
  =\sum_{i\in\cS_{k-1}^{>}}\frac{3R^2}{16\sqrt{A_k\Lambda(A_k)}}
  \leq\frac38\sum_{i\in\cS_{k-1}^{>}}A'_i\eta_i\norm{s_i}^2
  \leq\frac38\sum_{i=1}^{k-1}A'_i\eta_i\norm{s_i}^2
  \leq\frac{R^2}{2}.
\end{equation*}
Therefore,
\begin{equation}
  |\cS_{k-1}^{>}|\leq\frac83\sqrt{A_k\Lambda(A_k)}.
  \label{eq:ahpe_damped_up_count}
\end{equation}
It remains to bound the number $(k-1)-|\cS_{k-1}^{=}|-|\cS_{k-1}^{>}|$ of iterations outside $\cS_{k-1}^{=}\cup\cS_{k-1}^{>}$.
For these iterations, there are two cases:
\begin{itemize}
  \item First, if $\tau_i=\eta_i$, then $\eta_i>2\Lambda(A_i)$, $\eta_{i+1}=\frac{\eta_i}{2}$.
  \item Second, if $\tau_i>\eta_i$, then $\eta_i<\frac{\Lambda(A_i)}{2}, \eta_{i+1}=2\eta_i$.
\end{itemize}
To measure the change relative to $\Lambda(A_i)$, set
$V_i\eqdef\left|\log\frac{\eta_i}{\Lambda(A_i)}\right|$.
In the first case, $\eta_{i+1}>\Lambda(A_i)$, so
\begin{equation*}
  \left|\log\frac{\eta_{i+1}}{\Lambda(A_i)}\right|
  =\log\frac{\eta_i}{2\Lambda(A_i)}
  =\log\frac{\eta_i}{\Lambda(A_i)}-\log2
  =V_i-\log2.
\end{equation*}
In the second case, $\eta_{i+1}<\Lambda(A_i)$, so
\begin{equation*}
  \left|\log\frac{\eta_{i+1}}{\Lambda(A_i)}\right|
  =-\log\frac{2\eta_i}{\Lambda(A_i)}
  =-\log\frac{\eta_i}{\Lambda(A_i)}-\log2
  =V_i-\log2.
\end{equation*}
For $i\in\cS_{k-1}^{=}\cup\cS_{k-1}^{>}$, the update of $\eta_i$ gives
$\left|\log\frac{\eta_{i+1}}{\Lambda(A_i)}\right|\leq V_i+\log2$.

Since $\Lambda(A_{i+1})\leq\Lambda(A_i)$,
\begin{gather*}
  V_{i+1}
  =\left|\log\frac{\eta_{i+1}}{\Lambda(A_i)}
        +\log\frac{\Lambda(A_i)}{\Lambda(A_{i+1})}\right|
  \leq\left|\log\frac{\eta_{i+1}}{\Lambda(A_i)}\right|
       +\log\frac{\Lambda(A_i)}{\Lambda(A_{i+1})}.
\end{gather*}
Summing over $i=1,\ldots,k-1$ and telescoping the logarithms gives
\begin{gather*}
  V_k-V_1
  \leq\ls|\cS_{k-1}^{=}|+|\cS_{k-1}^{>}|\rs\log2
       -\ls(k-1)-|\cS_{k-1}^{=}|-|\cS_{k-1}^{>}|\rs\log2
       +\sum_{i=1}^{k-1}\log\frac{\Lambda(A_i)}{\Lambda(A_{i+1})}\\
  =\ls2\ls|\cS_{k-1}^{=}|+|\cS_{k-1}^{>}|\rs-(k-1)\rs\log2
       +\log\frac{\Lambda(A_1)}{\Lambda(A_k)}.
\end{gather*}
Initialization gives $\eta_1=1/A_1$. Since $A_1\Lambda(A_1)\geq1$,
\begin{equation*}
  V_1=\left|\log\frac1{A_1\Lambda(A_1)}\right|
  =\log\ls A_1\Lambda(A_1)\rs.
\end{equation*}
Also, monotonicity of $A\Lambda(A)$ implies
\begin{equation*}
  \log\frac{\Lambda(A_1)}{\Lambda(A_k)}
  =\log\frac{A_k}{A_1}
   +\log\frac{A_1\Lambda(A_1)}{A_k\Lambda(A_k)}
  \leq\log\frac{A_k}{A_1}.
\end{equation*}
Using $V_k\geq0$, combining the logarithms, and dividing by $\log2$, we obtain
\begin{equation}
  k-1\leq2\ls|\cS_{k-1}^{=}|+|\cS_{k-1}^{>}|\rs
  +\log_2\ls A_k\Lambda(A_1)\rs.
  \label{eq:ahpe_damped_iteration_balance}
\end{equation}
Substituting~\eqref{eq:ahpe_damped_down_count} and~\eqref{eq:ahpe_damped_up_count}
into~\eqref{eq:ahpe_damped_iteration_balance} gives
\begin{equation}
  k-1\leq6\int_{A_1}^{A_k}\sqrt{\frac{\Lambda(t)}{t}}dt
  +\frac{16}{3}\sqrt{A_k\Lambda(A_k)}
  +\log_2\ls A_k\Lambda(A_1)\rs.
  \label{eq:ahpe_damped_iteration_integral}
\end{equation}
We now estimate the integral and logarithmic term in~\eqref{eq:ahpe_damped_iteration_integral}.
By~\eqref{eq:ahpe_damped_threshold} and~\eqref{eq:ahpe_damped_transition_weight},
\begin{equation*}
  \sqrt{A\Lambda(A)}=
  \begin{cases}
    B, & 0<A\leq A_{\tr},\\
    \max\lb(96M_0RA)^{2/7},\sqrt{32\delta A}\rb,
       & A>A_{\tr}.
  \end{cases}
\end{equation*}
Consider an interval $[u,v]\subset[A_{\tr},+\infty)$ on which the same term attains the maximum.
If $\sqrt{t\Lambda(t)}=(96M_0Rt)^{2/7}$ on this interval, then
\begin{gather*}
  \int_u^v\sqrt{\frac{\Lambda(t)}t}dt
  =(96M_0R)^{2/7}\int_u^v t^{-5/7}dt
  =\frac72\ls\sqrt{v\Lambda(v)}-\sqrt{u\Lambda(u)}\rs,\quad 
  \log\frac vu=\frac72\log\frac{\sqrt{v\Lambda(v)}}{\sqrt{u\Lambda(u)}}.
\end{gather*}
If $\sqrt{t\Lambda(t)}=\sqrt{32\delta t}$ on this interval, then
\begin{gather*}
  \int_u^v\sqrt{\frac{\Lambda(t)}t}dt
  =\sqrt{32\delta}\int_u^v t^{-1/2}dt
  =2\ls\sqrt{v\Lambda(v)}-\sqrt{u\Lambda(u)}\rs,\quad 
  \log\frac vu=2\log\frac{\sqrt{v\Lambda(v)}}{\sqrt{u\Lambda(u)}}.
\end{gather*}
Both coefficients are at most $7/2$.
Summing over the intervals above $A_{\tr}$, the intermediate endpoint terms cancel in both bounds.
Below $A_{\tr}$, we have $\sqrt{t\Lambda(t)}=B$ and integrate $B/t$.
Thus, for every $A\geq A_1$,
\begin{gather}
  \int_{A_1}^{A}\sqrt{\frac{\Lambda(t)}{t}}dt
  \leq B\log_+\frac{\min\lb A,A_{\tr}\rb}{A_1}
       +\frac72\ls\sqrt{A\Lambda(A)}-B\rs,
  \label{eq:ahpe_damped_clipped_integral}\\
  \log\frac{A}{A_1}
  \leq\log_+\frac{\min\lb A,A_{\tr}\rb}{A_1}
       +\frac72\log\frac{\sqrt{A\Lambda(A)}}B.
  \label{eq:ahpe_damped_clipped_log}
\end{gather}
These bounds also hold when $A_1>A_{\tr}$, since $\sqrt{A_1\Lambda(A_1)}\geq B$.
Using $\log_2t\leq\frac32\log t$ for $t\geq1$,~\eqref{eq:ahpe_damped_clipped_log},
and the monotonicity of $A\Lambda(A)$, we obtain
\begin{gather}
    \log_2\ls A_k\Lambda(A_1)\rs
    =\log_2\frac{A_k}{A_1}+\log_2\ls A_1\Lambda(A_1)\rs
    \leq\frac32\log\frac{A_k}{A_1}
         +2\log_2B+3\log\frac{\sqrt{A_k\Lambda(A_k)}}B \notag \\
    \leq2\log_2B+\frac32\log_+\frac{\min\lb A_k,A_{\tr}\rb}{A_1}
         +\frac{33}{4}\log\frac{\sqrt{A_k\Lambda(A_k)}}B.
  \label{eq:ahpe_damped_clipped_log_term}
\end{gather}
Substituting~\eqref{eq:ahpe_damped_clipped_integral} and~\eqref{eq:ahpe_damped_clipped_log_term}
into~\eqref{eq:ahpe_damped_iteration_integral}, we obtain
\begin{gather*}
  k-1\leq\frac{16}{3}B+2\log_2B
       +\ls6B+\frac32\rs\log_+\frac{\min\lb A_k,A_{\tr}\rb}{A_1}\\
  +\frac{79}{3}\ls\sqrt{A_k\Lambda(A_k)}-B\rs
       +\frac{33}{4}\log\frac{\sqrt{A_k\Lambda(A_k)}}B.
\end{gather*}
Finally, since $\sqrt{A_k\Lambda(A_k)}\geq B\geq1$,
\begin{equation*}
  \log\frac{\sqrt{A_k\Lambda(A_k)}}B
  \leq\frac{\sqrt{A_k\Lambda(A_k)}-B}{B}, \quad 
  \frac{79}{3}+\frac{33}{4B}
  \leq\frac{415}{12}<35.
\end{equation*}
Substituting these inequalities proves~\eqref{eq:ahpe_damped_iteration_count}.
\end{proof}

We restate Theorem~\ref{thm:ahpe_damped_convergence} with explicit constants.

\begin{theorem}
\label{thm:ahpe_damped_convergence_appendix}
Suppose that $f$ is convex,~\eqref{as:inexact_hessian} holds, and
either~\eqref{as:gh} or~\eqref{as:lh} holds.
Let $\overline\tau_0\eqdef\overline\tau(x_0)$, where $\overline\tau(x)$ is defined in~\eqref{eq:backtrack_upper_bound}.
For every $\e>0$, the method generates an iterate $w_N$ satisfying
$f(w_N)-f^\ast\leq\e$ within
\begin{equation}
\begin{aligned}
  N\leq{}&2+\frac{16}{3}B+2\log_2B
    +\ls6B+\frac32\rs\log_+\frac{\max\lb\eta_0,\overline\tau_0\rb R^2}{\max\lb\e,\e_{\tr}\rb}\\
    &+35\ls\max\lb B,\ls\frac{48M_0R^3}{\e}\rs^{2/7},
      4\sqrt{\frac{\delta R^2}{\e}}\rb-B\rs
\end{aligned}
  \label{eq:app_ahpe_damped_outer_complexity}
\end{equation}
outer iterations, including initialization, each using at most one Hessian evaluation.
Moreover, set
\begin{equation}
  K\eqdef2+\left\lfloor
    \frac{16}{3}B+2\log_2B
    +\ls6B+\frac32\rs\log_+\frac{\max\lb\eta_0,\overline\tau_0\rb R^2}{\e_{\tr}}
  \right\rfloor.
  \label{eq:app_ahpe_damped_transition_index}
\end{equation}
For this $K$,
\begin{gather}
  f(w_k)-f^\ast
  \leq\max\lb\eta_0,\overline\tau_0\rb R^2
    \exp\ls-\frac{\max\lb0,k-1-\frac{16}{3}B-2\log_2B\rb}{6B+\frac32}\rs,
  \qquad 1\leq k<K,
  \label{eq:app_ahpe_damped_far_rate}\\
  f(w_{K+n})-f^\ast
  \leq\max\lb\frac{48M_0R^3}{(B+n/35)^{7/2}},
    \frac{16\delta R^2}{(B+n/35)^2}\rb,
  \qquad n=0,1,\ldots.
  \label{eq:app_ahpe_damped_near_rate}
\end{gather}
\end{theorem}

\begin{proof}
Lemmas~\ref{lem:pd_ms_backtracking} and~\ref{lem:ahpe_damped_initialization}
show that all backtracking loops terminate.
If initialization returns a minimizer, all claims are immediate.
By~\eqref{eq:ahpe_damped_energy_budget},
\begin{equation}
  f(w_k)-f^\ast\leq\frac{\ECal_k}{A_k}\leq\frac{R^2}{2A_k}.
  \label{eq:ahpe_damped_error_certificate}
\end{equation}
Since $a_{k+1}>0$, the $A_k$ strictly increase.

Fix $\e>0$ and set $A_\e\eqdef R^2/(2\e)$.
If $A_1\geq A_\e$, initialization suffices by~\eqref{eq:ahpe_damped_error_certificate}.
Otherwise, by~\eqref{eq:ahpe_damped_iteration_count}, unless a minimizer is returned earlier,
there exists a first index $N$ with $A_N\geq A_\e$.
Then $N\geq2$, $A_{N-1}<A_\e\leq A_N$, and~\eqref{eq:ahpe_damped_error_certificate}
gives $f(w_N)-f^\ast\leq\e$.
Applying~\eqref{eq:ahpe_damped_iteration_count} at $N-1$
and using monotonicity of its right-hand side, we obtain
\begin{gather*}
  N\leq2+\frac{16}{3}B+2\log_2B
  +\ls6B+\frac32\rs\log_+\frac{\min\lb A_\e,A_{\tr}\rb}{A_1}
  +35\ls\sqrt{A_\e\Lambda(A_\e)}-B\rs.
\end{gather*}
By~\eqref{eq:app_ahpe_damped_scales} and~\eqref{eq:ahpe_damped_transition_weight},
$A_{\tr}=R^2/(2\e_{\tr})$. Hence,
\begin{gather*}
  \log_+\frac{\min\lb A_\e,A_{\tr}\rb}{A_1}
  =\log_+\frac{R^2}{2A_1\max\lb\e,\e_{\tr}\rb},\\
  \sqrt{A_\e\Lambda(A_\e)}
  =\max\lb B,\ls\frac{48M_0R^3}{\e}\rs^{2/7},
    4\sqrt{\frac{\delta R^2}{\e}}\rb.
\end{gather*}
Next, by Lemma~\ref{lem:ahpe_damped_initialization} and~\eqref{eq:ahpe_damped_initialization} we have 
\begin{equation}
  \frac{1}{A_1} = \eta_1 \leq \max\lb\eta_0, 2 \overline{\tau}_0\rb \leq 2\max\lb\eta_0,\overline\tau_0\rb.
  \label{eq:ahpe_damped_initial_weight}
\end{equation}
Hence,
\begin{equation}
  \frac{R^2}{2A_1}\leq\max\lb\eta_0,\overline\tau_0\rb R^2,
  \qquad
  2\max\lb\eta_0,\overline\tau_0\rb A_{\tr}=\frac{\max\lb\eta_0,\overline\tau_0\rb R^2}{\e_{\tr}},
  \label{eq:ahpe_damped_initial_scales}
\end{equation}
and
\begin{equation*}
  \log_+\frac{R^2}{2A_1\max\lb\e,\e_{\tr}\rb}
  \leq\log_+\frac{\max\lb\eta_0,\overline\tau_0\rb R^2}{\max\lb\e,\e_{\tr}\rb}.
\end{equation*}
Substituting these estimates proves~\eqref{eq:app_ahpe_damped_outer_complexity}.
If a zero gradient is queried before reaching $A_\e$, let $N$ count the
outer iterations up to that return. The last completed weight still satisfies
$A_{N-1}<A_\e$, so the same estimate applies, and the returned point is a minimizer.
There is one Hessian evaluation during initialization and at most one on each
main iteration, since all inner trials reuse that Hessian.

We now prove~\eqref{eq:app_ahpe_damped_far_rate} and~\eqref{eq:app_ahpe_damped_near_rate}
for iterates produced before a zero-gradient return.
The right-hand side of~\eqref{eq:ahpe_damped_iteration_count} is nonincreasing in $A_1$,
and $1/A_1\leq2\max\lb\eta_0,\overline\tau_0\rb$ by~\eqref{eq:ahpe_damped_initial_weight}.
Hence,~\eqref{eq:ahpe_damped_iteration_count} remains valid with $1/A_1$ replaced by $2\max\lb\eta_0,\overline\tau_0\rb$.
First, let $1\leq k<K$. If $A_k\leq A_{\tr}$, then
$\sqrt{A_k\Lambda(A_k)}=B$, and~\eqref{eq:ahpe_damped_iteration_count} gives
\begin{equation*}
  k-1\leq\frac{16}{3}B+2\log_2B
  +\ls6B+\frac32\rs\log\ls2\max\lb\eta_0,\overline\tau_0\rb A_k\rs.
\end{equation*}
If $A_k>A_{\tr}$, then~\eqref{eq:app_ahpe_damped_transition_index},~\eqref{eq:ahpe_damped_initial_scales},
and $2\max\lb\eta_0,\overline\tau_0\rb A_k\geq1$ give
\begin{gather*}
  k-1\leq K-2
  \leq\frac{16}{3}B+2\log_2B
    +\ls6B+\frac32\rs\log_+\ls2\max\lb\eta_0,\overline\tau_0\rb A_{\tr}\rs\\
  \leq\frac{16}{3}B+2\log_2B
    +\ls6B+\frac32\rs\log\ls2\max\lb\eta_0,\overline\tau_0\rb A_k\rs.
\end{gather*}
Thus, in both cases, since $\log\ls2\max\lb\eta_0,\overline\tau_0\rb A_k\rs\geq0$,
\begin{equation*}
  \log\ls2\max\lb\eta_0,\overline\tau_0\rb A_k\rs
  \geq\frac{\max\lb0,k-1-\frac{16}{3}B-2\log_2B\rb}{6B+\frac32}.
\end{equation*}
Combining this bound with~\eqref{eq:ahpe_damped_error_certificate}
proves~\eqref{eq:app_ahpe_damped_far_rate}.

Next, let $k=K+n$ with $n\geq0$.
By~\eqref{eq:app_ahpe_damped_transition_index} and~\eqref{eq:ahpe_damped_initial_scales},
\begin{equation*}
  \frac{16}{3}B+2\log_2B
  +\ls6B+\frac32\rs\log_+\ls2\max\lb\eta_0,\overline\tau_0\rb A_{\tr}\rs+n<k-1.
\end{equation*}
Since $\min\lb A_k,A_{\tr}\rb\leq A_{\tr}$,~\eqref{eq:ahpe_damped_iteration_count} also gives
\begin{equation*}
  k-1\leq\frac{16}{3}B+2\log_2B
  +\ls6B+\frac32\rs\log_+\ls2\max\lb\eta_0,\overline\tau_0\rb A_{\tr}\rs
  +35\ls\sqrt{A_k\Lambda(A_k)}-B\rs.
\end{equation*}
Combining these inequalities yields
\begin{equation*}
  \sqrt{A_k\Lambda(A_k)}>B+\frac n{35}\geq B.
\end{equation*}
Therefore, by~\eqref{eq:ahpe_damped_threshold},
\begin{gather*}
  \max\lb(96M_0RA_k)^{2/7},\sqrt{32\delta A_k}\rb>B+\frac n{35},\quad
  \frac1{A_k}\leq
  \max\lb\frac{96M_0R}{(B+n/35)^{7/2}},
    \frac{32\delta}{(B+n/35)^2}\rb.
\end{gather*}
Multiplying the last inequality by $R^2/2$ and using~\eqref{eq:ahpe_damped_error_certificate}
proves~\eqref{eq:app_ahpe_damped_near_rate}.
After a zero-gradient return, the output sequence consists of minimizers,
so both rate bounds continue to hold.
\end{proof}

\subsection{Oracle Complexity}
\label{sec:ahpe_damped_oracle}

In this section, we bound the numbers of Hessian and gradient evaluations
of Algorithm~\ref{alg:ahpe_damped_appendix}. Let
\begin{equation}
  \overline g\eqdef\max_{\norm{x-x^\ast}\leq(1+2/\sqrt3)R}\norm{\grad f(x)}.
  \label{eq:ahpe_damped_gradient_bound}
\end{equation}
By~\eqref{eq:ahpe_damped_localization}, $\norm{\grad f(v_i)}\leq\overline g$ for any $v_i$.
Recall that $r_0$ is defined in~\eqref{eq:ahpe_damped_movement_radius}.

\begin{corollary}
\label{cor:ahpe_damped_oracle_complexity}
Under the assumptions of Theorem~\ref{thm:ahpe_damped_convergence_appendix},
fix $\e>0$. Let $N$ be the number of outer iterations, including initialization,
until $A_N\geq R^2/(2\e)$ or an earlier zero-gradient return.
Then Algorithm~\ref{alg:ahpe_damped_appendix} uses at most $N$ Hessian evaluations and at most
\begin{equation}
\begin{aligned}
  &7N+\log_{2,+}\frac{\overline\tau_0}{\eta_0}
    +30\sqrt{\frac{2\delta R^2}{\e}}
    +120\ls\frac{2M_0R^3}{\e}\rs^{2/7}\\
  &+8B\log_{2,+}\frac{\delta+\overline g/r_0}{\delta+M_0r_0}
    \ls1+\log_+\frac{4\max\lb\eta_0,\overline\tau_0\rb R^2}{\e}\rs^{2/3}
\end{aligned}
  \label{eq:ahpe_damped_gradient_count}
\end{equation}
gradient evaluations.
\end{corollary}

We prove Corollary~\ref{cor:ahpe_damped_oracle_complexity} after three auxiliary lemmas.
Lemma~\ref{lem:ahpe_damped_rejection_sum} bounds the total number of rejected trials.
Its proof uses the bounds on the sums of $\eta_i^{-1/2}$ and $(\norm{s_i}/\eta_i)^{2/7}$
from Lemma~\ref{lem:ahpe_damped_summability} and the bound on the number of
rejected iterations with $\norm{s_i}>r_0$ from Lemma~\ref{lem:ahpe_damped_long_count}.
By~\eqref{eq:ahpe_damped_energy_budget} and $A_i\leq A'_i$,
\begin{equation}
  \sum_{i=1}^{m}A_i\eta_i\norm{s_i}^2\leq\frac43R^2,\qquad m\geq1.
  \label{eq:ahpe_damped_step_energy}
\end{equation}

\begin{lemma}
\label{lem:ahpe_damped_summability}
Under the assumptions of Theorem~\ref{thm:ahpe_damped_convergence_appendix},
every completed main iteration of Algorithm~\ref{alg:ahpe_damped_appendix} satisfies
\begin{equation}
 \eta_iA_i\geq\frac12,\qquad
 A_{i+1}\leq(2+\sqrt3)A_i<4A_i.
 \label{eq:ahpe_damped_weight_overshoot}
\end{equation}
For every $m\geq1$ and $A=A_{m+1}$,
\begin{gather}
 \sum_{i=1}^{m}\eta_i^{-1/2}\leq5\sqrt A,\qquad
 \sum_{i=1}^{m}\ls\frac{\norm{s_i}}{\eta_i}\rs^{2/7}
 \leq20(RA)^{2/7},
 \label{eq:ahpe_damped_weighted_sums}\\
 \sum_{i=1}^{m}\frac1{\sqrt{\eta_iA_i}}
 \leq8\log6+4\log\frac A{A_1}.
 \label{eq:ahpe_damped_logarithmic_sum}
\end{gather}
\end{lemma}
\begin{proof}
Initialization gives $\eta_1A_1=1$. If $\tau_i>\eta_i$, then
$\eta_{i+1}=2\eta_i$ and $A_{i+1}>A_i$.
If $\tau_i=\eta_i$, then $A_{i+1}=A'_i$, $a'_i\geq1/\eta_i$, and
\begin{equation*}
 \eta_{i+1}A_{i+1}=\frac12\eta_iA'_i
 =\frac12(\eta_i a'_i)^2\geq\frac12.
\end{equation*}
Induction proves the first bound in~\eqref{eq:ahpe_damped_weight_overshoot}.
Since $\eta_iA_i\geq1/2$ and $0<\gamma_i\leq1$,
\begin{equation*}
 \frac{A_{i+1}}{A_i}\leq\frac{A'_i}{A_i}
 =1+\frac{1+\sqrt{1+4\eta_iA_i}}{2\eta_iA_i}
 \leq2+\sqrt3<4.
\end{equation*}
On every iteration with $\tau_i=\eta_i$, we have $a'_i=\sqrt{A_{i+1}/\eta_i}$ and
\begin{equation*}
 \sqrt{A_{i+1}}-\sqrt{A_i}
 =\frac{a'_i}{\sqrt{A_{i+1}}+\sqrt{A_i}}
 =\frac{\sqrt{A_{i+1}}}{\sqrt{\eta_i}\ls\sqrt{A_{i+1}}+\sqrt{A_i}\rs}
 \geq\frac{1}{2\sqrt{\eta_i}}.
\end{equation*}
All weight increments are nonnegative, so
\begin{equation*}
 \sum_{\substack{1\leq i\leq m\\\tau_i=\eta_i}}\eta_i^{-1/2}
 \leq2(\sqrt A-\sqrt{A_1}).
\end{equation*}
The update of $\eta_i$ gives $\eta_{i+1}^{-1/2}=\sqrt2\eta_i^{-1/2}$ when $\tau_i=\eta_i$,
and $\eta_{i+1}^{-1/2}=\eta_i^{-1/2}/\sqrt2$ otherwise.
Telescoping and $\eta_1^{-1/2}=\sqrt{A_1}$ yield
\begin{gather*}
 (1-1/\sqrt2)\sum_{\substack{1\leq i\leq m\\\tau_i>\eta_i}}\eta_i^{-1/2}
 =(\sqrt2-1)\sum_{\substack{1\leq i\leq m\\\tau_i=\eta_i}}\eta_i^{-1/2}
 +\eta_1^{-1/2}-\eta_{m+1}^{-1/2},\\
 \sum_{i=1}^{m}\eta_i^{-1/2}
 \leq(1+\sqrt2)\sum_{\substack{1\leq i\leq m\\\tau_i=\eta_i}}\eta_i^{-1/2}
 +(2+\sqrt2)\sqrt{A_1}
 \leq2(1+\sqrt2)\sqrt A-\sqrt2\sqrt{A_1}
 \leq5\sqrt A.
\end{gather*}
Set $S_0=\sqrt{A_1}$ and $S_i=S_0+\sum_{j=1}^{i}\eta_j^{-1/2}$.
Applying the preceding bound to each shorter prefix gives
\begin{equation}
 S_m\leq6\sqrt A,\qquad
 S_{i-1}\leq6\sqrt{A_i},\qquad
 S_i\leq(6+\sqrt2)\sqrt{A_i}\leq8\sqrt{A_i}.
 \label{eq:ahpe_damped_inclusive_sums}
\end{equation}
Here $\eta_i^{-1/2}\leq\sqrt{2A_i}$ follows from~\eqref{eq:ahpe_damped_weight_overshoot}.
Since $\eta_i^{-1/2}=S_i-S_{i-1}$, integration of the decreasing function $1/t$
yields
\begin{equation*}
 \sum_{i=1}^{m}\frac{1}{\sqrt{\eta_iA_i}}
 \leq8\sum_{i=1}^{m}\frac{S_i-S_{i-1}}{S_i}
 \leq8\int_{S_0}^{S_m}\frac{dt}{t}
 \leq8\log6+4\log\frac A{A_1}.
\end{equation*}
This proves~\eqref{eq:ahpe_damped_logarithmic_sum}.

For the second bound in~\eqref{eq:ahpe_damped_weighted_sums},
H\"older's inequality gives
\begin{gather*}
 \sum_{i=1}^{m}\ls\frac{\norm{s_i}}{\eta_i}\rs^{2/7}
 =\sum_{i=1}^{m}\ls A_i\eta_i\norm{s_i}^2\rs^{1/7}
       \ls\frac{1}{\eta_i^{1/2}A_i^{1/6}}\rs^{6/7}
 \leq\ls\sum_{i=1}^{m}A_i\eta_i\norm{s_i}^2\rs^{1/7}
      \ls\sum_{i=1}^{m}\frac{1}{\eta_i^{1/2}A_i^{1/6}}\rs^{6/7},\\
 \sum_{i=1}^{m}\frac{1}{\eta_i^{1/2}A_i^{1/6}}
 \leq2\sum_{i=1}^{m}\frac{S_i-S_{i-1}}{S_i^{1/3}}
 \leq2\int_{S_0}^{S_m}t^{-1/3}dt
 \leq3S_m^{2/3}
 \leq3\cdot6^{2/3}A^{1/3}.
\end{gather*}
The last display uses~\eqref{eq:ahpe_damped_inclusive_sums}.
Substituting these bounds and~\eqref{eq:ahpe_damped_step_energy} gives
\begin{equation*}
 \sum_{i=1}^{m}\ls\frac{\norm{s_i}}{\eta_i}\rs^{2/7}
 \leq(4R^2/3)^{1/7}
       (3\cdot6^{2/3}A^{1/3})^{6/7}
 \leq20(RA)^{2/7}.
\end{equation*}
\end{proof}

\begin{lemma}
\label{lem:ahpe_damped_long_count}
Under the assumptions of Theorem~\ref{thm:ahpe_damped_convergence_appendix},
for $m\geq1$ and $A=A_{m+1}$, define
\begin{equation*}
 \mathcal I_m\eqdef
 \lb i\in\{1,\ldots,m\}:\tau_i>\eta_i,\ \norm{s_i}>r_0\rb.
\end{equation*}
Then
\begin{equation}
 |\mathcal I_m|\leq8B\ls1+\log\frac A{A_1}\rs^{2/3}.
 \label{eq:ahpe_damped_long_count}
\end{equation}
If $r_0=+\infty$, then $\mathcal I_m$ is empty.
\end{lemma}
\begin{proof}
Let $r_0<+\infty$. By~\eqref{eq:app_ahpe_damped_scales}, $(R/r_0)^{2/3}\leq B$.
H\"older's inequality,~\eqref{eq:ahpe_damped_step_energy},
and~\eqref{eq:ahpe_damped_logarithmic_sum} give
\begin{gather*}
 |\mathcal I_m|
 \leq r_0^{-2/3}\sum_{i\in\mathcal I_m}\norm{s_i}^{2/3}
 =r_0^{-2/3}\sum_{i\in\mathcal I_m}
 \ls A_i\eta_i\norm{s_i}^2\rs^{1/3}\ls\frac{1}{\sqrt{\eta_iA_i}}\rs^{2/3}\\
 \leq r_0^{-2/3}
 \ls\sum_{i=1}^{m}A_i\eta_i\norm{s_i}^2\rs^{1/3}
 \ls\sum_{i=1}^{m}\frac{1}{\sqrt{\eta_iA_i}}\rs^{2/3}
 \leq(4/3)^{1/3}B
 \ls8\log6+4\log\frac A{A_1}\rs^{2/3}.
\end{gather*}
Since $(4/3)^{1/3}(8\log6)^{2/3}<8$ and $8\log6>4$,
this proves~\eqref{eq:ahpe_damped_long_count}.
\end{proof}

\begin{lemma}
\label{lem:ahpe_damped_rejection_sum}
Under the assumptions of Theorem~\ref{thm:ahpe_damped_convergence_appendix},
let $b_i\eqdef\log_2(\tau_i/\eta_i)$ be the number of rejected trials
in main iteration $i$, including failed $\PD$ tests.
For every $m\geq1$ and $A=A_{m+1}$,
\begin{equation}
  \sum_{i=1}^{m}b_i\leq{}5m+30\sqrt{\delta A}+120(M_0RA)^{2/7}
  +8B\log_{2,+}\frac{\delta+\overline g/r_0}{\delta+M_0r_0}
    \ls1+\log\frac{A}{A_1}\rs^{2/3},
  \label{eq:ahpe_damped_rejection_sum}
\end{equation}
where the last term is absent if $r_0=+\infty$.
\end{lemma}
\begin{proof}
If $\tau_i=\eta_i$, then $b_i=0$. Let $\tau_i>\eta_i$.
If $\norm{s_i}\leq r_0$, then Lemma~\ref{lem:ahpe_damped_movement} gives
$\tau_i\leq8\delta+24M_0\norm{s_i}$.
If $\norm{s_i}>r_0$, then~\eqref{as:inexact_hessian}, convexity,
and~\eqref{eq:ahpe_damped_gradient_bound} give
\begin{gather*}
  (\tau_i-\delta)\norm{s_i}^2
  \leq\la(H(v_i)+\tau_iI)s_i,s_i\ra
  =-\la\grad f(v_i),s_i\ra
  \leq\overline g\norm{s_i},\\
  \tau_i\leq\delta+\frac{\overline g}{r_0}
  \leq\max\lb1,\frac{\delta+\overline g/r_0}{\delta+M_0r_0}\rb
  \ls\delta+M_0\norm{s_i}\rs.
\end{gather*}
Since $8\delta+24M_0\norm{s_i}\leq32\eta_i\ls1+\frac{\delta}{\eta_i}\rs\ls1+\frac{M_0\norm{s_i}}{\eta_i}\rs$,
in both cases
\begin{equation*}
  b_i\leq5+\log_2\ls1+\frac{\delta}{\eta_i}\rs
  +\log_2\ls1+\frac{M_0\norm{s_i}}{\eta_i}\rs
  +\mathbf1_{\{i\in\mathcal I_m\}}\log_{2,+}\frac{\delta+\overline g/r_0}{\delta+M_0r_0}.
\end{equation*}
Since $\log(1+x)=\int_0^x\frac{dt}{1+t}\leq\frac{x^p}{p}$ for $x\geq0$ and $p\in(0,1]$,
we have $\log_2(1+x)\leq6x^p$ for $p\in\lb\frac12,\frac27\rb$. Hence,
\begin{equation*}
  \sum_{i=1}^{m}b_i\leq5m+6\sqrt\delta\sum_{i=1}^{m}\eta_i^{-1/2}
  +6M_0^{2/7}\sum_{i=1}^{m}\ls\frac{\norm{s_i}}{\eta_i}\rs^{2/7}
  +|\mathcal I_m|\log_{2,+}\frac{\delta+\overline g/r_0}{\delta+M_0r_0}.
\end{equation*}
Substituting~\eqref{eq:ahpe_damped_weighted_sums} and~\eqref{eq:ahpe_damped_long_count}
proves~\eqref{eq:ahpe_damped_rejection_sum}.
\end{proof}

\begin{proof}[Proof of Corollary~\ref{cor:ahpe_damped_oracle_complexity}]
By Lemma~\ref{lem:ahpe_damped_initialization}, initialization rejects at most
$1+\log_{2,+}\frac{\overline\tau_0}{\eta_0}$ trials.
If $A_1\geq R^2/(2\e)$, then $N=1$ and there are no main iterations.
Otherwise, the minimality of $N$ and~\eqref{eq:ahpe_damped_weight_overshoot} give
$A_N<4A_{N-1}<2R^2/\e$. Together with~\eqref{eq:ahpe_damped_initial_weight}, this gives
\begin{gather*}
  \sqrt{\delta A_N}\leq\sqrt{\frac{2\delta R^2}{\e}},\quad
  (M_0RA_N)^{2/7}\leq\ls\frac{2M_0R^3}{\e}\rs^{2/7}, \quad
  \frac{A_N}{A_1}\leq\frac{4\max\lb\eta_0,\overline\tau_0\rb R^2}{\e}.
\end{gather*}
Substituting these bounds into~\eqref{eq:ahpe_damped_rejection_sum} with $m=N-1$
and adding the initialization rejections, we obtain that the total number of rejected
trials is at most the right-hand side of~\eqref{eq:ahpe_damped_gradient_count} with $7N$ replaced by $5N$.
Each outer iteration makes one call $\MSBacktrack(x,\eta)$, which uses $H(x)$, $\grad f(x)$,
and one gradient for each trial with a positive definite matrix.
The updates reuse the gradient at the accepted point.
Hence, the method uses at most $N$ Hessian evaluations, and the number of gradient
evaluations is at most $2N$ plus the number of rejected trials.
This proves~\eqref{eq:ahpe_damped_gradient_count}.

\end{proof}

\section{Extra Experiments}
\label{app:extra_experiments}

\subsection{Experiment Details}
\label{app:experiment_details}
Table~\ref{tab:exp_constants} gives the initial parameters for experiments on Figure~\ref{fig:exp}. ACRN and A-ACRN use $\delta_0=0$ and adapt $\M$ by doubling at each rejected trial~\citep{kamzolov2026accelerated}.

\begin{table}[h]
  \centering
  \caption{Initial parameters of the runs in Figure~\ref{fig:exp}, from $x_0=-\mathbf 1$.}
  \label{tab:exp_constants}
  \setlength{\tabcolsep}{5pt}
  \begin{tabular}{lcccc}
    \toprule
    & \texttt{gisette} & \texttt{german\_numer} & \texttt{ijcnn1} & \texttt{abalone} \\
    \midrule
    $L(x_0)$ & $0.1$ & $0.1$ & $0.1$ & $789.5$ \\
    $\norm{\nabla f(x_0)}$ & $0.4233$ & $0.3268$ & $0.07210$ & $19.48$ \\
    $\norm{x_0}$ & $70.71$ & $4.899$ & $4.690$ & $3$ \\
    $\eta_0$ (Algorithms~\ref{alg:basic_pd_ms}, \ref{alg:ahpe_damped}, \ref{alg:adaptive_ahpe}) & $0.2057$ & $0.1808$ & $0.08491$ & $124.0$ \\
    $H_0$ (AdaN) & $0.05$ & $0.05$ & $0.05$ & $394.8$ \\
    $\gamma$ (Extra-Newton) & $141.4$ & $9.798$ & $9.381$ & $6$ \\
    $\beta_0$ (Extra-Newton) & $846.6$ & $3.138$ & $0.6345$ & $5.537\cdot10^5$ \\
    $\M_0$ (ACRN, A-ACRN, Algorithm~\ref{alg:fully_adaptive_cubic}) & $0.1$ & $0.1$ & $0.1$ & $789.5$ \\
    \bottomrule
  \end{tabular}
\end{table}

\subsection{Sensitivity to Initial Parameters}
\label{app:sensitivity}
\begin{figure}[ht]
  \centering
  \includegraphics[width=\textwidth]{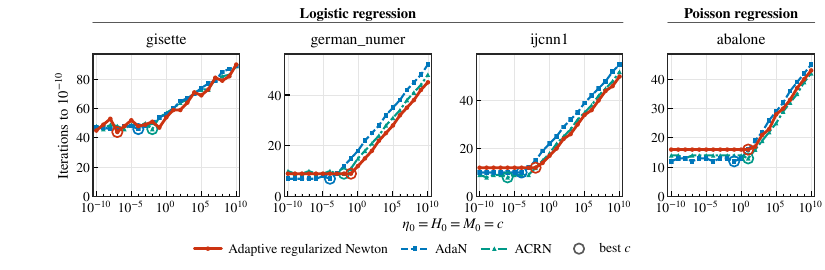}
  \caption{Iterations to reach $f(x_k)-f^\ast\leq10^{-10}$ for $\eta_0=H_0=M_0=c$. Circles mark the best $c$ of each method.}
  \label{fig:sensitivity_iterations}
\end{figure}

\begin{figure}[ht]
  \centering
  \includegraphics[width=\textwidth]{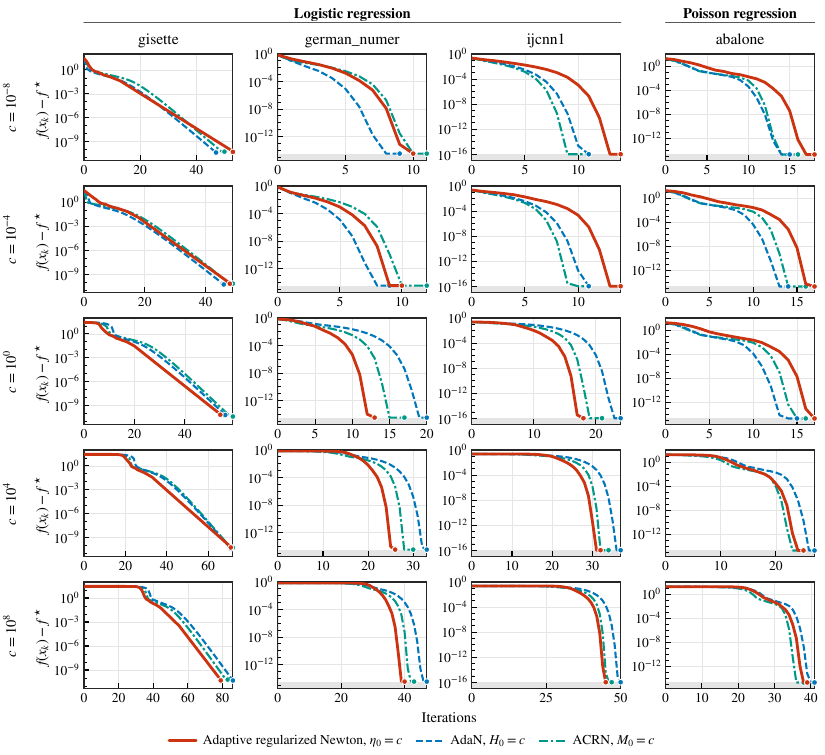}
  \caption{Convergence of Algorithm~\ref{alg:basic_pd_ms}, AdaN and ACRN for $\eta_0=H_0=M_0=c$, one row per value of $c$.}
  \label{fig:sensitivity_curves}
\end{figure}

We start from $x_0=-\mathbf1$, as in Figure~\ref{fig:exp}, and run Algorithm~\ref{alg:basic_pd_ms}, AdaN and ACRN with exact Hessians. Their initial parameters $\eta_0$, $H_0$ and $M_0$ take the same value $c$ from the grid $\lb10^{-10},10^{-9},\ldots,10^{10}\rb$. Figure~\ref{fig:sensitivity_iterations} shows the number of iterations to reach $f(x_k)-f^\ast\leq10^{-10}$ for every $c$, and Figure~\ref{fig:sensitivity_curves} shows the convergence for five values of $c$.

\end{document}